\documentclass[10pt,twoside]{article}
\usepackage[all]{xy}
\usepackage{stmaryrd}
\usepackage[usenames, dvipsnames]{color}
\usepackage[pdfborder={0 0 0}, colorlinks = true, linkcolor=blue, citecolor=OliveGreen]{hyperref}
\usepackage{graphicx}
\usepackage{amssymb}
\usepackage{amsthm}
\usepackage{amsmath, amscd}
\usepackage{mathrsfs} 
\usepackage{paralist}

\DeclareMathOperator{\Q}{\mathbb{Q}}
\DeclareMathOperator{\Z}{\mathbb{Z}}
\DeclareMathOperator{\C}{\mathbb{C}}
\DeclareMathOperator{\A}{\mathbb{A}}
\DeclareMathOperator{\R}{\mathbb{R}}
\DeclareMathOperator{\F}{\mathbb{F}}

\DeclareMathOperator{\D}{\mathcal{D}}

\newcommand{\la}{\langle} 	
\newcommand{\ra}{\rangle} 	
\newcommand{ \lie}[1]{\mathfrak{#1} } 
\newcommand{\isomto}{\overset{\sim}{\longrightarrow}}

\theoremstyle{plain}
\newtheorem{theorem}{Theorem}[section]
\newtheorem*{mainthm*}{Main Theorem}
\newtheorem*{theorem*}{Theorem}
\newtheorem{lemma}[theorem]{Lemma}
\newtheorem{proposition}[theorem]{Proposition}
\newtheorem{definition}[theorem]{Definition}
\newtheorem{corollary}[theorem]{Corollary}

\numberwithin{equation}{section}

\theoremstyle{remark}
\newtheorem{remark}[theorem]{Remark}

\DeclareMathOperator{\OB}{\mathcal{O}_B}
\DeclareMathOperator{\OBx}{\mathcal{O}_B^{\times}}
\DeclareMathOperator{\OBp}{\mathcal{O}_{B,p}}

\DeclareMathOperator{\Zed}{\mathcal{Z}}

\DeclareMathOperator{\CB}{\mathcal{C}_{B}}
\newcommand{\CBok}{ \ensuremath{ \mathcal{C}_{B   {\scriptstyle  / o_k }}  }}
\DeclareMathOperator{\X}{\mathbb{X}}
\DeclareMathOperator{\Y}{\mathbb{Y}}
\DeclareMathOperator{\V}{\mathbb{V}}

\DeclareMathOperator{\Fpbar}{ \mathbb{F}}

\DeclareMathOperator{\Spec}{Spec}
\DeclareMathOperator{\Spf}{Spf}
\DeclareMathOperator{\Hom}{Hom}
\DeclareMathOperator{\End}{End}

\def\YEAR{\year}\newcount\VOL\VOL=\YEAR\advance\VOL by-1995
\def\firstpage{1}\def\lastpage{1000}
\def\received{}\def\revised{}
\def\communicated{}

\makeatletter
\def\magnification{\afterassignment\m@g\count@}
\def\m@g{\mag=\count@\hsize6.5truein\vsize8.9truein\dimen\footins8truein}
\makeatother

\font\eightrm=cmr8
\font\caps=cmcsc10                    
\font\Caps=cmcsc10 scaled \magstep1   

\makeatletter
\@specialpagefalse\headheight=8.5pt
\def\DocMath{}
\renewcommand{\@evenhead}{%
    \ifnum\thepage>\lastpage\rlap{\thepage}\hfill%
    \else\rlap{\thepage}\slshape\leftmark\hfill{\caps\SAuthor}\hfill\fi}%
\renewcommand{\@oddhead}{%
    \ifnum\thepage=\firstpage{\DocMath\hfill\llap{\thepage}}%
    \else{\slshape\rightmark}\hfill{\caps\STitle}\hfill\llap{\thepage}\fi}%
\makeatother

\def\TSkip{\bigskip}
\newbox\TheTitle{\obeylines\gdef\GetTitle #1
\ShortTitle  #2
\SubTitle    #3
\Author      #4
\ShortAuthor #5
\EndTitle
{\setbox\TheTitle=\vbox{\baselineskip=20pt\let\par=\cr\obeylines%
\halign{\centerline{\Caps##}\cr\noalign{\medskip}\cr#1\cr}}%
	\copy\TheTitle\TSkip\TSkip%
\def\next{#2}\ifx\next\empty\gdef\STitle{#1}\else\gdef\STitle{#2}\fi%
\def\next{#3}\ifx\next\empty%
    \else\setbox\TheTitle=\vbox{\baselineskip=20pt\let\par=\cr\obeylines%
    \halign{\centerline{\caps##} #3\cr}}\copy\TheTitle\TSkip\TSkip\fi%
\centerline{\caps #4}\TSkip\TSkip%
\def\next{#5}\ifx\next\empty\gdef\SAuthor{#4}\else\gdef\SAuthor{#5}\fi%
\ifx\received\empty\relax
    \else\centerline{\eightrm Received: \received}\fi%
\ifx\revised\empty\TSkip%
    \else\centerline{\eightrm Revised: \revised}\TSkip\fi%
\ifx\communicated\empty\relax
    \else\centerline{\eightrm Communicated by \communicated}\fi\TSkip\TSkip%
\catcode'015=5}}\def\Title{\obeylines\GetTitle}
\def\Abstract{\begingroup\narrower
    \parskip=\medskipamount\parindent=0pt{\caps Abstract. }}
\def\EndAbstract{\par\endgroup\TSkip}

\long\def\MSC#1\EndMSC{\def\arg{#1}\ifx\arg\empty\relax\else
     {\par\narrower\noindent%
     2010 Mathematics Subject Classification: #1\par}\fi}

\long\def\KEY#1\EndKEY{\def\arg{#1}\ifx\arg\empty\relax\else
	{\par\narrower\noindent Keywords and Phrases: #1\par}\fi\TSkip}

\newbox\TheAdd\def\Addresses{\vfill\copy\TheAdd\vfill
    \ifodd\number\lastpage\vfill\eject\phantom{.}\vfill\eject\fi}
{\obeylines\gdef\GetAddress #1
\Address #2 
\Address #3
\Address #4
\EndAddress
{\def\xs{4.3truecm}\parindent=0pt
\setbox0=\vtop{{\obeylines\hsize=\xs#1\par}}\def\next{#2}
\ifx\next\empty 
     \setbox\TheAdd=\hbox to\hsize{\hfill\copy0\hfill}
\else\setbox1=\vtop{{\obeylines\hsize=\xs#2\par}}\def\next{#3}
\ifx\next\empty 
     \setbox\TheAdd=\hbox to\hsize{\hfill\copy0\hfill\copy1\hfill}
\else\setbox2=\vtop{{\obeylines\hsize=\xs#3\par}}\def\next{#4}
\ifx\next\empty\ 
     \setbox\TheAdd=\vtop{\hbox to\hsize{\hfill\copy0\hfill\copy1\hfill}
                \vskip20pt\hbox to\hsize{\hfill\copy2\hfill}}
\else\setbox3=\vtop{{\obeylines\hsize=\xs#4\par}}
     \setbox\TheAdd=\vtop{\hbox to\hsize{\hfill\copy0\hfill\copy1\hfill}
	        \vskip20pt\hbox to\hsize{\hfill\copy2\hfill\copy3\hfill}}
\fi\fi\fi\catcode'015=5}}\gdef\Address{\obeylines\GetAddress}

\begin{document}
\Title Unitary cycles on Shimura curves and the Shimura lift I
\ShortTitle Unitary cycles on Shimura curves I
\SubTitle   
\Author Siddarth Sankaran
\ShortAuthor 
\EndTitle
\Abstract This paper concerns two families of divisors, which we call the ``orthogonal'' and ``unitary'' special cycles, defined on integral models of  Shimura curves. The orthogonal family was studied extensively by Kudla-Rapoport-Yang, who showed that they are closely related to  the Fourier coefficients of modular forms of weight 3/2, while the `unitary' divisors are analogues of cycles appearing in more recent work of Kudla-Rapoport on unitary Shimura varieties. Our main result shows that these two families are related by (a formal version of) the Shimura lift.
\EndAbstract
\MSC 14G35, 11G18, 11F30.
\EndMSC
\KEY 
\EndKEY
\Address  Mathematisches Institut der Universit\"at Bonn \\ Endenicher Allee 60 \\ 53115  Bonn, Germany \\ sankaran@math.uni-bonn.de
\Address
\Address
\Address
\EndAddress

\tableofcontents
\begin{section}{Introduction}

In a series of works leading up to the monograph \cite{KRYbook}, Kudla, Rapoport and Yang study a family of arithmetic divisors that lie in the first arithmetic Chow groups of integral models of Shimura curves. One of their main results states that if one assembles these divisors into a formal generating series, the result is \emph{modular} of weight 3/2; in effect, pairing this generating series with a suitable linear functional yields the $q$-expansion of a weight 3/2 modular form. In this paper, together with its forthcoming sequel, we take up the study of the \emph{Shimura lift} of this generating series: our main result relates the Shimura lift to a generating series comprised of `unitary' divisors, that are analogues of the cycles constructed by Kudla-Rapoport in their more recent work \cite{KRunn} on unitary Shimura varieties.

In the present work, we focus on the geometry of the two families of divisors, and prove that their members satisfy relations that mirror the relations between the Fourier coefficients of a holomorphic modular form and those of its Shimura lift; as we discuss in more detail below, the key step is determining the local structure of the cycles in a formal neighbourhood of a prime of bad reduction.

Let $B$ be a rational indefinite quaternion algebra with discriminant $D_B$, and fix a maximal order $\OB$. Consider the moduli stack $\CB$ over $\Spec (\Z)$ that parametrizes pairs $\underline A = (A, \iota)$; here $A$ is an abelian scheme over some base scheme $S$, equipped with an $\OB$-action 
\[ \iota\colon \OB \to End(A) \] 
that satisfies a determinant condition, cf.\  Definition \ref{ShCurveDefn}. This stack is then an integral model for the classical Shimura curve associated to $B$. 

For a positive integer $n$, Kudla, Rapoport and Yang define (more or less, see Definition \ref{OrthCycleDef}) an \emph{`orthogonal' special cycle} $\Zed^o(n)$ as the moduli space of diagrams
\[ \xi \colon A \to A, \]
where $A$ is a $\OB$-abelian surface (i.e.\ a point of $\CB$), and $\xi$ is a traceless $\OB$-linear endomorphism of $A$ such that $\xi^2 = -n$. We may view $\Zed^o(n)$ as a cycle on $\CB$ via the natural forgetful map, and form the \emph{orthogonal generating series} (which, to emphasize the connection to modular forms, we have suggestively written as a $q$-expansion in terms of the variable $\tau \in \lie{H}$):
\[ \Phi^o(\tau) \ := \Zed^o(0) \ + \ \sum_{n>0} \Zed^o(n) \ q_{\tau}^n \ \in \ Div(\CB) \llbracket q_{\tau} \rrbracket, \]
here $Div(\CB)$ is the group of divisors on $\CB$, i.e., the free abelian group generated by the closed irreducible substacks of $\CB$ that are \'etale-locally Cartier, and $\Zed^o(0)$ is an appropriate constant term, cf.\ \eqref{orthGenSeriesDef}.

The `unitary' cycles are also defined by a moduli problem. Let $k = \Q(\sqrt{\Delta})$ be an imaginary quadratic field, where $\Delta<0$ is squarefree, and denote its ring of integers by $o_k$. Throughout this paper, we assume that (i) $\Delta$ is even and (ii) every prime $p |D_B$ is inert in $k$. The second condition implies the existence of an embedding $\phi \colon o_k \to \OB$, and we fix one such embedding for the moment.

Consider the moduli stack $\mathcal E^+$ over $\Spec(o_k)$ that parametrizes tuples 
\[ \underline E = (E, i_E, \lambda_E).\]
 Here $E$ is an elliptic curve over some base scheme $S$ over $\Spec(o_k)$, endowed with an $o_k$-action 
\[ i_E \colon o_k \to \End(E) \]
and a compatible principal polarization $\lambda_E$. We also assume that  on the Lie algebra $Lie(E)$, the action induced by $i_E$ agrees with the structural morphism $o_k \to \mathcal O_S$.  Let $\mathcal E^-$ denote the moduli space of tuples $\underline E = (E, i_E, \lambda_E)$ as before, except now we insist that on $Lie(E)$, the action of $i_E$ is given by the \emph{conjugate} of the structural morphism. Finally we take
\[ \mathcal E \ = \ \mathcal E^+ \ \coprod \mathcal E^- \] 
to be the disjoint union of these two stacks over $\Spec(o_k)$.

Fix a base scheme $S$ over $o_k$, and suppose we are given a pair of $S$-points $\underline E \in \mathcal E(S)$ and  $\underline A = (A, \iota) \in \CB(S)$. Following \cite{KRunn}, we form the space of \emph{special homomorphisms}
\begin{equation*} Hom_{o_k, \phi}( E,  A) := \{ y \in Hom_S(E, A) \ | \ y \circ i_E(a) = \iota \left( \phi(a) \right) \circ y \ \text{for all } a \in o_k  \},
\end{equation*}
which is equipped with an $o_k$-hermitian form $h^{\phi}_{E, A}$, cf.\ \eqref{hEADefEqn}.

For an integer $m>0$, the \emph{unitary special cycle} $\Zed(m, \phi)$ is the moduli stack over $\Spec(o_k)$ that parametrizes tuples $ ( \underline E, \underline A, y)$, where $\underline E$ and $\underline A$ are as above and 
\begin{equation} \label{uCycleDiagEqn}
 y \in \Hom_{o_k, \phi}(E, A) \qquad \text{such that }h^{\phi}_{E,A}(y, y) = m. 
\end{equation}
Again, we may view $\Zed(m, \phi)$ as a cycle on $\CBok := \CB \times \Spec(o_k)$ via the natural forgetful map, which in fact turns out to be a divisor.

We also define the following rescaled version:
\[ \Zed^*(m, \phi) \ := \ \Zed \left( \frac{m}{\gcd(m, D_B) } , \phi \right) \]
Finally, we form the \emph{unitary generating series}:
\begin{align*}
 \Phi^u(\tau)  :=  \Zed(0) +\frac{1}{2 h(k)} \sum_{m > 0}\left(  \sum_{[\phi] \in Opt / \OBx}  \Zed(m, \phi) \ + \ \Zed^*(m , \phi ) \right) q_{\tau}^m  \\
 \in Div(\CBok) \ \llbracket q_{\tau} \rrbracket \otimes_{\Z} \Q
 \end{align*}
for an appropriate constant term $\Zed(0)$, cf.\ \eqref{unGenSeriesDef}, and the sum on $[\phi]$ runs over an equivalence class of optimal embeddings (cf.\ the notation section below). 

Our main theorem describes the relationship between $\Phi^o$ and $\Phi^u$ in terms of the \emph{Shimura lift}, which is a classical operation on modular forms. It takes as its input a modular form $F$ of half-integral weight together with a squarefree integer parameter $t$, and yields a modular form $Sh(F)$ of \emph{even} integral weight. Moreover, when $F = \sum a(n) q^n$ is holomorphic, there are explicit formulas (depending on $t$ together with the weight, level, and character) for the Fourier coefficients of $Sh(F)$ in terms of those Fourier coefficients of $F$ that are of the form $a(tm^2)$. We define the \emph{formal Shimura lift} to be the operator on formal power series determined by these formulas, cf. \S \ref{ShimLiftSect}. 

\begin{mainthm*}[see Theorem \ref{mainThm}] \label{MainThmIntro} Suppose $k =  \Q(\sqrt{\Delta})$, where $\Delta < 0$ is a squarefree even integer and assume further that every prime dividing $D_B$ is inert in $k$.  Let $\Phi^o_{/o_k}$ denote the generating series
\[ \Phi^o_{/o_k} (\tau) \ = \ \Zed^o(0)_{/o_k} \ + \ \sum_{n>0} \Zed^o(n)_{/o_k} \ q_{\tau}^n \ \in \ Div(\CBok) \llbracket q_{\tau} \rrbracket \]
obtained by taking the base change to $o_k$ of each coefficient. Then we have an equality of formal generating series
\begin{equation} \label{MainThmIntroEqn}
Sh ( \Phi^o_{ /o_k}) (\tau) = \Phi^u(\tau), 
\end{equation}
in $Div(\CBok) \llbracket q_{\tau} \rrbracket \otimes \Q$. 
\end{mainthm*}
%

We now give an outline of the proof. It turns out that both the orthogonal and unitary divisors have vertical components only at primes $p|D_B$.   This leads us to study the \emph{p-adic uniformization} of the Shimura curve $\CB$, which relates the formal completion along its fibre at such a prime $p$ to (a formal model  of) the Drinfeld $p$-adic upper half-plane $\D$. We also have $p$-adic uniformizations for the orthogonal and unitary special cycles, which are expressed in terms of linear combinations of analogous `local' cycles defined on $\D$.

 Let $\F = \F_p^{alg}$ be an algebraic closure of $\F_p$, and let $W= W(\F)$ the ring of Witt vectors. Denote by $\mathbf{Nilp}$ the category of $W$-schemes such that $p$ is locally nilpotent, and for a scheme $S \in \mathbf{Nilp}$, let $\overline S := S \times_W \F$. Then the Drinfeld upper half-plane $\D$ parametrizes tuples $\underline X = (X, \iota_X, \rho_X)$, where $X$ is a $p$-divisible group of height 4 and dimension 2, over a base scheme $S \in \mathbf{Nilp}$, together with a ``special'' action $\iota_X: \OBp \to End(X)$  of the maximal order $\OBp$ (cf. \S\ref{DrinfeldSec}); finally $\rho_X$ is an $\OBp$-linear quasi-isogeny
\[ \rho_X: X \times_S \overline S \to \X \times_{\F} \overline S \]
of height 0, where $\X$ is some fixed $p$-divisible group over $\F$ endowed with a special $\OBp$-action. In this picture, the local analogues of the orthogonal and unitary cycles are described as the \emph{deformation loci} of homomorphisms of $p$-divisible groups. A complete description of the orthogonal cycles can be found in \cite{KRpadic}, and so our aim in Section \ref{locSect} is provide the same for the unitary cycles. 
 
Recent work of Kudla and Rapoport \cite{KRdrin} describes the special fibre of $\D$ in terms of the Bruhat-Tits tree for $SU(C)$, where $C$ is the split 2-dimensional hermitian space over $k_p$ (recall that by assumption, $k_p$ is an unramified quadratic extension of $\Q_p$). Combining this description with a healthy dose of Grothendieck-Messing theory, we are able to write down explicit equations for a local unitary cycle as a closed formal subscheme of $\D$, and we consequently obtain a precise description of its irreducible components. This description parallels the one found in \cite{KRpadic} for the orthogonal cycles, and by comparing the two formulas, we obtain the following key result (Theorems \ref{locCycCompThm} and \ref{locCycCompHorThm}):  any local orthogonal cycle that appears in the $p$-adic uniformization of an orthogonal cycle $\Zed^o(|\Delta|n^2)$ can be expressed as sum of two (explicitly determined) \emph{local} unitary cycles. 

We now turn to $p$-adic uniformizations. If we fix an $\F$-valued point $\underline{\mathbf A} = (\mathbf A, \iota_{\mathbf A}) \in \CB(\F)$, then the space of $\OB$-linear quasi-endomorphisms 
\[ B' := \End_{\OB}(\mathbf A)_{\Q} \]
is a definite quaternion algebra over $\Q$ with discriminant $D_B/p$. The $p$-adic uniformization of the Shimura curve $\CB$ can be expressed in the following way: there is a finite subgroup $\Gamma' \subset (B')^{\times}$ acting on $\D$, such that if we let $\widetilde \CB$ denote the base change to $W$ of the formal completion of $\CB$ along its fibre at $p$, then there is an isomorphism
\[ \widetilde \CB \ \simeq \ \left[ \Gamma' \big\backslash \D \right]\]
of formal stacks\footnote{Here and throughout this paper, we use the term `isomorphism of formal stacks' in a rather cavalier fashion. What we mean in all cases is the following: upon fixing sufficiently deep level structure away from $p$, one obtains a formal scheme on each side, together with a covering map from the corresponding stack. Our assertion is that there is an isomorphism of formal schemes compatible with the automorphisms of the covering maps, and the morphisms induced by varying the level structure, in the natural way.} over $W$. 

Similarly, let $\widetilde \Zed{}^o(n)$ and $\widetilde \Zed(m, \phi)$ denote the base change to $W$ of the formal completions of the special cycles $\Zed^o(n)$ and $\Zed(m, \phi)$ along their  fibres at $p$. Viewed as a cycle on $\widetilde \CB$, we may express an orthogonal cycle as a sum 
\[ \widetilde \Zed{}^o(n) = \sum_{\substack{\xi \in \Omega^o(n) \\ \text{mod } \Gamma'} } [Z^o(\xi[p^{\infty}])] \]
where $[Z^o(\xi[p^{\infty}])]$ is the projection to $[\Gamma' \backslash \D]$ of a \emph{local orthogonal cycle} on $\D$, and 
\[ \Omega^o(n) \subset \left\{ b' \in B' \ | \ Tr(b') = 0 \right\} \]
is a $\Gamma'$-invariant set of vectors of reduced norm $n$ satisfying a certain integrality property, cf.\ Theorem \ref{OrthUnifThm}.  

For the unitary cycles, we fix a triple $\underline{\mathbf E}\in \mathcal E(\F)$. A unitary cycle then decomposes as 
\[ \widetilde \Zed(m, \phi) = \frac{1}{|o_k^{\times}|} \sum_{[\lie{a}] \in Cl(k)} \left( \sum_{ \substack{\beta \in \Omega^+(m, \lie{a}, \phi) \\ \text{mod } \Gamma' }} [Z(\beta[p^{\infty}])]  +
\sum_{ \substack{\beta' \in \Omega^-(m, \lie{a}, \phi) \\ \text{mod } \Gamma' }} [Z(\beta'[p^{\infty}])] \right),  \]
where $Cl(k)$ is the class group of $k$, and as before $[Z(\beta[p^{\infty}])]$ and $[Z(\beta'[p^{\infty}])]$ are local unitary cycles. The sets $\Omega^{\pm}(m, \lie{a}, \phi)$ appearing above are subsets
\[ \Omega^{\pm}(m, \lie{a}, \phi) \subset \Hom(\mathbf E, \mathbf A) \otimes_{\Z} \Q, \]
consisting of quasi-morphisms of a specified norm, a linearity (or anti-linearity, in the case of $\Omega^-$) condition with respect to the action of $o_k$, and again satisfying an integrality condition, cf.\ Theorem \ref{ZuLocalUnifThm}.

Thus, in order to compare the unitary and orthogonal cycles,  we need to compare the indexing sets $\Omega^o(n)$ and $\Omega^{\pm}(m, \lie{a}, \phi)$, at least in the case that the squarefree part of $n$ is equal to $|\Delta|$, and as $\phi$ varies among classes of optimal embeddings. This task, which amounts to a study in the arithmetic of quaternion algebras, is carried out in \S \ref{CalcSec}. Together with the description of the local cycles discussed previously, we arrive at a relationship between the two families of cycles, in a formal neighbourhood at $p|D_B$, that matches exactly the formula for the Fourier coefficients of the Shimura lift. Since the vertical components of the cycles only occur at such primes, the main theorem follows immediately, cf.\ Theorem \ref{mainThm}.

I would like to conclude the introduction by placing this result in the context of the sequel \cite{San2} to the present work, whose aim is to establish the same Shimura lift formula in the first \emph{arithmetic} Chow group  of $\CB$ (in the sense of Gillet-Soul\'e). Recall that in \cite{KRYbook}, the authors prove that the generating series
\[ \widehat{\Phi}{}^o(\tau) := \sum_{n \in \Z} \ \widehat{\Zed}{}^o(n, v) \ q_{\tau}^n , \qquad \Im(\tau) = v \]
is a non-holomorphic modular form of weight 3/2; here the coefficients are arithmetic classes
\[ \widehat{\Zed}{}^o(n, v)  \ = \ \left(\Zed^o(n), Gr^o(n, v) \right) \in \widehat{CH}{}^1(\CB) \] 
 for appropriate choices of Green functions $Gr^o(n,v)$, and constant term $\widehat{\Zed}{}^o(0,v)$. 
The modularity of $\widehat{\Phi}{}^o(\tau)$ means that in particular, applying a linear functional to its coefficients yields the $q$-expansion of a modular form in the usual sense (in fact, their proof of modularity amounts to showing that this is true for a well-chosen set of linear functionals). 

In \cite{San2}, we augment the unitary cycles with  Green functions $Gr(m,v,\phi)$ and obtain classes
\[ \widehat{\Zed}(m, v, \phi) \ = \ (\Zed(m,\phi), Gr(m,v,\phi) ) \in \widehat{CH}{}^1(\CB) \]
and, together with an appropriate constant term $\widehat \Zed(0,v)$, define
\[ \widehat{\Phi}{}^u(\tau) := \widehat\Zed(0,v) + \frac{1}{2h(k)} \ \sum_{m \in \Z} \sum_{[\phi]} \left[ \widehat{\Zed}(m,v,\phi) + \widehat{\Zed}{}^*\left(m, v,\phi \right) \right] q^m. \]
The main result of \cite{San2} identifies this generating series with the Shimura lift of $\widehat\Phi{}^o$, by combining the geometric results of the present work with calculations that relate the Green functions. In particular, applying a linear functional to both series yields an identity involving the \emph{classical Shimura lift} of modular forms.

Some linear functionals are `geometric' in nature, i.e., they only depend on the cycles and not the Green functions. When this is the case, the main result of the present work already implies that the Shimura lift relation holds for the corresponding modular forms, and we give a few examples at the end of this paper that will play a role in the sequel.

\vspace{10pt}
\noindent \textbf{Acknowledgments:} This paper is an extension of the work carried out in my Ph.\ D.\ thesis; I am grateful to my advisor S. Kudla for introducing me to the subject, and his consistent encouragement. I would also like to thank B. Howard, M. Rapoport, and the anonymous referee for helpful comments. This work was supported by the SFB/TR 45 `Periods, Moduli Spaces, and Arithmetic of Algebraic Varieties' of the DFG (German Research Foundation).

\subsection*{Notation:} 
\begin{compactitem}
\item B is an indefinite quaternion algebra over $\Q$, with discriminant $D_B$, and with reduced trace and norm denoted by $Trd$ and $Nrd$ respectively. We fix a maximal order $\OB \subset B$. 
\item $k = \Q(\sqrt{\Delta})$ is an imaginary quadratic field, with ring of integers $o_k$. We denote the non-trivial Galois automorphism by $a \mapsto a'$. We also assume throughout this paper that 
	\begin{compactenum}[(i)]
		\item $\Delta<0$ is a squarefree \emph{even} integer;
		\item and every prime $p|D_B$ is inert in $k$.
	\end{compactenum}
\item By definition, an embedding $\phi \colon o_k \hookrightarrow \OB$ is \emph{optimal} if 
\begin{equation} \label{optEmbDefEqn} \phi(o_k) \ = \ \phi(k) \cap \OB . \end{equation}
Let $Opt$ denote the set of optimal embeddings. Note that $\OBx$ acts on $Opt$ by conjugation: for $\xi \in \OBx$, set
\[  (\xi \cdot \phi)(a) \ := \ (Ad_{\xi} \circ \phi)(a) \ = \ \xi \cdot \phi(a) \cdot \xi^{-1}, \qquad a \in o_k. \]
Denote by $Opt / \OBx$ the set of equivalence classes of optimal embeddings under this action.  
\item $\chi_k$ is the quadratic character associated to $k$, so that for a prime $p$, 
\begin{equation} \label{chiKNotSecEqn}
\chi_k(p) = \begin{cases} 1, & \text{if } p \text{ splits,} \\ 0, & \text{if } p \text{ is ramified}, \\ -1, & \text{if } p \text{ is inert}. \end{cases} 
\end{equation}
\end{compactitem}
\end{section}


\begin{section}{Shimura curves and global special cycles} \label{GlobalDefsSec}
In this section, we recall the construction of Shimura curves, and the global orthogonal and unitary special cycles on them. 

\begin{definition}[Shimura curve] \label{ShCurveDefn} Let $\CB$ denote the moduli problem which associates to a scheme $S$ over $\Spec (\Z)$ the category whose objects are pairs 
\[ \CB(S) = \{ (A, \iota) \},\]
 where (i) $A$ is an abelian surface over $S$, and (ii) $\iota\colon \OB \to \End_S(A)$ is an action of $\OB$ on $A$. We also require that for every $b \in \OB$, 
\[ det(T - \iota(b)|_{Lie(A)}) = T^2 - Trd(b)T + Nrd(b) \in \mathcal O_S[T]. \]
\end{definition}

\begin{proposition}[cf.\ \cite{KRYbook}, Proposition 3.1.1] The moduli problem $\CB$ is representable by a Deligne-Mumford (DM) stack, which we also denote by $\CB$. It is regular, proper and flat over $Spec(\Z)$ of relative dimension 1, and smooth over $Spec \Z[ D_B^{-1}] $.
 \qed
\end{proposition}

The orthogonal special cycles, as constructed in e.g.\ \cite{KRYbook}, are also defined by a moduli problem:
\begin{definition} \label{OrthCycleDef}
For $n \in \Z_{>0}$, let $\Zed^o(n)^{\sharp}$ denote the DM stack which represents the following moduli problem over $Spec \Z$: to a scheme $S/ \Z$, we let $\Zed^o(n)^{\sharp}(S)$ be the category of tuples $(A, \iota, \xi)$, where
\begin{compactenum}
\item $(A, \iota) \in \CB(S)$, and 
\item $\xi \in End_{\OB}(A)$ is an $\OB$-linear endomorphism of $A$ with $Tr(\xi) = 0$ and $\xi^2 = -n$.
\end{compactenum}
We view $\Zed^o(n)^{\sharp}$ as a cycle on $\CB$ via the natural forgetful map $\Zed^o(n)^{\sharp} \to \CB$, which is finite and unramified by \cite[\S3.4]{KRYbook}. Let $\Zed^o(n)$ denote the Cohen-Macauleyfication of $\Zed^o(n)^{\sharp}$, as in p.\ 55 of loc.\ cit.,  so that $\Zed^o(n)$ is of pure codimension 1 in $\CB$. 
\end{definition}

We now turn to the unitary special cycles, following \cite{KRunn}. Fix once and for all an element $\theta \in \OB$ such that $\theta^2 = - D_B$. 
Given a point $(A, \iota) \in \CB(S)$,  there exists a unique principal polarization $\lambda^0_A$ on $A$ such that the Rosati involution $\varphi \mapsto \varphi^{\dagger}$ satisfies 
\[ \iota(b)^{\dagger} = (\lambda^0_A)^{-1} \circ \iota(b)^{\vee} \circ \lambda^0_A = \iota(\theta^{-1}  \ b^{\iota} \ \theta  ) , \qquad \text{for all } b \in \OB,\]
cf.\ \cite[\S 3.1]{How}. 

\begin{lemma} Let $\phi \colon o_k \to \OB$ be an embedding. Then
	\begin{compactenum}[(i)]
		\item $Trd(\theta \phi( \sqrt{\Delta})) \neq 0$;
		\item Assume, without loss of generality, that $Trd(\theta \phi( \sqrt{\Delta}) )>0$, by replacing $\theta$ by $-\theta$ if necessary. Then the isogeny 
\begin{equation} \label{lambdaADefEqn}
  \lambda_{A, \phi} \  :=  \ \lambda_A^0 \ \circ \ \iota \Big( \theta  \ \phi(\sqrt{\Delta}) \Big).
\end{equation}		is a (non-principal) polarization.
	\end{compactenum}
\begin{proof}
(i) Since $\theta^2$ and $\Delta$ are both negative, the statement follows immediately from the assumption that $B$ is indefinite. 

(ii) We need to check that on geometric points, the map $\lambda_{A, \phi}$ is induced by an ample line bundle; therefore we assume $A$ is an abelian variety over an algebraically closed field. The endomorphism $\iota( \theta \phi(\sqrt{\Delta}))$ is symmetric with respect to $\dagger$, and so lies in the image of  the map
\[ NS(A) \to \End(A) \]
where $NS(A)$ is the N\'eron-Severi group. An element in this image comes from an \emph{ample} bundle if and only if it is totally positive, i.e.\ its characteristic polynomial has positive real roots, cf.\ \cite[\S 21]{Mum}. The roots of the characteristic polynomial are the same as those of its minimal polynomial, cf. \emph{loc.\ cit.} p.~203, and  $\iota(\theta \phi(\sqrt{\Delta}))$ has minimal polynomial
\[ P(x) =  x^2 \ - \  Trd(\theta \phi(\sqrt{\Delta}))\cdot x \  + \ Nrd(\theta \phi(\sqrt{\Delta})). \]
We may fix an isomorphism $B \otimes \R \simeq M_2(\R)$ such that $Trd$ and $Nrd$ are identified with the trace and determinant respectively, and the positive involution 
\[ b \ \mapsto \ b^{\dagger} \ = \  \theta^{-1} \cdot  b^{\iota} \cdot \theta \]
is identified with the transpose operator. Then $\theta \phi(\sqrt{\Delta})$ is identified with a symmetric matrix $A = A^{t}$; the conditions $tr(A)>0$ and $det(A) = |\Delta|D_B > 0 $ then imply that the roots of $P(x)$ are real and positive. 
\end{proof}
\end{lemma}Note that the Rosati involution $*$ associated to the polarization $\lambda_{A, \phi}$ satisfies
\[ \iota_A \left( \phi(a) \right)^{*} = \iota_A \left( \phi(a') \right), \qquad \text{ for all } a \in o_k. \]

Let $\mathscr E^+$ denote the moduli stack over $\Spec(o_k)$ such that for a base scheme $S/o_k$ , the $S$-points parametrize tuples
\[ \mathscr E^+(S) \ = \ \left\{ \underline E =  (E, i_E, \lambda_E \right) \}; \]
here $E$ is an elliptic curve over $S$ endowed with an action $i_E\colon o_k \to \End(E)$, and a principal polarization $\lambda_E$ such that the induced Rosati involution $^*$ satisfies
\[ i_E(a)^* = i_E(a'). \] 
We further impose the condition that the action of $o_k$ on $Lie(E)$ induced by $i_E$ agrees with the action given via the structural morphism $o_k \to \mathcal O_S$. 

Similarly, we define $\mathscr E^-$ to be the moduli space of tuples $\underline E = (E, i_E, \lambda_E)$ as above, except we require that the action on $Lie(E)$ induced by $i_E$ is equal to the \emph{conjugate} of the structural morphism.
Finally, we take
\[ \mathscr E \ = \ \mathscr E^+ \coprod \mathscr E^- \] 
to be the disjoint union of these stacks.

Suppose $S$ is a scheme over $o_k$, and we are given two points $\underline A \in \CB(S)$ and $\underline E \in \mathscr E(S)$. We form the space of \emph{special homomorphisms}:
\[ \Hom_{\phi}(\underline E, \underline A)  :=  \left\{ {y} \in \Hom(E, A) \ | \  y \circ i_E(a) \ = \ \iota_A( \phi (a) ) \circ  y, \ \text{for all } a \in o_k \right\}\]
This space comes equipped with an $o_k$-hermitian form $h_{E,A}^{\phi}$ defined by 
\begin{equation} \label{hEADefEqn}
 h_{E,A}^{\phi} ( s,  t) \ := \ (\lambda_E)^{-1} \circ { t}^{\vee} \circ \lambda_{A, \phi} \circ {s} \ \in \End(E, i_E) \simeq o_k. 
\end{equation}

\begin{definition}[Unitary special cycles.] Suppose $m \in \Z_{>0}$ and $\phi\colon o_k \to \OB$ is an optimal embedding. Let $\Zed(m, \phi)$ denote the DM stack over $Spec (o_k)$ representing the following moduli problem: for a scheme $S / o_k$, we define $\Zed(m, \phi)(S)$ to be the category of tuples
\[ \Zed(m, \phi)(S) = \Big\{ (\underline E, \underline A, y)  \Big\} \]
where (i) $\underline E \in \mathscr{E}(S)$, 
(ii) $\underline A \in \CB(S)$, and
(iii) $y \in Hom_{o_k, \phi}(\underline E, \underline A)$ such that $h_{E,A}^{\phi}(y, y) = m$. 
\end{definition}

By the proof of \cite[Proposition 2.10]{KRunn}, which applies verbatim to the present setting, the forgetful map  
\[ \Zed(m, \phi) \to \CBok := \CB \times_{\Z} \Spec(o_k) \]
is finite and unramified, and so we may view $\Zed(m, \phi)$ as a cycle on $\CBok$; abusing notation, we shall refer to the cycle by the same symbol, and hope that context will suffice to clarify which instance of the notation is intended. 

\begin{proposition}  \label{UCycleIsDivisorProp}
The cycle $\Zed(m,\phi)$ is a divisor on  $\CBok$, i.e., each irreducible component is of codimension 1. 
  Moreover, it has vertical components in characteristic $p$ if and only if (i)  $p$ divides $D_B$ and (ii) $ord_pm > 0$. 

\begin{proof} 
By the complex uniformization \cite[Theorem 3.3.3]{San}, the irreducible components of the generic fibre  $\Zed(m, \phi)_k$ are $0$-dimensional, and so their closures in $\CBok$ (i.e.\  the horizontal components) are 1-dimensional.

Let $p|D_B$; in Section \ref{locEqnsSec} below, we show that  $\Zed(m,\phi)$ is a divisor in a formal neighbourhood of the fibre at  $p$ by writing down explicit equations. In particular, Theorem \ref{cycleDecompThm} asserts that $\Zed(m,\phi)$ contains vertical irreducible components in this fibre if and only if $ord_p(m)>0$. 

It remains to consider the fibre at primes $\lie{p}\subset o_k$ with $(\lie{p}, D_B) = 1$. Suppose $z \in \Zed(m, \phi)(\F)$ is a geometric point, for an algebraically closed field $\F$ over $o_k / \lie{p}$, corresponding to a triple $(\underline E, \underline A, y)$ over $\F$. Let $x \in \CBok(\F)$ denote the point below $z$. 
We need to prove that (i) the cycle $\Zed(m, \phi)$ is a divisor at $x$ and (ii) there are no vertical components, i.e., that the data $(\underline E, \underline A, y)$ can be lifted to characteristic zero and $\Zed(m, \phi)$ does not contain the entire fibre at $\lie{p}$.  

Let $\mathcal{W}$ be the completion of the maximal unramified extension of the local ring $o_{k,\lie{p}}$, and  denote by $\mathbf{ART}$  the category of local Artinian $\mathcal W$-algebras with residue field $\F$. The natural morphism
\[ o_{k,\lie{p}} \to \mathcal W\]
endows any $\mathcal W$-algebra with an $o_{k, \lie{p}}$ structure, and so in particular every object of $\mathbf{ART}$ can be viewed as an $o_k$-algebra in this way. 

Let $R_z$ and $R_x$ denote the \'etale local rings of $z$ and $x$ respectively. Then $R_x$ pro-represents the functor of deformations of the data $\underline A = (A, \iota_A)$ to objects of $\mathbf {ART}$. By the Serre-Tate theorem, giving a deformation of $\underline A$ is equivalent to giving a deformation of the underlying $p$-divisible group (with induced $\OB \otimes \Z_p$-action):
\[\underline X =  (X, \ \iota_X) \ : = \ (A[p^{\infty}] , \ \iota_A \otimes \Z_p). \]
Similarly, the local ring $R_z$ pro-represents deformations of the data $(\underline E, \underline A, y)$, which in turn are equivalent to deformations of 
\[ (\underline Y,  \ \underline X, \ b) \ := \ (\underline E[p^{\infty}], \ \underline A[p^{\infty}], \ y[p^{\infty}] ),  \]
where $\underline Y = (Y, i_Y, \lambda_Y)$ is the $p$-divisible group attached to $E$, together with its induced $o_k \otimes \Z_p$ action, and principal polarization. 
More precisely, the $\OB$-action on $A$ and the $o_k$-action on $E$ induce actions
\[ \iota_{X} \colon \OB \otimes {\Z_p} \to \End_{\Z_p}(X) , \qquad \text{and} \qquad i_{Y}\colon o_k \otimes \Z_p \to End_{\Z_p}(Y) \]
respectively, such that 
\begin{equation} \label{bActionEqn}
 b \circ i_Y(a) = \iota_X (\phi(a))  \circ  b, \qquad \text{for all } a  \in o_{k,p}  :=  o_k \otimes \Z_p. 
\end{equation}
Recalling that $p \nmid D_B$, we may fix an isomorphism $\OB \otimes \Z_p \simeq M_2(\Z_p)$ such that $\phi(\sqrt{\Delta})$ is identified with the matrix $(\begin{smallmatrix} & 1 \\ \Delta & \end{smallmatrix})$; note that when $p=2$, this is possible on account of the assumption that $|\Delta|$ is even. The two idempotents 
\[ e_1 = \begin{pmatrix} 1&0\\0&0 \end{pmatrix} \qquad \text{and} \qquad  e_2 = \begin{pmatrix} 0&0\\0&1 \end{pmatrix} \]
determine a splitting $X  =  X^o \times X^o$, where $X^o$ is a $p$-divisible group of  dimension 1 and height 2; keep in mind that $X^o$ itself does not inherit any additional endomorphism structure. In a similar manner, any deformation of $\underline X = (X, \iota_X)$ decomposes into a product of two copies of a deformation of  $X^o$. Hence, ${Def}(X, \iota_X)  =  {Def}(X^o)$, and it is a well-known fact, cf.\ e.g.\ \cite[Theorem 3.8]{MierendorffZiegler}, that the latter deformation problem is pro-represented by $\Spf (\mathcal{W}  \llbracket t \rrbracket )$. 

Our aim is to analyze the deformation locus $ Def(\underline Y, \underline X, b)$ as a closed formal subscheme of $Def(\underline X) = Def( X^o) = \Spf (\mathcal{W}  \llbracket t \rrbracket )$. Under the splitting $X = X^o \times X^o$, we have a decomposition  $b = (b_1, b_2)$, where $b_i \in Hom(Y, X^o)$ are non-zero morphisms. Morever, because of \eqref{bActionEqn},  
\[ b_2 = b_1 \circ i_{Y}(\sqrt{\Delta}). \]
Thus, the problem of deforming $(\underline Y, \underline X, b)$ is equivalent to deforming $(\underline Y, X^o, b_1)$. 

The key tool is the theory of \emph{canonical lifts}: given a $p$-divisible group $\lie{g}$ over $\F$ of height 2 and dimension 1 together with an action 
\[ i_{\lie{g}} \colon o_{k, p} \to End(\lie{g}), \]
there exists a \emph{canonical lift} $\underline{\lie{G}} = (\lie{G}, i_{\lie{G}})$ over $\mathcal W$. This fact is a consequence of the existence of Serre-Tate canonical coordinates in the case when $\lie{g}$ is ordinary, cf.\ \cite[Appendix]{Messing}, and of Gross' theory of (quasi-)canonical liftings when $\lie{g}$ is supersingular, cf.\ \cite{Gross}. Moreover, for any object $R \in \mathbf{ART}$, there exists a unique lift of $\underline{\lie{g}}$ to $R$: namely, the base change $\underline{\lie{G}}_{/R}$ of $\underline{\lie{G}}$ to $R$. 
%
%

In particular, if we let $\underline{\lie{Y}}$  be the canonical lift of $\underline Y$ as above, then the problem of deforming $(\underline Y, X^o, b_1)$ to $R \in \mathbf{ART}$ is equivalent to finding deformations $\lie{X}^o$ of $X^o$ such that $b_1$ lifts to a morphism $\lie{Y}_{/R} \to \lie{X}^o$. 
It is not hard to show that the locus in $Def(X^o)$ where $b_1$ lifts is cut out by a single equation, cf.\ the proof of \cite[Prop. 5.1]{Wew}. Note that in our case, this equation cannot be 0. If it were, then $R_x = R_z$, and so $\Zed(m, \phi)$ would contain a positive-dimensional component in the generic fibre; this would contradict the description of the complex points above.  

Moreover, we claim $\Zed(m, \phi)$ cannot contain the entire special fibre $(\CB)_{\lie{p}}$. Note that if $X^o$ and $Y$ are isogenous, they are either both super-singular or both ordinary. But $Y$ is necessarily ordinary when $\lie{p}$ is a split prime, and supersingular otherwise, as it admits an action of $o_{k,p}$;  thus taking $X^o$ to be supersingular in the first case, or ordinary in the second, yields a point in the special fibre of $(\CB)_{\lie{p}}$ that cannot lie on $\Zed(m, \phi)$. 

It remains to show that there always exists a lift to characteristic $0$.  Suppose first that $X^o$ and $Y$ are ordinary. Then $p$ splits in $k$ and $End(X^o) = End(Y) = o_{k,p}$. Fixing an isomorphism $X^o \simeq Y$, we have 
\[ b_1 \in Hom(Y, X^o) \simeq End(X^o) \simeq o_{k,p}, \]
which lifts to an endomorphism of the canonical lift $\lie{X}^o$ of $X^o$. 

Next, we consider the case where $X^o$ and $Y$ are supersingular; this necessarily implies that $p$ is non-split in $k$. Let $D$ be the division quaternion algebra over $\Q_p$ with maximal order denoted by $ O_D$, and fix a uniformizer $\Pi \in  O_D$. 

Suppose first that $p = \lie{p}^2$ is ramified in $k$, and let $\varpi \in o_{k,p}$ be a uniformizer. Then we may identify 
\[ End(Y) \simeq  O_D \] 
such that $i_Y(\varpi)$ is identified with $\Pi$. We may also fix an isomorphism 
\begin{equation} \label{XoYisom} X^o \ \simeq  \ Y , \end{equation}
which induces isomorphisms
\[ End(X^o) \ \simeq Hom(Y, X^o) \simeq End(Y) \simeq  O_D \] 
in the natural way. Without loss of generality, we may normalize \eqref{XoYisom} such that $b_1$ is identified with $\Pi^t \in  O_D$, for some $t$. 
We let $i_{X^o} \colon o_{k, p} \to End(X^o)$ denote the action given by composing $i_Y$ with the isomorphism \eqref{XoYisom}. Then, appealing again to Gross' results, this action determines a canonical lifting  $\underline{\mathcal X}^o  = (\mathcal X^o, i_{\mathcal X^o})$ of $(X, i_X)$ to $\mathcal W$, and it is clear by construction that $b_1$ also lifts. 

Finally, we suppose that $p = \lie{p}$ is inert in $k$. We may fix an identification $\End(Y) \simeq  O_D$ such that 
\[ \Pi \cdot i_Y(a) = i_Y(a') \cdot \Pi , \qquad \text{for all } a \in o_{k,p} .\]
As before, we may also fix an isomorphism $\alpha\colon X^o \isomto Y$ such that $b_1$ is identified with $\Pi^t$ for some $t$. We then define an action $i_{X^o}: o_{k,p} \to End(X^o)$ by the formula:
\[ i_{X^o}(a) \ := \begin{cases} \alpha^{-1} \circ i_Y(a) \circ \alpha, \qquad \text{if } t \text{ is even,} \\   \alpha^{-1} \circ i_Y(a') \circ \alpha, \qquad \text{if } t \text{ is odd.} 
\end{cases}
\]
Then, by construction, the morphism $b_1$ will lift to a morphism $\lie{Y} \to \lie{X}^o$, where $\lie{Y}$ and $\lie{X}^o$ are the canonical lifts of $Y$ and $X^o$ to $\mathcal W$ determined by their respective $o_{k,p}$-actions as above. 

\end{proof}
\end{proposition}

\end{section}

\begin{section}{Local cycles on the Drinfeld upper half-plane} \label{locSect}
In this section, we discuss the local analogues of the unitary and orthogonal special cycles, and relate them to the global cycles via $p$-adic uniformization. 
Suppose $p |D_B$, so that in particular $p$ is inert in $k$ and $p \neq 2$. Let $k_p$ denote the completion at $p$, and let $\delta \in k_p$ denote the image of $\sqrt{\Delta}$. Throughout this section, we fix an embedding
\begin{equation} \label{fixedPhiLocEqn}
 \phi\colon o_{k,p}\hookrightarrow \OBp, 
\end{equation}
where $o_{k,p}$ and $\OBp$ are the maximal orders in $k_p$ and $B_p$ respectively. We also fix a uniformizer $\Pi \in \OBp$ such that $\Pi \phi(a) = \phi(a') \Pi$, for all $a \in o_{k,p}$. 

Let $\F = \F_p^{alg}$ denote a fixed algebraic closure of $\F_p$, and $W = W(\F)$ the ring of Witt vectors, with Frobenius endomorphism $\sigma\colon W \to W$. 
We fix an embedding $\tau_0 \colon o_{k,p} / (p) \to \F$, which lifts uniquely to an embedding $\tau_0 \colon o_{k,p} \to W$. In addition, we let $\tau_1$ denote the conjugate embedding (or its lift to characteristic 0). Note that these maps allow us to view $\F$ and $W$ as $o_k$-algebras.

Let $\mathbf{Nilp}$ denote the category of $W$-schemes for which the ideal sheaf generated by $p$ is locally nilpotent, and for a scheme $S \in \mathbf{Nilp}$, we set
 \[ \overline{S} := S \times_{W} \Spec(\F). \]
 
Additionally, we fix a trivialization of the prime-to-p roots of unity in $\F$: 
\[ \A_f^p(1): = \left( \prod_{\ell \neq p} \varprojlim_n \mu_{\ell^n}(\F) \right) \otimes_{\Z} \Q \simeq \A_f^p \]
Finally, we set $\widehat{\Z}{}^p := \prod_{\ell \neq p} \Z_{\ell}$, and if $M$ is any $\Z$-module, we define $\widehat{M}{}^p := M \otimes_{\Z} \widehat{\Z}{}^p$. 
 
\begin{subsection}{Drinfeld space} \label{DrinfeldSec}
In this section, we recall the definition of Drinfeld space as a moduli space of $p$-divisible groups, as well as an ``alternative'' description of its special fibre as a union of projective lines indexed by \emph{hermitian} lattices, as in \cite{KRdrin}.

Following \cite{BC}, we define a \emph{special formal} $\OBp$-\emph{module} over a scheme $S$ to be a pair $(X, \iota_X)$ consisting of a $p$-divisible group $X$ of height 4 and dimension 2 over $S$, and a map 
\[ \iota_X: \OBp \to End(X), \]
which satisfies the following \emph{special condition}: the Lie algebra $Lie(X)$ is (locally on $S$) a free $\mathcal{O}_S \otimes_{W, \tau_0} o_{k,p}$ module of rank 1. Here $Lie(X)$ is viewed as an $\mathcal{O}_S \otimes o_{k,p}$-module via the action of $o_{k,p}$ on $X$ induced by the composition $\iota_X \circ \phi$.

Fix once and for all a special formal $\OBp$-module $(\X, \iota_{\X})$ over $\F$; note such a pair is unique up to isogeny, by the classification of Dieudonn\'e isocrystals over $\F$ \cite[Proposition II.5.2]{BC}. The pair $(\X, \iota_{\X})$ serves as a ``base point'' for the following moduli problem:

\begin{definition} \label{drinfeldDef} Let $\D$  denote the moduli problem on $\mathbf{Nilp}$, which associates to a scheme $S$ the category of isomorphism classes of tuples
\[ \D(S) = \left\{ (X, \iota_{X}, \rho_X) \right\} / \simeq, \]
consisting of
\begin{compactitem}
\item  $(X, \iota_X)$ a special formal $\OB$-module over $S$, 
\item and an $\OBp$-linear quasi-isogeny
\[ \rho_X: X \times_S \overline{S} \ \to \ \X \times_{Spec(\Fpbar)} \overline S, \]
of height 0. 
\end{compactitem}
\end{definition}

An isomorphism between two tuples $(X, \iota, \rho)$ and $ (X', \iota', \rho')$ is an isomorphism $\alpha: X \to X'$ which is $\OBp$-equivariant, and such that $\rho = \rho' \circ (\alpha \times_S {S_0})$.

The functor $\D$ is then represented by  \emph{(a formal model of) the Drinfeld upper-half plane}, and in particular is a formal scheme over $Spf \ W$; see \cite{BC} for a discussion of this result. 

Next, we recall the `alternative' description of the reduced locus of $\D$ as given in \cite{KRdrin}. Crucial to this description is the following theorem:
\begin{theorem}[Drinfeld, cf. \S III.4 of \cite{BC}] \label{DrinfeldPolThm}
Suppose $(X, \iota_X)$ is a special formal $\OBp$-module, over any base $S \in \mathbf{Nilp}$. Then there exists a principal polarization $\lambda_X^0$ on $X$ such that 
\begin{equation} \label{PrinPolRosatiEqn}
 (\lambda^0_X)^{-1}  \circ \iota_X(b)^{\vee}  \circ \lambda^0_X \ = \ \iota_X \left( \Pi \ b^{\iota}  \ \Pi^{-1} \right) \qquad \text{for all } b \in \OBp;
 \end{equation}
 here the map $b \mapsto b^{\iota}$ is the involution of $B$. Moreover, $\lambda^0_X$ is unique up to multiplication by $\Z_p^{\times}$. 
\end{theorem}

For the base point $\X$, we shall fix once and for all a polarization $\lambda_{\X}^0$ as in Drinfeld's theorem. Then, for any point $(X, \iota_X, \rho_X) \in \D(S)$, there is a \emph{unique} principal polarization $\lambda_X^0$ satisfying \eqref{PrinPolRosatiEqn}, and such that the diagram 
\[
\begin{CD}
X \times_S \overline{S} @> \lambda_X^0 \times \overline{S} >> X^{\vee} \times_S \overline{S} \\
@V \rho_X VV @VV \rho_X^{\vee}V \\ 
\X \times_{\Fpbar} \overline{S}@>> \lambda_{\X}^0 \times\overline{S}> {\X}^{\vee} \times_{\Fpbar}\overline{S}
\end{CD}
\]
commutes. Thus, for any point $(X, \iota_X, \rho) \in \D(S)$, we may associate another polarization $\lambda_{X, \phi}$ by the formula
\begin{equation} \label{hermPolDefEqn}
\lambda_{X, \phi} := \ \lambda_X^0  \circ  \iota_X \left( \Pi \ \phi( \delta) \right) . 
\end{equation}
Note that the Rosati involution $*$ induced by $\lambda_{X, \phi}$ has the property 
\[ \iota_X \left( \phi (a ) \right)^* = \iota_X \left( \phi ({a'}) \right), \qquad \text{for all } a \in o_{k,p}.  \]

Let $M(\X)$ be the Dieudonn\'e module of $\X$ over $W = W(\Fpbar)$, with Frobenius and Verschiebung operators denoted by $F$ and $V$ respectively. We have a  decomposition 
\[ M(\X) = M(\X)_0 \oplus M(\X)_1, \]
where 
\begin{equation} \label{DieuModGrading} 
	M(\X)_i := \left\{ m \in M(\X) \ | \ (\iota_{\X} \circ \phi)( a) \cdot m = \tau_i(a) m, \text{ for all } a \in o_{k,p} \right\}.
 \end{equation}
Let $N(\X) := M(\X) \otimes_{\Z} \Q$ be the rational Dieudonn\'e module, with induced decomposition
\[ N(\X)  = N(\X)_0 \oplus N(\X)_1 . \]
Note that the operator $p V^{-2}= V^{-1}F: N(\X)_0 \to N(\X)_0$ is a $\sigma^2$-linear operator, and hence the space of invariants 
\begin{equation} C := (N(\X)_0)^{V^{-1}F=1} \end{equation}
is a two-dimensional vector space over $k_p$. Here $k_p$ acts on $C$ via the embedding $\tau_0: k_p \to W$.

The polarization $\lambda_{\X, \phi}$, as defined in \eqref{hermPolDefEqn}, induces an alternating pairing
\[ \left\{ \cdot, \cdot \right\}_{\X}: N(\X) \times N(\X) \to W_{\Q} \] such that for all $x, y \in N(\X)$, we have
\[ \left\{ Fx, y \right\}_{\X} = \sigma \left( \left\{ x, Vy \right\}_{\X} \right). \]
Thus, if we define 
\begin{equation} \label{hDef}  h(x, y) := \frac{1}{p \delta} \left\{ x, F y \right\}_{\X}, \end{equation}
it is straightforward to verify that the restriction of $h$ to $C$ defines a hermitian form; we denote this restricted form again by $h$. We shall shortly see that $(C,h)$ is in fact  \emph{split}, cf.\ Remark \ref{epsilonRmk}. Let $q(x) = h(x,x)$ denote the corresponding quadratic form. 

\begin{definition} 
(i) If $L$ is a $W$-lattice in $N(\X)_0$, we let 
\[ L^{\sharp} := \{ n \in N(\X)_0 \ | \ h(n, L) \subset W \}. \]
 Note that $(L^{\sharp})^{\sharp} = p V^{-2} L.$
Similarly, if $\Lambda \subset C$ is an $o_{k,p}$-lattice, we set 
\[  \Lambda^{\sharp} := \{ v \in C \ | \ h(v, \Lambda) \subset o_{k,p} \}.\]

(ii) Suppose $\Lambda \subset C$. We say $\Lambda $ is a \emph{vertex lattice} of type 0 (resp. type 2) if $\Lambda^{\sharp} = \Lambda$ (resp. $\Lambda^{\sharp} = p \Lambda$). We shall use the term ``vertex lattice'' to mean a vertex lattice of type 0 or 2.

(iii) Let $\mathcal B$ denote the Bruhat-Tits tree for $SU(C)$, which is a graph with the following description. The vertices are vertex lattices, and edges can only occur between vertex lattices of differing type. Two vertex lattices $\Lambda$ and $\Lambda'$ of type 0 and 2 respectively are joined by an edge if and only if
\[ p \Lambda' \subset \Lambda \subset \Lambda', \]
where the successive quotients are $\F_{p^2}$ vector spaces of dimension 1. In particular, this graph is a $p+1$-regular tree.
\end{definition}

Suppose $x = ( X, \iota_X, \rho_X) \in \D(\F)$. We may use the quasi-isogeny $\rho_X$ to identify the Dieudonn\'e module $M(X)$ as a $W$-lattice inside of $N(\X)$. Furthermore,  as $\rho_X$ is $o_{k,p}$-linear, we have $M(X)_i = M(X) \cap N(\X)_i$ for $i = 0,1$, where $M(X)_i$ is defined in the same way as \eqref{DieuModGrading}. Hence, to any point $x$, we may associate a chain of $W$-lattices $B \subset A$, where
\[ B = M(X)_0, \qquad A = (VM(X)_1)^{\sharp}.  \]
By \cite[Corollary 2.3]{KRdrin}, we have either $ B^{\sharp}=B$ or $A^{\sharp} = pA$, or both. If both conditions are satisfied, then we say the point $x$ is \emph{superspecial}; otherwise, $x$ is \emph{ordinary}. We say a point is \emph{special} if both $A$ and $B$ are $pV^{-2}$-invariant, so in particular superspecial points are special.

This construction yields a bijection between $\D(\F)$ and pairs of $W$-lattices $B \subset A$ such that either $B^{\sharp} = B$ or $A^{\sharp} = p A$. 
Moreover, if $B^{\sharp} = B$, then $B = \Lambda \otimes_{o_{k,p}} W $ for some vertex lattice $\Lambda$ of type 0; on the other hand, if $A^{\sharp} = pA$, then $A = \Lambda' \otimes W$ for a vertex lattice $\Lambda'$ of type 2 , cf.\ \cite[Corollary 2.3]{KRdrin}.
 
Suppose $\Lambda$ is a vertex lattice of type 0. We may define a map
\[ \mathbb P_{\Lambda}(\F) := \mathbb P( p^{-1} \Lambda / \Lambda)(\F) \to \D(\F), \]
 by sending a line $\ell \subset (p^{-1} \Lambda / \Lambda) \otimes_{\F_{p^2}} \F$ to the pair of lattices $B \subset A$, where $B = \Lambda_W = \Lambda \otimes W$, and $A $ is the inverse image of $\ell$ in $p^{-1}\Lambda_W$.

If $\Lambda'$ is a vertex lattice of type 2, we obtain a map
\[ \mathbb P_{\Lambda'}(\F) := \mathbb P( \Lambda' / p\Lambda')(\F) \to \D(\F), \]
defined by sending a line $\ell' \subset (\Lambda' / p\Lambda') \otimes_{\F_{p^2}} \F$ to the pair of lattices $B \subset A$, where $A = \Lambda'_W$, and $B $ is the inverse image of $\ell'$ in $\Lambda'_W$.

Note that if $\Lambda$ and $\Lambda'$ are neighbours in $\mathcal B$, i.e if $p \Lambda' \subset \Lambda \subset \Lambda'$, then the lines 
\[ \ell = \Lambda' \otimes_{o_{k,p}, \tau_0} \F \ \in \ \mathbb P( p^{-1} \Lambda / \Lambda)(\F), \qquad \ell' = \Lambda \otimes \F \ \in \ \mathbb P(  \Lambda' /  p \Lambda')(\F)   \]
define the same point of $\D(\F)$; this point is superspecial, and all superspecial points arise in this way. 

By \cite[Proposition 2.4]{KRdrin}, the above maps are induced by embeddings of schemes over $\F$:
\[ \mathbb P_{\Lambda} \to \D_{red}, \qquad \Lambda \text{ a vertex lattice}, \]
where $\D_{red}$ is the underlying reduced subscheme of the formal scheme $\D$, and the collection of such maps, as $\Lambda$ varies among the vertex lattices, yield a cover of $\D_{red}$ by projective lines. 

\begin{remark} \label{epsilonRmk}
(i) In \cite[\S 1]{KRpadic}, there is a similar description of the special fibre of $\D$, but it is given in terms of homothety classes of $\Z_p$-lattices. These two descriptions are essentially the same; to start, note that the operator
 \[ \epsilon := \Pi_{\X}^{-1} \circ V : N(\X) \to N(\X), \qquad \Pi_{\X} := \iota_{\X}(\Pi) \]
is $0$-graded and commutes with $F^{-1} V$, and hence restricts to a Galois-semilinear operator $ \epsilon: C \to C$. 

Without loss of generality, we may assume $ \X = \Y \times \Y$, where $\Y$ is a supersingular $p$-divisible group of height 2 and dimension 1 over $\F$. 
Then the Dieudonn\'e module $M(\X)$ has a basis  $\{ e_0, e_1, f_0, f_1 \}$ consisting of vectors that are both $ F^{-1} V$- and $\epsilon$-invariant, and
such that 
\[ M(\X)_0 = W \cdot e_0 \oplus W \cdot f_1, \qquad M(\X)_1 = W \cdot e_1 \oplus W \cdot f_0, \]
cf. \cite[\S III.4.5]{BC}. 
The set $\{e_0, f_1 \}$ is an $\epsilon$-invariant $k_p$-basis for $C$ such that $h(e_0, f_1) = \delta$ and $h(e_0, e_0) = h(f_1, f_1) = 0$, and in particular, we see that $C$ is split.
Let $\Lambda_0 := span_{o_{k,p}}(e_0, f_1)$, which is a vertex lattice of type 0, i.e. $\Lambda_0^{\sharp} = \Lambda_0$.

Now suppose $\gamma \in  End_{\OBp}(\X) \otimes_{\Z_p} \Q_p$. Then $\gamma$ also acts on $ C$,
 and is described by a matrix of the form
\[ [\gamma]  = \begin{pmatrix} a & b \\ c & d \end{pmatrix} , \qquad a,b,c,d \in \Q_p. \]
with respect to the basis $\{e_0, f_1 \}$.
In particular (at least, when $\gamma$ is invertible), the matrix $[\gamma]$ lies in $ GU(C)$ with $[\gamma]^* = det(\gamma) \cdot [\gamma]^{-1}$. One can then verify directly that the map $\gamma \mapsto [\gamma]$ induces an isomorphism
\begin{equation} \label{SU(C)Eqn}
\left\{ \gamma \in End_{\OBp}(\X)_{\Q_p}, \ det(\gamma) = 1 \right\} \isomto SU(C). 
\end{equation}
Suppose $[L]$ is a homothety class of $\Z_p$-lattices in $M(\X)_0^{\epsilon = 1} = C^{\epsilon = 1} \simeq (\Q_p)^2$, and $L \in [L]$. Then there exists $\gamma \in End_{\OBp}(\X)_{\Q_p}$ such that $ \gamma(L) = L_0 := span_{\Z_p} \{ e_0, f_1 \} $. Set $\Lambda := L \otimes_{\Z_p} o_{k,p}$, so that 
\[ \Lambda^{\sharp} =( [\gamma]^*)^{-1} \cdot \Lambda_0^{\sharp} = det(\gamma)^{-1} [\gamma] \cdot \Lambda_0 = det(\gamma)^{-1} \Lambda.  \]
Hence, by scaling by a power of $p$, there is a unique representative $L \in [L]$ such that $L \otimes o_{k,p}$ is a vertex lattice, whose type depends on the parity of $ord_p det(\gamma)$.

Conversely, suppose $\Lambda\subset C$ is a vertex lattice. Since $SU(C)$ acts transitively on the set of lattices of a given type, there exists an element $[\gamma] \in SU(C)$ such that $[\gamma] \cdot \Lambda$ is equal to one of $\Lambda_0$ or $\Lambda_0':= span_{o_{k,p}} (p^{-1} e_0 , f_1)$, depending on the type of $\Lambda$. In either case, by \eqref{SU(C)Eqn}, the transformation $[\gamma]$ commutes with $\epsilon$, and so $\Lambda$ admits an $\epsilon$-invariant basis. Thus the $\Z_p$-lattice   $\Lambda^{\epsilon = 1}$  determines a homothety class of  $\Z_p$-lattices $ [L] := [\Lambda^{\epsilon = 1}]$. 

(ii) Suppose $x = (X, \iota_X, \rho_X) \in \D(\F)$ , and let $M = M_0 \oplus M_1$ denote its Diedonn\'e module, where each $M_i$ is viewed as a $W$-lattice in $N(\X)_i$. Each lattice $M_i$ is self-dual with respect to the pairing induced by the polarization $\lambda_X^0$ described in Theorem \ref{DrinfeldPolThm}, and it is then a straightforward calculation to show
\[ M_0^{\sharp} = \Pi^{-1} V M_0, \qquad \text{and} \qquad (V M_1)^{\sharp} = \Pi^{-1} M_1,\]
where the $\sharp$ denotes the dual with respect to the hermitian form $h$. 

Recall that $x$ is said to be ``0-critical" in the sense of \cite{BC} if and only if $\Pi M_0 = V M_0$, which is equivalent to the relation $M_0^{\sharp} = M_0$, i.e. $x \in \mathbb P_{\Lambda}(\F)$ for a vertex lattice $\Lambda$ of type 0. 

 Similarly $x$ is ``1-critical" if and only if $\Pi M_1 = V M_1$, which is equivalent to the relation $((V M_1)^{\sharp})^{\sharp} = p (V M_1)^{\sharp}$. This last condition is then equivalent to the condition $x\in \mathbb P_{\Lambda'}(\F)$ for a vertex lattice $\Lambda'$ of type 2.
 
In particular, this discussion implies that our use of the terms ``ordinary'', ``special'' and ``superspecial'' coincides with their use in \cite[\S1]{KRpadic}.  $\diamond$
\end{remark}
 
\end{subsection}

\begin{subsection}{The structure of a local unitary cycle} \label{locEqnsSec}
We begin this section by defining the \emph{(local) unitary special cycles}, whose construction is due to Kudla and Rapoport \cite{KRloc}. 

To start, fix a triple $\underline \Y = (\Y, i_{\Y}, \lambda_{\Y})$ over $\F$ consisting of 
\begin{compactenum}[(i)]
\item a (supersingular) $p$-divisible group $\Y$ of dimension 1 and height 2 over $\F$;
\item an action $i_{\Y}\colon o_{k,p} \to End(\Y)$ such that on $Lie(\Y)$, this action coincides with the action of $o_{k,p}$ via the embedding $\tau_0: o_{k,p} / (p) \to \F$;
\item and finally, a principal polarization $\lambda_{\Y}$ such that the induced Rosati involution * satisfies
\[ i_{\Y}(a)^* = i_{\Y}( a') , \qquad \text{for all } a \in o_{k,p}. \]
\end{compactenum}

We define two spaces of \emph{special homomorphisms}:
\begin{equation} 
\V^+_{\phi} := \left\{ \mathbf b \in Hom(\Y, \X) \otimes_{\Z_p} \Q_p \ | \ \mathbf b \circ i_{\Y} \left( a \right) = \iota_{\X} \left( \phi(a) \right) \circ \mathbf b \text{ for all }  a \in o_{k,p} \right\} 
\end{equation}
and
\begin{equation}
\V^-_{\phi} := \left\{ \mathbf b \in Hom(\Y, \X) \otimes_{\Z_p} \Q_p \ | \ \mathbf b \circ i_{\Y} \left( a \right) = \iota_{\X} \left( \phi( a') \right) \circ \mathbf b \text{ for all }  a \in o_{k,p} \right\} .
\end{equation}

Using the polarization $\lambda_{\Y}$, and the polarization $\lambda_{\X, \phi}$ as defined in \eqref{hermPolDefEqn}, we may construct natural hermitian forms $h^{+}$ and $h^-$ on $\V_{\phi}^+$ and $\V_{\phi}^-$ respectively; these are defined by the formulas 
\begin{align} \label{hpmDef} 
h^{+}(\mathbf b_1, \mathbf b_2) &:= \lambda_{\Y}^{-1} \circ \mathbf b_2^{\vee} \circ \lambda_{\X, \phi} \circ \mathbf b_1  \in End_{o_k}(\underline \Y^+) \otimes_{\Z_p} \Q_p \simeq k_p\\
h^{-}(\mathbf b_1, \mathbf b_2) &:= \lambda_{\Y}^{-1} \circ \mathbf b_1^{\vee} \circ \lambda_{\X, \phi} \circ \mathbf b_2 \in End_{o_k}(\underline \Y^+) \otimes_{\Z_p} \Q_p \simeq k_p.
\end{align}
Let $q^{\pm}(\mathbf b) :=  h^{\pm}(\mathbf b, \mathbf b)$ denote the corresponding quadratic forms.

Our next step is to relate the spaces $\V_{\phi}^{\pm}$ to the hermitian space $(C,h)$. Let $M(\Y)$ denote the Dieudonn\'e module over $W = W(\F)$ attached to $\Y$; this is a free $W$-module of rank 2.
As before, we have a grading $M(\Y) = M(\Y)_0 \oplus M(\Y)_1$, where 
\[M(\Y)_i := \{ m \in M(\Y) \ | \ i_{\Y}(a) \cdot m = \tau_i(a) m, \text{ for all } a \in o_{k,p} \} .\]
Moreover, we may choose generators $f_0$ and $f_1$ for $M(\Y)_0$ and $M(\Y)_1$ respectively, such that $V f_0 = f_1, V f_1 = p f_0$, and that the alternating form $\{ \cdot, \cdot \}_{\Y}$ defined by the polarization $ \lambda_{\Y}$ satisfies
\begin{equation} \label{pairingOnY}
 \{ f_0, f_1 \}_{\Y} = \delta ;
 \end{equation}
see \cite[Remark 2.5]{KRloc}.

Suppose that $\mathbf b \in \V_{\phi}^+$. Abusing notation,  we denote the corresponding map on (rational) Dieudonn\'e modules again by 
\[ \mathbf b: N(\Y) = M(\Y)\otimes \Q \to N(\X).\]
 Let $b := \mathbf b(f_0)$. Then since $\mathbf b$ is $o_{k,p}$-linear, we have that $b \in N(\X)_0$, and furthermore, 
\[ V^{-1}F b = p V^{-2} b = p \ \mathbf b( V^{-2} f_0 ) = b \]
so $b \in C$. Finally, we note that 
\[ \mathbf b( f_1) = \mathbf b (V f_0 ) = V \mathbf b(f_0) = Vb, \]
and so $\mathbf b$ is determined by $b$.
We therefore obtain an isomorphism
\begin{equation} \label{varphiPlusDef}
\varphi^+: \V_{\phi}^+ \to C, \qquad \mathbf b \mapsto b :=  \mathbf b(f_0).
\end{equation}
One readily checks that  $q( \varphi^+ \mathbf b ) = p^{-1} q^+(\mathbf b)$, where $q$ and $q^+$ are the quadratic forms on $C$ and $\V_{\phi}^+$ defined by \eqref{hDef} and \eqref{hpmDef} respectively.

In a similar manner, if $\mathbf b \in \V_{\phi}^{-}$, then $b := \mathbf b(f_1) \in C$, and $\mathbf b$ is again determined by $b$. Hence we obtain an isomorphism
\begin{align} \label{varphiMinusDef}
 \varphi^-: \V_{\phi}^- \to C, \qquad \mathbf b \mapsto b := \mathbf b(f_1)
\end{align}
such that  $q( \varphi^- \mathbf b) = q^-(\mathbf b)$. 

Finally, we come to the definition of the unitary special cycles, as in \cite{KRloc}. Note that the $p$-divisible group $\underline \Y = (\Y, i_{\Y})$, together with its $o_{k,p}$-action, admits a canonical lift $\underline \Y_{W} = (\Y_{W}, i_{\Y_W})$ to $W$. For any $S \in \mathbf{Nilp}$, we let $\underline \Y_S = (\Y_S, i_{\Y_S})$ denote the base change $\Y_S =  \Y_{W} \times_W S$, with induced $o_{k,p}$-action. 
\begin{definition} \label{locUSpCycDef}
For $\mathbf b \in \V_{\phi}^{\pm}$, we define the \emph{local unitary special cycle} $Z(\mathbf b)$ by the following moduli problem: for $S \in \mathbf{Nilp}$, let $ Z(\mathbf b)(S) $ denote the set of points $(X, \iota_X, \rho) \in \D(S)$ such that quasi-isogeny
\[ \rho^{-1}  \circ \mathbf b: \Y \times_{\F}  \ \overline S \to X \times_S \overline S \]
lifts to a morphism of $p$-divisible groups $\Y_S \to X$. 
\end{definition}
\begin{remark} 
(i) If such a lift exists, then it is unique, by rigidity for $p$-divisible groups. In particular, if the lift exists, then it is $o_{k,p}$-linear (resp. 		$o_{k,p}$-antilinear) whenever $\mathbf b \in \V_{\phi}^+$ (resp. $\mathbf{b} \in \V_{\phi}^-$). 

(ii) By \cite[Proposition 2.9]{RZ}, the moduli problem $Z(\mathbf b)$ is represented by a \emph{closed formal subscheme} of $\D$.	
%
$\diamond$
\end{remark}

Our first result is the following description of the $\F$-points of the special cycles.

\begin{proposition} \label{FPointsProp}
(i) Suppose $\mathbf b \in \V_{\phi}^-$, and $b = \varphi^- \mathbf b \in C$ is the corresponding vector, cf. \eqref{varphiMinusDef}. Let $\Lambda$ be any vertex lattice. Then
	\[ Z(\mathbf b)(\F) \cap \mathbb P_{\Lambda}(\F) = 
	\begin{cases} 	\emptyset,		 & \text{if } b \notin \Lambda, \\
				\text{a single point } \{ x \}, & \text{if } b \in \Lambda - p \Lambda, \\
				\mathbb P_{\Lambda}(\F) , & \text{if } b \in p \Lambda. 
	\end{cases}
	\]
	If the second case above occurs, then $\Lambda$ is necessarily of type 0, and the point $x$ is special. Moreover $x$ is superspecial if and 		only if $ord_p q^-(\mathbf b) > 0$. 
	
(ii) Similarly, for $\mathbf b \in \V_{\phi}^+$ with $b = \varphi^+ \mathbf b$, cf. \eqref{varphiPlusDef}, we have
	\[ Z(\mathbf b)(\F) \cap \mathbb P_{\Lambda}(\F) = 
	\begin{cases} 	\emptyset,		 & \text{if } b \notin \Lambda, \\
				\text{a single point } \{ x \}, & \text{if } b \in \Lambda - p \Lambda, \text{ with } \Lambda^{\sharp} = p \Lambda\\
				\mathbb P_{\Lambda}(\F), & \text{if } b \in \Lambda \text{ with } \Lambda^{\sharp} = \Lambda, \\
				\mathbb P_{\Lambda}(\F) , & \text{if } b \in p \Lambda, \ \Lambda  \text{ arbitrary.}
	\end{cases}
	\]
In the second case, the unique point $x$  is special, and is superspecial if and only if $ord_pq^+( \mathbf b )> 0$.

\begin{proof}
(i) Suppose $\mathbf b \in \V^-_{\phi}$. We observe that a point $x = (\underline X, \rho_X) \in \D(\F)$ is in $Z(\mathbf b)( \F)$ if and only if, upon identifying $M(X)$ with a $W$-lattice in $N(\X)$, we have 
\begin{align*}\mathbf b ( M(\Y) )  \subset M(X)  &\iff \mathbf b(f_0) \in M(X)_1, \text{ and } \mathbf b(f_1) \in M(X)_0 \\
&\iff b = \mathbf b(f_1) \in VM(X)_1;
\end{align*}
the last equivalence follows from the relation $V\mathbf b(f_0 ) =  \mathbf b(f_1)$.
Recall also that $\mathbf b(f_0) \in C = N(\X)_0^{FV^{-1}}$, so we have 
\[ x = (\underline X, \rho_X) \in Z(\mathbf b)(\F) \ \iff \ b \in VM(X)_1 \ \cap \ C.\]

Now suppose $x \in \mathbb P_{\Lambda}(\F) \cap Z(\mathbf b)(\F)$, with $\Lambda^{\sharp} = p \Lambda$. By construction, $x$ corresponds to a lattice pair $B \subset A$, where 
\[ A = VM(X)_1^{\sharp} = \Lambda \otimes_{\tau_0} W .\]
 
On the other hand, note
\[ p A = A^{\sharp} = (VM_1(X)^{\sharp})^{\sharp} = FM_1(X). \]

Thus, as $x \in Z(\mathbf b)(\F)$ as well, then 
\[ b \in VM_1(X) \cap C = FM_1(X) \cap C = pA \cap C. \]
However, as $pA = p \Lambda_W$, the above line is true for \emph{all} $x \in \mathbb P_{\Lambda}(\F)$, as soon as it is true for a single point. 
Hence we have $\mathbb P_{\Lambda}(\F) \subset Z(\mathbf b)(\F)$.

Now suppose $x \in \mathbb P_{\Lambda}(\F)$ with $\Lambda^{\sharp} = \Lambda$. By construction, this means $M(X)_0$ is identified with $\Lambda \otimes_{\tau_0} W$, and so if $\mathbb P_{\Lambda}(\F) \cap Z(\mathbf b) \neq \emptyset$, then we must have $b \in \Lambda$.  Furthermore, any $x \in \mathbb P_{\Lambda}(\F)$ is determined by the sequence of inclusions of $\F$-codimension 1
\[ \begin{matrix} p M(X)_0 & \subset & VM(X)_1 & \subset & M(X)_0 \\
 || & & & & || \\
 p \Lambda \otimes W & & & & \Lambda \otimes W.
 \end{matrix}
\]
Hence, if $b \in p \Lambda$, then $b \in VM(X)_1$ for all $(\underline X, \rho_X) \in \mathbb P_{\Lambda}(\F)$, and so
\[ Z(\mathbf b)(\F) \cap \mathbb P_{\Lambda}(\F) = \mathbb P_{\Lambda}(\F). \]

If on the other hand $b \in \Lambda \setminus p \Lambda$, and $\mathbb P_{\Lambda}(\F) \cap Z(\mathbf b)(\F) \neq \emptyset$, then this intersection necessarily consists of a single point $x = (X, \rho_X)$: namely, the unique point with 
\[ VM(X)_1 = W \cdot b + p \Lambda_W \subset \Lambda_W. \]
Note in this case, $ (p V^{-2}) VM(X)_1= VM(X)_1$ as both $\Lambda_W$ and $b$ are $p V^{-2}$-invariant.

Now by construction, we have 
\[ \Lambda_W \subset VM(X)_1^{\sharp} \subset p^{-1} \Lambda_W . \]
If $ord_p q(b) = 0$, then $p^{-1} b \notin VM(X)_1^{\sharp}$, and so the point is ordinary. On the other hand, if $ord_pq(b) > 0$, then 
\[ VM(X)_1^{\sharp} = W \cdot p^{-1} b + \Lambda_W = p^{-1} VM(X)_1 = p^{-1} \left( VM(X)_1^{\sharp} \right)^{\sharp}, \]
and so this point is superspecial.

(ii) Suppose $\mathbf b \in \V_{\phi}^+$, and let $b = \mathbf b( f_0) = \varphi^+(\mathbf b)$. A point $x = (\underline X, \rho_X)$ lies in $Z(\mathbf b)(\F)$ if and only if $b \in M(X)_0$. The lemma follows via a similar argument to the previous case. 

\end{proof}

\end{proposition}

\begin{lemma} \label{centralLatticeLemma}
Fix an $\epsilon$-invariant basis $ \{ v_0, v_1 \}$ of $C$ (cf. Remark \ref{epsilonRmk}) such that 
\[ h(v_0, v_0) = h(v_1, v_1) = 0, \qquad h(v_0, v_1) = -h(v_1, v_0) = \delta .\]
Suppose that $b = a_0 v_0 + a_1 v_1$, with $ord_pq(b)$ either $0$ or $-1$. Set $b' := \epsilon(b) =  a_0' v_0 +  a_1' v_1$ and
\[ \Lambda_b := span_{o_{k,p}} \left\{  b, b' \right\}.\]
If $ord_pq(b) = 0$,  then $\Lambda_b$ is of type 0 (i.e $\Lambda_b^{\sharp} = \Lambda_b$), and is the unique lattice of type 0 such that $b \in \Lambda_b - p \Lambda_b$. Similarly, if $ord_pq(b) = -1$, then $\Lambda_b$ is of type 2, and is the unique lattice of type 2 such that $b \in \Lambda_b - p \Lambda_b$. 
\begin{proof}
Suppose that $ord_p q(b) = 0$, and $\Lambda$ is a type 0 lattice with $b \in \Lambda - p \Lambda$. Then there exists an element $\gamma \in SU(C)$ such that 
\[ \gamma \cdot \Lambda =\Lambda_b, \]
as $SU(C)$ acts transitively on the set of type 0 lattices. 
We may assume without loss of generality that $q(b) = 1$. Note the vectors $b$ and $b'$ form an orthogonal basis for $C$, with $h(b,b) = -h(b',b') = 1$. Let 
\[ [\gamma] =  \begin{pmatrix} x &y \\ w& z \end{pmatrix} \] denote the matrix representation of $\gamma$ with respect to the basis $\{ b, b' \}$. Since $b = \left( \begin{smallmatrix} 1 \\ 0 \end{smallmatrix} \right) \in \Lambda$, we have $x, w \in o_{k,p}$. On the other hand, the equation $\gamma \cdot \gamma^* =1$ implies
\[ x  x' - y  y' = -w  w'  + z  z' = 1, \qquad y  w' = x  z' ,\]
which implies that $y, w \in o_{k,p}$ as well. Hence $\gamma$ and  $\gamma^{*}$ stabilize $\Lambda_b$, which yields the result $\Lambda = \Lambda_b$. 

The proof in the case $ord_p q(b) = -1$ is similar. 
\end{proof}

\end{lemma}

\begin{remark} \label{centralLatticeRemark}
Suppose $\mathbf b \in \V_{\phi}^{\pm}$, with $ord_p q ^{\pm} \mathbf b = 0$. Let $b = \varphi^{\pm} \mathbf b \in C$ denote the corresponding vector. Then by Lemma \ref{centralLatticeLemma}, there is a \emph{unique} lattice $\Lambda_{\odot}$ such that $b \in \Lambda_{\odot} - p \Lambda_{\odot}$, and $\Lambda_{\odot}$ is of type 0 (resp. of type 2) if $\mathbf b \in \V_{\phi}^-$ (resp. $\mathbf b \in \V_{\phi}^+$). 
Combining this observation with Proposition \ref{FPointsProp}, we find that 
\[ Z(\mathbf b)(\F) = \{ x \} \]
is a \emph{single ordinary point} lying in the component $\mathbb P_{\Lambda_{\odot}}(\F)$. $\diamond$
\end{remark}

Our goal for the remainder of this section is to give a complete description of the special cycles $Z(\mathbf b)$ as cycles on $\D$, as in Theorem \ref{cycleDecompThm} below. We shall do this by writing down equations using the (formal) affine open cover  described in \cite[\S1]{KRpadic}, which consist of affine schemes of two types: 
\begin{enumerate}
\item $ \widehat\Omega_{\Lambda}^{ord}\simeq Spf W[T, (T^p - T)^{-1}]^{\vee}$, for each vertex lattice $\Lambda$, and
\item $ \widehat\Omega_{[\Lambda, \Lambda']} \simeq Spf W[T_1, T_2, (T_1^{p-1} - 1)^{-1}(T_2^{p-1} - 1)^{-1} ] ^{\vee} $, for each pair of neighbouring vertex lattices $\Lambda$ and $\Lambda'$. 
\end{enumerate}
In both cases above, the superscript ${}^{\vee}$ denotes completion along the ideal generated by $p$, and the isomorphisms are determined by a choice of basis for $\Lambda$. The underlying set of $ \widehat\Omega_{\Lambda}^{ord} $ is the set of ordinary points $\mathbb P_{\Lambda}^{ord}$ in $\mathbb P_{\Lambda}$ (that is, the complement of the superspecial points), while the underlying set of $\widehat\Omega_{[\Lambda, \Lambda']}$ is the union 
\[ \mathbb P_{\Lambda}^{ord} \cup \mathbb P_{\Lambda'}^{ord} \cup \{ x \} \]
where $x$ is the superspecial point at the intersection of $\mathbb P_{\Lambda} \cap \mathbb P_{\Lambda'}$. 
For $\Lambda$ a type 0 lattice, with neighbour $\Lambda'$ of type 2, we have open immersions (cf.\ \cite[\S 1]{KRpadic}),
\[ \widehat\Omega^{ord}_{\Lambda} \to \widehat\Omega_{[\Lambda, \Lambda']}, \qquad \text{induced by } T_0 \mapsto T, \ T_1 \mapsto p T^{-1} 
\]
and 
\[ \widehat\Omega^{ord}_{\Lambda'} \to \widehat\Omega_{[\Lambda, \Lambda']}, \qquad \text{induced by } T_0 \mapsto pT^{-1}, \ T_1 \mapsto  T.  
\]


We first consider a type 0 lattice $\Lambda = \Lambda^{\sharp}$. Let $\mathbf b \in \V_{\phi}^{\pm}$, with $ord_p q^{\pm}(\mathbf b) \geq 0$ and such that $Z(\mathbf b) \cap \widehat \Omega_{\Lambda}^{ord} \neq \emptyset$.  Let $b = \varphi^{\pm} \mathbf b \in C$ the corresponding vector; by Proposition \ref{FPointsProp}, we must have $b \in \Lambda$. 
Fix an $\epsilon$-invariant basis $\{ v_0, v_1 \}$ of $\Lambda$, such that with respect to this basis, the hermitian form $h$ has matrix
\begin{equation} \label{hNiceBasisEqn}
	h \sim \begin{pmatrix} & \delta \\ - \delta & \end{pmatrix},
\end{equation}
and write
\begin{equation}
	b = a_0 v_0 + a_1 v_1, \qquad \qquad \text{where } a_0, a_1 \in o_{k,p}.
\end{equation}
%
%
\begin{proposition} \label{ordEqnType0Prop} Let $\mathbf b \in \V_{\phi}^{\pm}$, and suppose $Z(\mathbf b) \cap  \widehat  \Omega^{ord}_{\Lambda} \neq \emptyset$. Then, with notation as in the previous paragraph,
\[ Z(\mathbf b) \cap \widehat  \Omega^{ord}_{\Lambda} \simeq Spf W[T, (T^p - T)^{-1}]^{\vee} / (f), \]
where 
\begin{equation}
	f:= 
		\begin{cases} 	a_0T + a_1,	&	\text{if } \mathbf b \in \V_{\phi}^- \\
					p(a_0' T + a_1'),			&	\text{if } \mathbf b \in \V_{\phi}^+.
		\end{cases}
\end{equation}

\begin{proof}

The proof of this proposition is modelled on arguments by Terstiege found in the proofs of \cite[Propositions 2.8 and 4.5]{terstiegeAnti}. We consider the following cases:

\vspace{10pt}
\noindent\fbox{\textbf{Case 1:} $\mathbf b \in \V_{\phi}^-$, with $b = \varphi^- \mathbf b \in \Lambda$.}

Fix a point $x \in Z(\mathbf b)(\F) \cap \mathbb P_{\Lambda}(\F)$, corresponding to $(\overline X, \iota_{\overline X}, \rho) \in \D(\F)$. At present, we do not require that $x$ be ordinary (i.e. $\overline X$ may be superspecial); in fact the argument we are about to discuss applies equally well in the superspecial case, and so it will be expedient to consider this slightly more general setting here. Consider the (complete) local ring
\[ R = \mathcal O_{\D, x}. \] 
Let $m_R$ denote the maximal ideal of $R$, and $I$ denote the ideal corresponding to $Z(\mathbf b)$. 
For the purposes of the proposition, we need to prove that $I$ is generated by the image of $f$ in $R$. 
Note $R$ is Noetherian, which implies that $I$ is complete for the $p$-adic topology. Hence if we set 
\[ R_n := R / p^n R, \qquad I_n = I + p^n R \subset R_n, \]
it will suffice to prove that $I_n$ is generated by the image of $f$ in $R_n$, for each $n$. Let $m_n$ denote the maximal ideal of $R_n$. We set 
\[ A := R_n / m_n I_n, \qquad A' := R_n / I_n. \]
Then the kernel $J := I_n / m_n I_n$ of the projection $ A \to A'$ satisfies $J^2 = 0$, and hence is endowed with a PD structure. Moreover, to prove the proposition, it suffices (by Nakayama) to show that $J$ is generated by the image of $f$ in $A$ (which, abusing notation, we shall henceforth denote as $f$).

Finally, we note that both $A$ and $A'$ can be viewed as $o_{k,p}$-algebras via the fixed embedding $\tau_0: o_{k,p} \to W$ composed with the respective structural morphisms. 

Now associated to the  rings $A$ and $ A'$ are two points in the moduli space $\D$, which in turn correspond to formal $\OB$-modules $X$ and $X'$, both of whose special fibres are equal to $\overline X$. Moreover, by assumption, $X' $ is in $Z(\mathbf b)(A')$.
In other words, the map
\begin{equation*}
	\beta_{\overline X} := \rho^{-1} \circ \mathbf b \ \colon \Y \to \overline X 
\end{equation*}
is a morphism of $p$-divisible groups, which lifts to a morphism
\begin{equation} \label{betaX'}
	\beta_{X'}: \Y_{A'} \to X'. 
\end{equation}
Here  $\Y_{ A'} $ is the base change of the canonical lift $\Y_W$ to  $A'$. 

Let $\mathbb D(\Y_{A'} / \cdot)$ denote the Grothendieck-Messing crystal of $\Y_{A'}$, which carries a $\Z / 2\Z$ grading induced by the action of $o_{k,p}$. In particular, for a $PD$ extension $B \to A'$, the canonical lift of $\Y_{A'}$ (together with its $o_{k,p}$-action) to $B$ is determined by the Hodge filtration:
\begin{equation}
\mathbb D(\Y_{A'}/B)_1 = \mathcal F_{\Y_B} \ \subset  \ \mathbb D(\Y_{A'} / B).
\end{equation}
This is a consequence of the requirement that $o_{k,p}$ act on $Lie(\Y_{B})$  with signature $(1,0)$. 

Turning to the $p$-divisible groups $X$ and $ X'$, we let $\mathbb D( X/ \cdot)$ and $\mathbb D( X' / \cdot)$ denote their respective Grothendieck-Messing crystals. 
We then have a diagram
\begin{equation} \label{GrothMessDiag}
\begin{CD}
\mathbb D( X/ A) & \ \simeq \ & \mathbb D( X'/ A) \\
@. @VV \text{mod } J V \\
@.   \mathbb D( X'/ A') 
\end{CD}
\end{equation}

Let $B \to A'$ be a PD-extension as above. Since $X' \in Z(\mathbf b)(A')$, there is a map of $B$-modules 
\[ \mathbb D (\beta_{X'}) \colon \mathbb D(\Y_{A'} /B) \to \mathbb D(X' / B) \]
induced by \eqref{betaX'}.
By Grothendieck-Messing theory, if $\tilde X$ is a lift of $X'$ to $B$ which  corresponds to an $\OBp$-stable direct summand $\mathcal F \subset \mathbb D(X'/ B)$, then 
\begin{equation*}
\tilde X \in Z(\mathbf b)(B) \qquad \iff \qquad \mathbb D(\beta_{X'})(\mathcal F_{\Y_B}) \subset \mathcal F.			
\end{equation*}
Fixing any basis vector $f_B \in \mathcal F_{\Y_B} = \mathbb D(\Y_{A'}/B)_1$, and recalling that $\mathbf b$ is $o_{k,p}$-antilinear,  the above condition can be rewritten as 
\begin{equation} \label{condOnTildeXEqn}
\tilde X \in Z(\mathbf b)(B) \qquad \iff \qquad \mathbb D(\beta_{X'})(f_B) \in  \mathcal F_0.			
\end{equation}
For convenience, we fix a basis vector $f_1' \in \mathbb D(\Y_{ A'} /  A')_1$, and  a lift  to a basis vector  $f_1 \in \mathbb D(\Y_A / A)_1$. 

We want to express the condition \eqref{condOnTildeXEqn} in terms of the affine chart
\begin{equation} \label{ssCoordEqn}
 \widehat\Omega_{[\Lambda, \Lambda']} \simeq Spf \ W[T_0, T_1, (T_0^{p-1} -1 )^{-1}, (T_1^{p-1} - 1)^{-1} ] / (T_0 T_1 - p)^{\vee};
\end{equation}
note that as $p$ is nilpotent in $A$, the point $X \in \widehat\Omega_{[\Lambda, \Lambda']}(A)$ is simply determined by a map of $W$-algebras
\begin{equation} 
  W[T_0, T_1, (T_0^{p-1} -1 )^{-1}, (T_1^{p-1} - 1)^{-1} ] / (T_0 T_1 - p) \to A.
\end{equation}
Let $t \in A$ denote the image of $T_0$ in $A$, so that $ t^{p-1} - 1 \in A^{\times}$. 

If $\overline X$ is \emph{ordinary}, as in the hypotheses of the proposition, then we may equally well restrict to the smaller open affine neighbourhood
\begin{equation} \label{ordCoordEqn}  \widehat\Omega_{\Lambda}^{ord}  \simeq Spf W[T, (T^p - T)^{-1}]^{\vee},
\end{equation}
via the open immersion
\[ \widehat\Omega^{ord}_{\Lambda} \to \widehat\Omega_{[\Lambda, \Lambda']}, \qquad \text{induced by } T_0 \mapsto T, \ T_1 \mapsto p T^{-1} 
\]
The point $X \in \widehat\Omega_{\Lambda}^{ord}(A)$ then corresponds to the same element $t \in A$ as above, which in this case now satisfies $t \in A^{\times}$. 

%

In any case, by carefully tracing through the construction of the charts \eqref{ssCoordEqn} and \eqref{ordCoordEqn} as they are constructed in \cite{BC}, we find that there exists isomorphisms
\begin{equation} \label{DCrysIsomLambdaEqn}
 \mathbb D( X / A)_0 \simeq \Lambda \otimes_{o_{k,p}} A   , \qquad  \mathbb D( X' / A')_0 \simeq \Lambda \otimes_{o_{k,p}} A' 
\end{equation}
such that 
\begin{compactenum}[(i)]
\item the 0-th graded piece of the Hodge filtration for $X$ is identified with 
\[ (\mathcal F_{X})_0 = span_A \{ v_0 \otimes 1 -  v_1 \otimes t \} \subset  \mathbb D( X / A)_0 \simeq \Lambda \otimes A, \]
where $\{ v_0, v_1 \}$ is the fixed basis of $\Lambda$ as in \eqref{hNiceBasisEqn};
\item the element $\mathbb D(\beta_{X'})(f_1)$ is identified with $ b \otimes 1 \in \Lambda \otimes A$; 
\item the vertical map in \eqref{GrothMessDiag} is the natural map
\[ \Lambda \otimes A \to \Lambda \otimes A' , \qquad l \otimes a \mapsto l \otimes a',\]
where $a \mapsto a'$ is the projection $A \to A'$.
\end{compactenum}

Now suppose $\lie{b}$ is any ideal of $A$ contained in $J$, and let $B := A / \lie b$. Let $t_B$ denote the image of $t$ in $B$ and $ X_B$ the corresponding $B$-valued point of $\D$. We may view $ X_B$ as a deformation of $X'$, which by Grothendieck-Messing theory corresponds to the direct summand 
\[ \mathcal F_{ X_B} = \mathcal F_X \otimes_A B \subset \mathbb D(X'/B), \]
and in particular we have
\[ (\mathcal F_{ X_B})_0 = span_B \{ v_0 \otimes 1 - v_1 \otimes t_B \} \subset \mathbb D(X'/B)_0 = \Lambda \otimes B. \]
Then in light of \eqref{condOnTildeXEqn},  and noting that $\mathbf b $ is $o_{k,p}$-antilinear, we have 
\begin{align*}
	 X_B \in Z(\mathbf b)(B) \ \iff& \ \mathbb D(\beta_{X'}/B)(f_1) \otimes 1 \in  (\mathcal F_{\tilde X})_0 \\
								\iff&\ b \otimes 1 \in span\{ v_0 \otimes 1- v_1 \otimes t_B \} \text{ in } \Lambda \otimes B.
\end{align*}
Recall that we had written $b = a_0 v_0 + a_1 v_1$. Hence the last condition is equivalent to 
\[ a_0 t_B + a_1 = 0 \text{ in } B. \]
Applying the necessity of this condition to the case $\lie b = J$, so $ B = A/J$, we find that 
\[ f = a_0 t + a_1 \in J.\] On the other hand, by the sufficiency of the above condition, the map $Spf(A/ (f)) \to \widehat\Omega_{\Lambda}^{ord}$ factors through $Z(\mathbf b)$.
But by definition $J$ is the smallest ideal of $A$ such that $Spf(A/J) \to \widehat\Omega_{\Lambda}^{ord} $ factors through $Z(\mathbf b)$, and so $J = (f)$; under the hypotheses of the proposition (i.e.\  $\overline X$ is  ordinary), this concludes the proof. 

\vspace{10pt}
\noindent\fbox{\textbf{Case 2:} $ \mathbf b \in \V_{\phi}^+$ with $b = \varphi^+ \mathbf b \in \Lambda$. }

Let
$ r = r(b, \Lambda) := max \{ n \mid p^{-n} b \in \Lambda \}, $
and write
\[ b = p^r b_0 = p^r(\alpha_0 v_0 + \alpha_1 v_1), \]
where $b_0 = \alpha_0 v_0 + \alpha_1 v_1 \in \Lambda - p\Lambda$, and so in particular at least one of $\alpha_0, \alpha_1$ is a unit in $o_{k,p}^{\times}$.  

As a first step, we shall prove that every $W_{r+1} := W / (p^{r+1})$ valued point of $\widehat \Omega^{ord}_{\Lambda}$ belongs to the special cycle $Z(\mathbf b)$. 
To this end, suppose $(X, \rho_X) \in \widehat \Omega_{\Lambda}^{ord}(W_{r+1})$, and let $\overline X \in \Omega_{\Lambda}^{ord}(\F) $ denote its reduction modulo $p$. Let 
\[ \overline M = M(\overline X) = \overline M_0 \oplus \overline M_1 \]
denote the Dieudonn\'e module of $\overline X$, endowed with the grading induced by the action of $o_{k,p}$. As the projection $W_{r+1} \to \F$ has kernel generated by $p$, it is equipped with a PD structure. Hence, via Grothendieck-Messing theory, the point $X$ corresponds to an $\OBp$-stable summand
\[ \mathcal F \ \subset \ \mathbb D( \overline X / W_{r+1}) \ =  \ \overline M \otimes_W W_{r+1} = \overline M / p^{r+1} \overline M, \]
such that 
\[ \mathcal F \otimes_{W_{r+1}} \F = V\overline M / p \overline M. \]

By Proposition \ref{FPointsProp}, we have $\overline X \in Z(\mathbf b)( \F)$, and so the map $\beta:= \rho_X^{-1} \circ \mathbf b: \Y \to \overline X$ induces a morphism of Dieudonn\'e modules
\[ \beta:  M(\Y) = W \cdot f_0 \oplus W \cdot f_1 \  \to \  \overline M. \]
By definition, we have 
\[ \beta(f_0) = b  \in \Lambda \otimes W = \overline M_0, \qquad \beta(f_1) =  V \beta(f_0) = Vb \in \overline M_1. \]
The morphism $\beta$ also induces a map on crystals
\[ \mathbb D(\beta/W_{r+1}):  \ M(\Y) \otimes_W W_{r+1} = \mathbb D(\Y/W_{r+1}) \ \to \ \mathbb D(\overline X / W_{r+1}) = \overline M \otimes_W W_{r+1}. \]

Recall that the direct summand 
\[ \mathcal F_{\Y} \subset M(\Y) \otimes_W W_{r+1} \]
 corresponding to the lift $\Y_{W_{r+1}}$ of $\Y$ is simply the rank-1 module  $\mathcal F_{\Y} = span_{W_{r+1}} \{ f_1 \otimes 1 \}$. Hence, by Grothendieck-Messing, we have
\begin{align}
 X \in Z(\mathbf b)(W_{r+1}) 	\iff& \mathbb D(\beta/W_{r+1})(f_1) \in \mathcal F_1 \notag \\
 \label{WrPointsEqn}		\iff& (Vb) \otimes 1 \in \mathcal F_1 \text{   in } \overline M_1 \otimes W_{r+1}. 
\end{align}
We shall show that this latter condition always holds for any $X \in \widehat\Omega_{\Lambda}^{ord}(W_{r+1})$. 

Consider the element $b_0  = p^{-r} b \in \Lambda \setminus p \Lambda$. 
If $Vb_0 \in p \overline M_1$, then 
\[ Vb \otimes 1 = p^r \cdot Vb_0 \otimes 1 \  \in \ p^{r+1} \cdot(\overline M_1 \otimes W_{r+1}) = \{ 0 \},\]
 and so \eqref{WrPointsEqn} holds trivially. If on the other hand $Vb_0 \in \overline M_1 - p \overline M_1$, then 
\[ \text{span}_{\F} ( V b_0  + p \overline M_1) = V \overline M_0 / p \overline M_1 = \mathcal F_1 \otimes_{W_{r+1}} \F;\]
in other words, the image of $V b_0$ in $\overline M_1 / p \overline M_1$ is a basis vector for $\mathcal F_1 \otimes \F$. 
Hence, there exists $\alpha \in \overline M_1 \otimes W_{r+1}$ such that
\[ V b_0 \otimes 1 + p \alpha \in \mathcal F_1, \]
and so 
\[ p^r (V b_0 \otimes 1) + p^{r+1} \alpha = Vb \otimes 1 \in \mathcal F_1, \]
as required. Thus, we have proven
\begin{equation} \label{allWr+1PtsEqn}
Z(\mathbf b) \cap \widehat \Omega_{\Lambda}^{ord}(W_{r+1}) = \widehat \Omega_{\Lambda}^{ord}(W_{r+1}).
\end{equation}

We now determine the local equation of $Z(\mathbf b) \cap \widehat\Omega_{\Lambda}^{ord}$ in $ \widehat \Omega_{\Lambda}^{ord} \simeq Spf W[T, (T^p-T)^{-1}]^{\vee}$, assuming the intersection is non-empty. Note that as $\Lambda = \Lambda^{\sharp}$ is a lattice of type 0, then by Remark \ref{centralLatticeRemark}, we must have $ord_p q^+ (\mathbf b) >0$; otherwise, $Z(\mathbf b)$ would meet the special fibre only at a type-2 lattice. 

Set
\[ \mathbf b' := \Pi_{\X}^{-1} \circ \mathbf b \ \in \ \V_{\phi}^-, \]
so that $ord_pq^-(\mathbf b') = ord_p q^+(\mathbf b) - 1 \geq 0$, and note
\begin{align*}
 \varphi^- \mathbf b' = \mathbf b'(f_1) = \Pi_{\X}^{-1} \mathbf b (f_1) =& V \Pi_{\X}^{-1} \mathbf b(f_0) \\
 =&  \epsilon(b) \\
 =& p^r(\alpha_0' v_0 + \alpha_1' v_1),
\end{align*}
since by assumption the basis $\{ v_0, v_1 \}$ was taken to be $\epsilon$-invariant. 

Furthermore, we have evident inclusions 
\[ Z(\mathbf b') \subset Z(\mathbf b) \subset Z(p \cdot \mathbf b'). \]
If we let $ I \subset W[T, (T^p - T)^{-1}]^{\vee}$ denote the ideal defining $Z(\mathbf b) \cap  \widehat \Omega_{\Lambda}^{ord}$, we then have
\[ \left( p^{r+1}( \alpha_0' T +  \alpha_1') \right) \  \subset  \ I \  \subset \   \left( p^{r}( \alpha_0' T +  \alpha_1') \right), \]
by applying the result from Case 1 to the two antilinear homomorphisms $\mathbf b'$ and $p \cdot \mathbf b '$. 

Recall that we want to prove that $I$ is generated by $f:= p^{r+1}(\alpha_0' T + \alpha_1')$. Thus,  it suffices to prove that $I \subset \left( p^{r+1}( \alpha_0' T + \alpha_1') \right)$. 

To this end, suppose $g \in I$. We write 
\[ g = p^r ( \alpha_0' T +  \alpha_1') g_0, \]
and we need to show that $p $ divides $ g_0$. 

Suppose not. Then the reduction modulo $p$ of $ ( \alpha_0' T +  \alpha_1') g_0$ is a non-zero rational function over $\F$ and so there exists $t \neq 0 \in \F$ such that
\[ ( \alpha_0' t +  \alpha_1') g_0(t) \neq 0 \in \F. \]
Choose any preimage $\tilde t \in W_{r+1}^{\times}$ of $t$. 
Then the map
\[ W[T, (T^p - T)^{-1}]^{\vee} \to W_{r+1}, \qquad T \mapsto \tilde t \]
does \emph{not} factor through $W[T, (T^p - T)^{-1}]^{\vee} / I$, which contradicts the assertion \eqref{allWr+1PtsEqn}. 
Hence, we have that $p$ divides $g_0$, and so 
\[ I = \left( p^{r+1} ( \alpha_0' T +  \alpha_1' ) \right) = \left( p(a_0' T + a_1') \right) \]
 as required.

\end{proof}

\end{proposition}

For a type 2 lattice $\Lambda'$, we describe the analogous result: choose an $\epsilon$-invariant basis $\{ w_0, w_1 \}$ for $\Lambda'$ such that $q(w_0) = q(w_1) = 0 $ and $h(w_0, w_1) = p^{-1} \delta$. Suppose $\mathbf b \in \V_{\phi}^{\pm}$, and let $b = \varphi^{\pm} \mathbf b$  denote the corresponding vector. If $Z(\mathbf b) \cap \widehat\Omega_{\Lambda'}^{ord} \neq \emptyset$, then by Proposition \ref{FPointsProp}, we must have $b \in \Lambda'$, and so we may write
\[ b =a_0 w_0 + a_1 w_1, \qquad a_0, a_1 \in o_{k,p}. \]
The proof of the following proposition is completely analogous to Proposition \ref{ordEqnType0Prop}, and is therefore omitted. 
\begin{proposition} \label{ordEqnType2Prop}
Let the notation be as in the previous paragraph. Then if $Z(\mathbf b) \cap \widehat  \Omega^{ord}_{\Lambda'} \neq \emptyset$, we have 
\[ Z(\mathbf b) \cap \widehat  \Omega^{ord}_{\Lambda'} \simeq Spf \  W[T, (T^p - T)^{-1}]^{\vee} / (f), \]
where 
\[f = \begin{cases} a_0  + a_1 T & \text{if } \mathbf b \in \V_{\phi}^+ \\
a_0' + a_1'T  & \text{if } \mathbf b \in \V_{\phi}^-. 
\end{cases}
\]
\qed
\end{proposition}

\begin{lemma} \label{rFormulaLemma} For $b \in C$, write
\[ ord_p(q(b)) = \begin{cases} 2t, & \text{if } ord_p(q(b)) \text{ is even} \\ 2t -1 , & \text{otherwise}.\end{cases} \]
For any vertex lattice $\Lambda$, let
\begin{equation} \label{rDefEqn}
 r(b, \Lambda) := max\{ r \in \Z \ | \ p^{-r} b \in\Lambda \};
\end{equation}
note that here we do not assume $b \in \Lambda$, and so $r(b, \Lambda)$ may be negative.
Finally, let $\Lambda_{\odot}$ denote the unique vertex lattice containing $p^{-t}b$, as in Lemma \ref{centralLatticeLemma}. Recall that $\Lambda_{\odot}$ is type 0 (resp. type 2) if $ord_p q(b)$ is even (resp. odd). 

Then we have the formula
\[ r(b, \Lambda) = \begin{cases}  t - \big\lfloor \frac{d(\Lambda, \Lambda_{\odot})}{2} \big\rfloor, & ord_pq(b) \text{ even} \\
                         t - \big\lfloor \frac{d(\Lambda, \Lambda_{\odot})+1}{2} \big\rfloor, & ord_pq(b) \text{ odd}.     
                         \end{cases} \]
Here  $d(\Lambda, \Lambda_{\odot})$ is the distance function on $\mathcal B$, the Bruhat-Tits tree. 

\begin{proof}
By scaling by a power of $p$, it suffices to prove this lemma in the case $t=0$; that is, we may assume that either $ord_p q(b) =  0 $ or $ord_p q(b) = -1$.  

We proceed by induction on $d = d(\Lambda, \Lambda_{\odot})$. 
If $d(\Lambda, \Lambda_{\odot}) = 0$, i.e. $\Lambda = \Lambda_{\odot}$, then $r(b, \Lambda_{\odot}) = 0$ by the definition of $ \Lambda_{\odot}$. 

Next, suppose we have proven the claim for all lattices $L$ with $0<d(L, \Lambda_{\odot}) \leq d$, and let $\Lambda$ be a lattice with $d(\Lambda, \Lambda_{\odot}) = d$. We shall prove the desired formula holds for all the neighbours of $\Lambda$. There are four cases here to check, as $\Lambda$ can be either type 0 or 2, and $ord_p q(b)$ can be even or odd.

For example, suppose that $\Lambda = \Lambda^{\sharp}$ and $ord_p q(b) = 0$. Then $d = d(\Lambda, \Lambda_{\odot})$ is even.

 There exists an $o_{k,p}$ basis $\{ v_0, v_1 \}$ for $\Lambda$ with $(v_0, v_0 ) = (v_1, v_1) = 0$, and $ (v_0, v_1) = -(v_1 , v_0) = \delta$. With respect to this basis, a complete list of the $p+1$ neighbours of $\Lambda$ in the Bruhat-Tits tree are described as follows:
\begin{align*} \Lambda'_{\infty} &:= span_{o_{k,p}} \{ p^{-1} v_0,  \ v_1 \}  \\
\Lambda'_{\alpha} &:= span_{o_{k,p}} \{ v_0, \ p^{-1} \alpha v_0 + p^{-1} v_1 \}
\end{align*}
as $\alpha \in \Z_p$ ranges over a complete set of representatives for $\F_p = \Z_p / p \Z_p$. 

Without loss of generality, we may assume that $d( \Lambda'_{\infty}, \Lambda_{\odot}) = d -1$, and so for all the other neighbours, $d(\Lambda'_{\alpha}, \Lambda_{\odot}) = d+1$. We may write
\begin{align*}
 b = p^{r}(a_0 v_0 + a_1 v_1) &= p^{r} \left( p a_0 \cdot ( p^{-1}v_0) +  a_1 \cdot v_1 \right) \\
 &= p^{r+1}\left(  a_0 \cdot ( p^{-1}v_0) +  p^{-1}a_1 \cdot v_1 \right)
 \end{align*}
where $r = r(b, \Lambda) = - \lfloor d/2 \rfloor$.

Now, noting that  $d$ is even, the induction hypothesis applied to $\Lambda $ and $\Lambda'_{\infty}$ yields  
\[ r(b, \Lambda'_{\infty}) = - \lfloor (d-1)/2 \rfloor = - \lfloor d/2 \rfloor + 1 = r(b, \Lambda) + 1 = r+1, \]
and so we must have $a_0 \in o_{k,p}^{\times}$ and $ p|a_1 $. Hence, by inspecting the remaining neighbours $\Lambda'_{\alpha}$ of $\Lambda$, we immediately see 
\[ r(b, \Lambda'_{\alpha}) = r = - \Big \lfloor \frac{d + 1}{2} \Big\rfloor, \]
as required.
The remaining cases all follow in the same manner. 
\end{proof}

\end{lemma}

We now determine the equations of the special cycles $Z(\mathbf b)$ in the local ring at a superspecial point $x \in \mathbb P_{\Lambda}(\F) \cap \mathbb P_{\Lambda'}(\F)$, for a pair of neighbouring vertex lattices $\Lambda$ and $\Lambda'$ of type 0 and 2 respectively.
Recall that we have a (formal) affine open neighbourhood of $x$ 
\[ \widehat\Omega_{[\Lambda, \Lambda']} \simeq Spf \left( W[ T_0, T_1, (T_0^{p-1} - 1)^{-1}, (T_1^{p-1} - 1)^{-1} ] / (T_0 T_1 - p) \right)^{\vee}, \]
where the point $x$ corresponds to the maximal ideal $\lie{m}_x = (T_0, T_1)$. 

Without loss of generality, we suppose that there is a basis $\{v_0, v_1\}$ for $\Lambda$, such that $(v_0, v_1) = - (v_1, v_0) = \delta$, $(v_0, v_0) = (v_1, v_1) = 0$, and
\begin{equation} \label{choiceOfLambda}
 \Lambda = span \{ v_0, v_1 \}, \qquad \Lambda' = span \{ p^{-1} v_0, v_1 \}. 
 \end{equation}

\begin{proposition} \label{ssEqnProp}
Suppose $\mathbf b \in \V_{\phi}^{\pm}$, with $b= \varphi^{\pm} \mathbf b \in C$ the corresponding vector. Suppose $x \in Z(\mathbf b)(\F)$ is a superspecial point, with $x \in \mathbb P_{\Lambda}(\F) \cap \mathbb P_{\Lambda'}(\F)$ as above. If we set
\[ r = r(b, \Lambda), \qquad r' = r(b, \Lambda'), \]
as in \eqref{rDefEqn}, then the equation for the cycle $Z(\mathbf b)$ in the local ring $\mathcal O_{\D, x}$ is given by
\[
\begin{cases} 	( T_0)^{r'} (T_1)^r = 0, & \text{if } \mathbf b \in \V_{\phi}^- \\
			(T_0)^{r'} (T_1)^{r+1} = 0 , & \text{if } \mathbf b \in \V_{\phi}^+.
\end{cases}
 \]

\begin{proof}

We begin with the case $\mathbf b \in \V_{\phi}^-$.
%
%
 With respect to a basis $\{ v_0, v_1 \}$ of $\Lambda$ as in \eqref{choiceOfLambda}, we may write
\begin{align}
\notag	b \ = \ & p^r(\alpha_0 v_0 + \alpha_1 v_1) \\
\label{bprimeEqn}		 \ = \ & p^r \left( p \alpha_0 \cdot (p^{-1} v_0 ) \ + \alpha_1 v_1 \right)
\end{align}
where, as usual, $r = r(b, \Lambda)$, and at least one of $\alpha_0$ or $\alpha_1$ is a unit in $o_{k,p}$. Then in light of the choice of basis \eqref{choiceOfLambda}, we have
\[ r' = \begin{cases} r, & \text{if } a_1 \in o_{k,p}^{\times} \\ r + 1, & \text{if } p \text{ divides } \alpha_1. \end{cases} \]
%

By the proof of Proposition \ref{ordEqnType0Prop}, Case (i), we have that the special cycle $Z(\mathbf b)$ is defined at $x$ by the vanishing of the element
\[f:=  p^r(\alpha_0 T_0 + \alpha_1) \in \mathcal{O}_{\D,x}. \]
Now consider the case $r' = r$ (i.e $\alpha_1 \in o_{k,p}^{\times}$). Then the factor $\alpha_0T_0 + \alpha_1$ is a unit in $\mathcal{O}_{\D,x}$ and so $Z(\mathbf b)$ is given by 
\[ p^r =( T_0) ^r \cdot (T_1)^r = 0, \]
as required. If on the other hand, $r' = r+1$, then $p$ divides $\alpha_1$ and $\alpha_0 \in o_{k,p}^{\times}$. Thus
\[ p^r(\alpha_0 T_0 + \alpha_1) = p^r T_0 \left( \alpha_0 +  (\alpha_1/p) T_1 \right) = p^r T_0 \cdot u \]
with $u = (\alpha_0 + (\alpha_1/p)T_1) \in \mathcal{O}_{\D,x}^{\times}$, and so the cycle $Z(\mathbf b)$ is given by 
\[ (T_0)^{r+1} (T_1)^{r} = (T_0)^{r'} (T_1)^r = 0,\]
as required.

The proof for a homomorphism $\mathbf b \in \V_{\phi}^+$ is completely analogous. 
\end{proof}
\end{proposition}

\begin{theorem} \label{cycleDecompThm}
Suppose $\mathbf b \in \V_{\phi}^{\pm}$, with $ord_p q^{\pm}(\mathbf b) \geq 0$ and write
\[ ord_p q^{\pm}(\mathbf b) = \begin{cases} 2t, & \text{ if }  ord_p q^{\pm}(\mathbf b) \text{ is even,  } \\ 2t - 1, & \text{ if }  ord_p q^{\pm}(\mathbf b) \text{ is odd.} \end{cases}\]
Let $b := \varphi^{\pm} \mathbf b \in C$ be the corresponding vector, and write $b = p^k b_{\odot}$, where $ord_p q(b_{\odot})$ is either $0$ or $-1$. Then by Lemma \ref{centralLatticeLemma}, there is a unique vertex lattice $\Lambda_{\odot}$ (the ``central lattice'') such that $b_{\odot} \in \Lambda_{\odot} - p \Lambda_{\odot}$. 
 
Finally, for any vertex lattice $\Lambda$ we define 
\begin{equation} \label{multDefEqn}
 m(\mathbf b, \Lambda) := \begin{cases} 0, & \text{if } b \notin \Lambda \\ t - \lfloor d(\Lambda, \Lambda_{\odot}) / 2 \rfloor, & \text{if } b \in \Lambda \text{ and } ord_p q^{\pm}(\mathbf b) \text{ is even, } \\
 t - \lfloor (d(\Lambda, \Lambda_{\odot})+1) / 2 \rfloor, & \text{if } b \in \Lambda \text{ and } ord_p q^{\pm}(\mathbf b) \text{ is odd. }  \end{cases}
\end{equation}
Then we have the equality of cycles on $\D$:
\[ 
Z(\mathbf b) = Z(\mathbf b)^{hor} \ + \ \sum_{\Lambda} \ m(\mathbf b, \Lambda) \ \mathbb P_{\Lambda},
\]
where $Z(\mathbf b)^{hor} \simeq Spf \ W$ is a horizontal divisor meeting the special fibre of $\D$ at a single ordinary special point in the component $\mathbb P_{\Lambda_{\odot}}$.

\begin{proof}
To start, we have the following relations:
\[ 
\begin{matrix}
m(\mathbf b, \Lambda) =	 \begin{cases} 	0 , 				& \text{if } b \notin \Lambda \\
								r(b, \Lambda )+ 1, 	& \text{if } b \in \Lambda, \ \Lambda^{\sharp} = \Lambda \\
								r(b, \Lambda), 		& \text{if } b \in \Lambda, \ \Lambda^{\sharp} = p \Lambda
					\end{cases}
&
\qquad \text{for } \mathbf b \in \V_{\phi}^+ , \ b = \varphi^+ \mathbf b \\
\end{matrix}
\]
and
\[
\begin{matrix}
m(\mathbf b, \Lambda) =	 \begin{cases} 	0 , 				& \text{if } b \notin \Lambda \\
								r(b, \Lambda), 		& \text{if } b \in \Lambda
					\end{cases}
&
\qquad \text{for } \mathbf b \in \V_{\phi}^-, 	\  b = \varphi^- \mathbf b
\end{matrix}
\]
which are easily verified by comparing the definition of $m(\mathbf b, \Lambda)$ above with Lemma \ref{rFormulaLemma}.

The proof of this theorem amounts to collating the information contained in the local equations given by Propositions \ref{ordEqnType0Prop}, \ref{ordEqnType2Prop}, and \ref{ssEqnProp}. We shall illustrate this argument in the case $\mathbf b \in \V_{\phi}^-$ with $ord_p q^-(\mathbf b) = 2t $ even; the remaining cases follow in an identical manner.

Suppose $b = \varphi^- \mathbf b \in C$ is the vector corresponding to $\mathbf b \in \V_{\phi}^-$, so by assumption $ord_p q (b) = ord_pq^-(\mathbf b) = 2t$ is also even. In particular, we have $b_{\odot} = p^{-t} b$. 

Let $\Lambda$ denote a type 0 lattice, and choose a basis $\{ v_0, v_1 \}$ such that $h(v_0, v_1) = \delta$ and $q(v_0) = q(v_1) = 0$. If we write
\[ b = p^r (\alpha_0 v_0 + \alpha_1 v_1) = p^r b_0, \]
where $r = r(b, \Lambda)$ and $b_0 \in \Lambda - p \Lambda$ as usual, then by Proposition \ref{ordEqnType0Prop}, we have that the intersection $Z(\mathbf b) \cap \widehat \Omega_{\Lambda}^{ord}$ with the \emph{ordinary locus} $\widehat\Omega_{\Lambda}^{ord}$ of $\widehat\Omega_{\Lambda}$ is defined by the vanishing of the element
\begin{equation} \label{cycleLambdaOrdEqn}
 p^r (\alpha_0 T + \alpha_1) \ \in \ W[T, (T^p - T)^{-1} ]^{\vee}.
\end{equation}
Note that by Lemma \ref{rFormulaLemma}, we have the following equivalencies:
\begin{align*}
ord_p q(b_0) = 0 \ \ \iff \ \ b_0 = b_{\odot} \ \ \iff \ \  r = r(b, \Lambda) = t  \ \ \iff \ \  \Lambda = \Lambda_{\odot}
\end{align*}
By the choice of basis $\{v_0, v_1\}$, we also have
\[ q(b_0) = \delta \left( \alpha_0 {\alpha_1'}  - \alpha_0' \alpha_1 \right). \]
If $\Lambda \neq \Lambda_{\odot}$, so that $ord_p q(b_0) > 0$, then $\alpha_0 \alpha_1^p \equiv \alpha_0^p \alpha_1 \mod p$. A short calculation implies that the term $(\alpha_0 T + \alpha_1 )$ is a unit in $W[T, (T^p - T)^{-1} ]^{\vee}$. Hence, for $\Lambda \neq \Lambda_{\odot}$, we have that $Z(\mathbf b) \cap \widehat \Omega_{\Lambda}^{ord}$ is determined by the equation
\[ p^r = 0. \]
As the component $\mathbb P_{\Lambda}^{ord}$ is given by the equation $p = 0$, we then obtain the equality of cycles
\begin{equation} \label{ZcycleLambdaEqn}
 Z(\mathbf b) \cap \widehat \Omega_{\Lambda}^{ord} \ = \ m(b, \Lambda)  \ \mathbb P_{\Lambda}^{ord}, \qquad \text{for } \Lambda \text{ of type } 0 , \ \Lambda \neq \Lambda_{\odot}. 
\end{equation}
On the other hand, if $\Lambda = \Lambda_{\odot}$, we have $\alpha_0, \alpha_1 \in o_{k,p}^{\times}$, and so the cycle defined by \eqref{cycleLambdaOrdEqn} is
\begin{equation} \label{ZcycleHorEqn}
 Z(\mathbf b) \cap \widehat \Omega_{\Lambda_{\odot}}^{ord} \ =\  Z(\mathbf b)^{hor} \ + \  m(\mathbf b, \Lambda_{\odot} )\ \mathbb P_{\Lambda_{\odot}}^{ord} \end{equation}
where 
\[ Z(\mathbf b)^{hor} := Spf \ W[T, (T^p - T)^{-1}]^{\vee} / (\alpha_0 T + \alpha_1) \simeq Spf \ W. \]

Now suppose $\Lambda'$ is a type 2 lattice. We choose a basis $\{ w_0, w_1 \}$ of $\Lambda'$ such that $h(w_0, w_1) = p^{-1} \delta$ and $q(w_0) = q(w_1) = 0$. Writing
\[ b = p^{r'} \left( \alpha_0 w_0 + \alpha_1 w_1 \right), \]
with $r' = r(b, \Lambda')$ and $b_0 := \alpha_0 w_0 + \alpha_1 w_1$, 
we have
\[ q(b_0) = p^{-1} \delta \left( \alpha_0 \alpha_1' - \alpha_0' \alpha_1 \right). \]
By Proposition \ref{ordEqnType2Prop}, the cycle $Z(\mathbf b) \cap \widehat\Omega^{ord}_{\Lambda'}$ is given by the vanishing of the element
\begin{equation} \label{cycleLambda'OrdEqn}
  p^{r'} (\alpha_0' + \alpha_1'T) \ \in \ W[T, (T^p - T)^{-1} ]^{\vee}.
\end{equation}
Since we have assumed at the outset that $ord_p q(b)$ is even, it follows that $ord_p q(b_0) > -1$. As before, this implies that the factor $(\alpha_0' +\alpha_1'T)$ is a unit, and so the cycle $Z(\mathbf b) \cap \widehat\Omega^{ord}_{\Lambda'}$ is given by the equation $p^{r'} = 0$. 
Hence, we have the equality of cycles
\begin{equation} \label{ZcycleLambda'Eqn}
  Z(\mathbf b) \cap \widehat\Omega^{ord}_{\Lambda'} \ = \ r' \ \mathbb P_{\Lambda'}^{ord}  \ = \ m(\mathbf b, \Lambda') \ \mathbb P_{\Lambda'}^{ord}, \qquad \text{ for } \Lambda' \text{ type } 2. 
\end{equation}

Combining \eqref{ZcycleLambdaEqn}, \eqref{ZcycleHorEqn} and  \eqref{ZcycleLambda'Eqn}, we have
\[ Z(\mathbf b)^{ord} = Z(\mathbf b)^{hor} \ + \ \sum_{\Lambda} \ m(\mathbf b, \Lambda) \ \mathbb P_{\Lambda}^{ord}, \]
where $Z(\mathbf b)^{ord}$ denotes the restriction of $Z(\mathbf b)$ to the ordinary locus of $\D$ (i.e. the open formal subscheme formed by the complement of the superspecial points).  

Now suppose $x$ is a superspecial point lying in the intersection $\mathbb P_{\Lambda} \cap \mathbb P_{\Lambda'}$, for a type 0 lattice $\Lambda$ and its type 2 neighbour $\Lambda'$. Then Proposition \ref{ssEqnProp} tells us that in a neighbourhood of $x$, the special cycle $Z(\mathbf b)$ is determined by the equation 
\[ (T_0)^{r'} \cdot( T_1)^{r}= 0. \]
Recall that in such a neighbourhood, the components $\mathbb P_{\Lambda}$ and $\mathbb P_{\Lambda'}$ are given by the equations $T_1 = 0$ and $T_0 = 0$ respectively, and so meet $Z(\mathbf b)$ with multiplicities $r = m(\mathbf b, \Lambda)$ and $r' = m(\mathbf b, \Lambda')$ respectively. Therefore, we have 
\[  Z(\mathbf b)= Z(\mathbf b)^{hor} \ + \ \sum_{\Lambda} \ m(\mathbf b, \Lambda) \ \mathbb P_{\Lambda}
\]
as required. 

We have proven the theorem in the case of an anti-linear special homomorphism $\mathbf b \in \V_{\phi}^-$ with $ord_p q^- (\mathbf b)$ even; the other cases all follow in a completely analogous manner.

\end{proof}
\end{theorem}

\end{subsection}

\begin{subsection}{Relation to local orthogonal special cycles}

In this section, we briefly recall definition of the \emph{local orthogonal special cycles} constructed by Kudla and Rapoport \cite{KRpadic}, and relate them to the unitary special cycles of the previous section.

\begin{definition} \label{locOrthCycDef} Suppose 
\[ j \in \End_{\OBp}(\X)_{\Q_p} := \End_{\OBp}(\X)\otimes_{\Z_p} \Q_p \]
is an $\OBp$-linear quasi-endomorphism of $\X$. Define $Z^o(j)^{\sharp}$ to be the closed formal subscheme of $\D$ which represents the following moduli problem: if $S \in \mathbf{Nilp}$, we take $Z^o(j)^{\sharp}(S)$ to be the locus of points $(X, \iota_X, \rho_X) \in \D(S)$ such that the quasi-morphism
\[ \rho_X^{-1} \circ j \circ \rho_X: \ X \times \overline S \to X \times \overline S \]
lifts to an endomorphism of $X$. 

We also define $Z^o(j)$ to be the Cohen-Macauleyfication of $Z^o(j)^{\sharp}$, namely the closed subscheme of $Z^o(j)^{\sharp}$ defined by the ideal sheaf of sections with finite support, cf. \cite[\S 4]{KRpadic}.
\end{definition}

These cycles (and their global counterparts, discussed in the next section) have been studied extensively by Kudla-Rapoport \cite{KRpadic} and Kudla-Rapoport-Yang \cite{KRYbook}; of immediate interest to us is the following special case of a result of Kudla and Rapoport, which is the counterpart to Theorem \ref{cycleDecompThm}.

\begin{proposition}[\cite{KRpadic} Proposition 4.5] \label{KROrthProp}
Suppose $j \in \End_{\OBp}(\X)_{ \Q_p}$ such that $j^2 = u^2 p^{2\alpha}\delta^2$ for some $u \in \Z_p^{\times}$ (so in particular, $Tr(j) = 0$, and $\Q_p(j)$ is isomorphic to $k_p$). Then  $Z^o(j)$ is a divisor on $\D$. 

Moreover, let  $\Lambda_0$ be the unique vertex lattice such that $j(\Lambda_0) = p^{\alpha} \Lambda_0$. For any vertex lattice $\Lambda$, define 
\[ \mu^o(j, \Lambda)  \ := \ \max \left\{ \ \alpha - d(\Lambda, \Lambda_0), \ 0  \ \right\}.  \]
Then there is an equality of cycles on $\D$:
\[ Z^o(j)= Z^o(j)^{hor} + \sum_{\Lambda} \mu^o(j, \Lambda) \mathbb P_{\Lambda}, \]
where $Z^o(j)^{hor}$ is a disjoint sum of two divisors, each of which is isomorphic to $Spf W$ and meets the special fibre of $\D$ at an ordinary special point in the component $\mathbb P_{\Lambda_0}$. 
\qed
\end{proposition}

Suppose $ \gamma \in End_{\OBp}(\X)_{\Q_p}$. Then $\gamma$ induces an endomorphism of $ N (\X)_0$, the 0-graded component of the rational Dieudonn\'e module of $\X$, and commutes with the operators $F$ and $V$. Therefore $\gamma$ defines an  $o_{k,p}$-linear endomorphism on $C := N(\X)_0^{V^{-1}F}$, which, abusing notation, we also denote by $\gamma$. 

Now we observe that if $j$ is as in Proposition \ref{KROrthProp}, then as an endomorphism of $C$, its characteristic polynomial splits. In particular, we may find an eigenvector $b_0 \in C $ with eigenvalue $p^{\alpha} \delta$. By scaling, we may assume $ord_pq(b_0) $ is either $0$ or $-1$. 
Note that $j$ also commutes with the Galois-semilinear operator $\epsilon := V^{-1} \circ \Pi_{\X}$, and so 
\[( j \circ \epsilon)(b_0) = (\epsilon \circ j) (b_0) = - p^{\alpha}  \delta \cdot\epsilon (b_0). \]
In other words, $b_0$ and $\epsilon (b_0)$ are eigenvectors for $j$, with eigenvalues $p^{\alpha} \delta$ and $- p^{\alpha} \delta$ respectively. Moreover, the ``central'' lattice $\Lambda_0$ attached to $j$, as in Proposition \ref{KROrthProp} (ii), is $\Lambda_0 = \text{span} \{ b_0, \epsilon(b_0) \}$.

In anticipation of the next section, we make the following definition: for $j$ as above, set
\begin{equation} \label{frobLocDefEqn}
\nu_p(j) = \nu_p(j, \phi) := \begin{cases} 1, & \text{ if } \exists \text{ an eigenvector } b\in C \text{ with } ord_pq(b) \text{ odd} \\ p, & \text{otherwise}. \end{cases}
\end{equation}
In other words, we have $\nu_p(j) = 1$ if and only if there exists a homomorphism $\mathbf b \in \V_{\phi}^+$ such that 
\[ j \circ \mathbf b = p^{\alpha} \cdot \left( \mathbf b \ \circ \ i_{\Y}(\delta) \right), \qquad \text{ with } \ ord_pq^+(\mathbf b) \text{ even.} \]
The following two theorems relate the orthogonal and unitary special cycles on $\D$. 
\begin{theorem} \label{locCycCompThm}
Suppose $j \in End_{\OBp}(\X)_{\Q_p}$, with $j^2 = u^2 p^{2\alpha} \delta^{2}$ for some $\alpha >0$ and $u \in \Z_p^{\times}$. Abbreviate $\nu = \nu_p(j)$, and fix $b_0 \in C$  an eigenvector with eigenvalue $p^{\alpha} \delta$ such that $q(b_0) = p^{-1} \nu$.
\begin{enumerate}[(i)]
\item If $\alpha$ is even, define $\mathbf b^+ \in \V_{\phi}^+$ and $\mathbf b^- \in \V_{\phi}^-$ by the relations:
\[ \varphi^+ (\mathbf b^+) = \nu^{-1} \cdot p^{\alpha/2}  \cdot b_0, \qquad \varphi^- (\mathbf b^- )= p^{\alpha/2}\cdot b_0. \]
\item If $\alpha$ is odd,  define $ \mathbf b^+$ and $\mathbf b^-$ by the relations:
\[\varphi^+(\mathbf b^+ ) = p^{\frac{\alpha-1}{2}} \cdot b_0, \qquad \varphi^-(\mathbf b^-) =\nu^{-1} \cdot  p^{\frac{\alpha +1}{2} }\cdot b_0\]
\end{enumerate}
Then we have an equality of cycles on $\D$
\begin{equation} \label{locCycCompEqn}
Z^o(j)  \ = \ Z(\mathbf b^+)\ + \ Z(\mathbf b^-). 
\end{equation}

\begin{remark} The key feature of this formula is that in every case, exactly one of the special homomorphism $\mathbf b^{\pm}$ appearing on the right hand side has norm $p^{\alpha}$, while the other has norm $p^{\alpha-1}$ (up to units in $\Z_p^{\times}$). $\diamond$
\end{remark}

\begin{proof}
Note that the central lattices for the cycles $Z^o(j)$, $Z(\mathbf b^+)$, and $Z(\mathbf b^-)$ are all the same lattice, namely $\Lambda_{\odot} := span_{o_{k,p}} \{ b, \epsilon(b)\}$. One can easily verify that in every case in the statement of the theorem, we have
\begin{align*}
 m(\mathbf b^+, \Lambda) + m(\mathbf b^-, \Lambda) &= \alpha - \lfloor d(\Lambda, \Lambda_{\odot}) / 2 \rfloor - \lfloor (d(\Lambda, \Lambda_{\odot}) + 1)/2\rfloor \\
 &= \alpha - d(\Lambda, \Lambda_{\odot}) \\
 &= \mu^o(j, \Lambda)
 \end{align*}
 for any vertex lattice $\Lambda$ with $d(\Lambda, \Lambda_{\odot}) \leq \alpha$. Indeed, the vectors $\mathbf b^+$ and $\mathbf b^-$ were scaled precisely so that this relation holds. On the other hand, if $d(\Lambda, \Lambda_{\odot}) > \alpha$, then neither side of \eqref{locCycCompEqn} meets the component $\mathbb P_{\Lambda}$.  Hence, any component $\mathbb P_{\Lambda}$ always occurs with the same multiplicity on both sides of \eqref{locCycCompEqn}.

It therefore suffices to show the horizontal components of both sides, which live in the open affine $\widehat\Omega_{\Lambda_{\odot}}^{ord}$, are equal. Suppose $\Lambda_{\odot}$ is type 0, i.e. 
\[ ord_p q(b_0) = 0 = ord_p \nu_p(j). \]
 Let $\{v_0, v_1 \}$ denote an $\epsilon$-invariant basis for $\Lambda_{\odot}$, with $q(v_0) = q(v_1) = 0$, and $h(v_0, v_1) = \delta$, and
write
\[ b_0 = a_0 v_0 + a_1 v_1, \qquad \epsilon(b_0) = a_0' v_0 + a_1' v_1. \]
Note that $q(b_0) = - q(\epsilon (b_0) )= \delta(a_0 a_1' - a_0' a_1)$ is a unit in $\Z_p^{\times}$, by assumption.

Let $j_0 := p^{-\alpha} j$. Then $(j_0)^2 = \delta^{2}$, and by \cite[Proposition 3.2]{KRpadic}, we have
\[ Z^o(j)^{hor} = Z^o(j_0). \]
As an endomorphism of $C$, the element $j_0$ is determined by the fact that $b_0$ and $\epsilon(b_0)$ are eigenvectors with eigenvalues $\delta$ and $-\delta$ respectively. It is therefore defined by the matrix
\[ [j] = \frac{\delta}{a_0a_1' - a_0'a_1} \begin{pmatrix} a_0a_1' + a_0' a_1 & -2 n(a_0) \\ 2 n(a_1) & -a_0a_1' - a_0'a_1 \end{pmatrix} \ \in \ M_2(\Z_p) \]
with respect to the basis $\{v_0, v_1 \}$.
Now in the affine neighbourhood 
\[ \widehat\Omega_{\Lambda_{\odot}}^{ord} \simeq Spf W[T, (T^p - T)^{-1}]^{\vee} ,\]
the proof of  Proposition \ref{ordEqnType0Prop} says that the cycle  $Z(\mathbf b^+)^{hor} + Z(\mathbf b^-)^{hor}$ is given by the vanishing of 
\begin{equation} \label{unHorEqn}
 (a_0 T + a_1) \cdot (a_0' T + a_1')  = n(a_0) T^2 + (a_0 a_1' + a_0'a_1)T + n(a_1).
\end{equation}
On the other hand, by \cite[Equation (3.5)]{KRpadic}, the equation for $Z^o(j_0)$ in $\widehat\Omega_{\Lambda_{\odot}}^{ord}$ is given by
\[ \frac{-2 \delta}{a_0a_1' - a_0'a_1} \left( n(a_0) T^2 + (a_0 a_1' + a_0' a_1) T + n(a_1) \right), \]
which differs from \eqref{unHorEqn} by a scalar in $\Z_p^{\times}$, and hence defines the same divisor. 

If $\Lambda_{\odot} = span\{ b_0, \epsilon(b_0) \}$ is a type 2 lattice, fix an $\epsilon$-invariant basis $\{ w_0, w_1 \}$ for $\Lambda_{\odot}$, such that $h(w_0, w_1) = p^{-1} \delta$ and $q(w_0) = q(w_1) = 0$. As before, we write
\[ b_0 = a_0 w_0 + a_1 w_1, \qquad \epsilon(b_0) = a_0' w_0 + a_1' w_1, \]
and note that $a_0 a_1' - a_0' a_1 \in o_{k,p}^{\times}$. 

The  cycle $Z(\mathbf b^+)^{hor} + Z(\mathbf b^-)^{hor}$ in $\widehat\Omega_{\Lambda_0}^{ord}$ is then given by the vanishing of the element
\[ (a_0 + a_1T) \cdot (a_0' + a_1' T) = n(a_0) + (a_0a_1' + a_0' a_1 )T + n(a_1)T^2. \]

Turning to the special endomorphism $j$, we note that it acts on $C$ via  the matrix
\[ 
 [j] = \frac{\delta}{a_0a_1' - a_0'a_1} \begin{pmatrix} a_0a_1' + a_0' a_1 & -2 n(a_0) \\ 2 n(a_1) & -a_0a_1' - a_0'a_1 \end{pmatrix} 
\]
with respect to the basis $\{ w_0, w_1 \}$. One can then check (e.g.\ by exploiting the action of $GL_2(\Q_p)$ on $\D$ cf.\ \cite[\S1]{KRpadic}), that the equation for $Z^o(j_0)$ in this case is given by 
\[   \frac{-2 \delta}{a_0a_1' - a_0'a_1} \left( n(a_0)  + (a_0 a_1' + a_0' a_1) T + n(a_1)T^2 \right), \]
and hence defines the divisor $Z(\mathbf b^+)^{hor} + Z(\mathbf b^-)^{hor}$, as required. 


\end{proof} 
\end{theorem}
We also have the corresponding theorem in the case $\alpha = 0$; note that all the cycles involved are horizontal. 
\begin{theorem}  \label{locCycCompHorThm}
Suppose $j \in End_{\OBp}(\X)_{\Q_p}$ with $j^2 = u^2 \delta^2 $ with $u \in \Z_p^{\times}$. Let $b_0 \in C$ denote an eigenvector with eigenvalue $\delta$, and suppose that $q(b_0) = p^{-1} \nu_p(j, \phi)$, so $ord_p q(b_0)$ is $0$ or $-1$.

Then if $\nu_p(j, \phi) = 1$, and we define $\mathbf b_1, \mathbf b_2 \in \V^+_{\phi}$ by the formulas
\[ \mathbf b_1 = \varphi^+(b_0), \qquad \mathbf b_2 = \varphi^+(\epsilon(b_0)), \]
then 
\begin{equation} \label{locCycCompHor+Eqn}
 Z^o(j)  = Z(\mathbf b_1) + Z(\mathbf b_2). 
\end{equation}
Similarly, if $\nu_p(j, \phi) = p$ and we define $\mathbf b_1, \mathbf b_2 \in \V^-_{\phi}$ by
\[ \mathbf b_1 = \varphi^-(b_0), \qquad \mathbf b_2 = \varphi^-(\epsilon(b_0)), \]
then
\begin{equation} \label{locCycCompHor-Eqn}
Z^o(j) = Z(\mathbf b_1) + Z(\mathbf b_2). 
\end{equation} 
\begin{proof}
The proof is completely analogous to the proof of Theorem \ref{locCycCompThm}, and is omitted. 
\end{proof}
\end{theorem}
\end{subsection}

\begin{subsection}{$p$-adic uniformizations} \label{padicUnifsSec}
Recall that we have fixed maps
\[ \tau_0\colon o_k / p \to \F \qquad \text{and} \qquad \tau_0 \colon o_{k,p} \to W;\]
 in this section, when we take the $\F$- or $W$-valued points of a stack defined over $o_k$, we shall implicitly do so via the $o_k$-structure induced by $\tau_0$.

If $\mathcal{X}$ denotes any of the stacks we have defined so far, we let $\widetilde{\mathcal{X}}$ denote the base change to $W$ of the formal completion of its fibre at $p$. The $p$-adic uniformizations we are about to describe relate these completions to the moduli spaces of $p$-divisible groups discussed in the previous section.

We begin with the Shimura curve $\CB$. Fix a pair $\underline{ \mathbf{A} }= (\mathbf A, \iota_{\mathbf A}) \in \CB(\F)$, such that the $p$-divisible group 
\[ \underline{\mathbf A}[p^{\infty}] = (\mathbf A[p^{\infty} ], \iota_{\mathbf A} \otimes \Z_p)\in  \D(\F), \]
together with the induced $\OBp$-action, 
is equal to the ``base-point'' $(\X, \iota_X)$ that we had fixed in defining the Drinfeld moduli space $\D$, cf.\ Definition \ref{drinfeldDef}. 
Let $B' = \End_{\OB}(\mathbf A) \otimes_{\Z} \Q$. Then $B'$ is a quaternion algebra over $\Q$ whose invariants differ from those of $B$ at exactly $p$ and $\infty$. We let $H' = (B')^{\times}$, viewed as an algebraic group over $\Q$, and define
\begin{align*} H'(\Q)^{\star} :=& \{ h' \in H'(\Q) \ | \ ord_p Nrd'(h') = 0 \} ,
\end{align*}
where $Nrd'$ is the reduced norm on $B'$. Let $K^p \subset H'(\A_f^p)$ denote the image of $(\widehat{\OB}{}^p)^{\times}$ under the natural identification $B'(\A_f^p) \simeq B^{op}(\A_f^p)$, and set 
\[ \Gamma' := K^p \cap H'(\Q)^{\star}.\]
 Then $\Gamma'$ acts on the Drinfeld upper-half plane $\D$ as follows: if $\underline X = (X, \iota_X, \rho_X) \in \D(S)$ for some $S \in \mathbf{Nilp}$, and $\gamma \in \Gamma'$, then we set
\[ \gamma \cdot \underline X \ := \ (X, \ \iota_X, \ \gamma_p \circ \rho_X), \]
where $\gamma_p$ denotes the image of $\gamma$ in $B' \otimes_{\Q} \Q_p = End_{\OBp}(\X) \otimes_{\Z_p} \Q_p$. We then obtain the following $p$-adic uniformization of $\CB$:

\begin{theorem} \label{ShCurveUnifThm}
There is an isomorphism of formal stacks over $\mathbf{Nilp}$:
\begin{equation}
 \widetilde{\CB} \ \simeq \ \left[ \Gamma' \backslash \D \right] .
 \end{equation}

\begin{proof} This is a restatement of \cite[Theorem III.5.3]{BC}, see also \cite[\S 8]{KRpadic}. Primarily to set up notation for the sequel, we indicate the idea of the proof. We begin by fixing  an $\OB$-linear isomorphism
\[ \eta_{\mathbf A}: Ta^p(\mathbf A) \isomto \widehat{\OB}{}^p, \]
where $Ta^p(\mathbf A) := \prod_{\ell \neq p} Ta_{\ell}(\mathbf A)$ is the prime-to-$p$ Tate module, and $\widehat{\OB}{}^p := \OB \otimes \widehat{\Z}{}^p$, which is viewed as an $\OB$-module via left-multiplication.

Let $(A, \iota) \in {\CB}(S)$, for some $S \in \mathbf{Nilp}$. Then there exists an $\OB$-linear quasi-isogeny
\begin{equation} \label{psiAEqn}
\psi_A: A  \to \mathbf A
\end{equation}
such that (i) the induced map on $p$-divisible groups (over $\overline S = S \times_W \F$ )
\[ \psi_A[p^{\infty}]_{\overline S}: \ A[p^{\infty}] \times \overline S  \ \to  \ \mathbf A[p^{\infty}] \times \overline S  = \X \times \overline S \]
is a quasi-isogeny of height 0, and (ii) the composition
\[ \eta_{\mathbf A} \circ \psi_A^p: Ta^p(A)_{\Q} \to B(\A_f^p) \]
maps $Ta^p(A)$ isomorphically onto $\widehat{\OB}{}^p$. 
Noting that any two quasi-isogenies satisfying both these properties necessarily differ by an element of $\Gamma'$, it follows that the map 
\[ \widetilde{\CB} \to [ \Gamma' \backslash \D], \qquad (A, \iota) \mapsto \Big[  A[p^{\infty}],  \ \iota \otimes \Z_p, \ \psi_A [p^{\infty}]_{\overline S} \Big] \]
is an isomorphism. 
\end{proof}
%
\end{theorem}

We have a similar statement for the orthogonal special cycles (Definition \ref{OrthCycleDef}).
\begin{theorem} \label{OrthUnifThm}
For $n \in \Z_{>0}$, let 
\begin{align} 
& \Omega^o(n) :=  \notag \\
&\ \ \left\{ \xi \in B'(\Q) \ | \ Trd'(\xi) = 0, \xi^2 = -n, \text{ and } \eta_{\mathbf A} \circ Ta^p(\xi) \circ \eta_{\mathbf A}^{-1} \in\End(\widehat{\OB}{}^p) \right\},\label{OmegaOrthDefEqn}
\end{align}
and note that $\Gamma'$ acts on $\Omega^o(n)$ by conjugation.
Then there is an isomorphism
\[ \widetilde{\Zed}{}^o(n) \simeq \left[ \Gamma' \ \Big\backslash \ \coprod_{\xi \in \Omega^o(n)} Z^o \Big(  \xi[p^{\infty}] \Big) \right] \]
as formal stacks over $W$. Here $\xi[p^{\infty} ] \in \End_{\OB}(\X)_{\Q_p}$ is the endomorphism of $\X = \mathbf A[p^{\infty}]$ induced by $\xi$, and $Z^o(\xi[p^{\infty}])$ is a local orthogonal special cycle as in Definition \ref{locOrthCycDef}.
\begin{proof}
This follows by rewriting \cite[(8.17)]{KRpadic} in terms of the uniformization as in Theorem \ref{ShCurveUnifThm}. 
\end{proof}
\end{theorem}

We now turn to the $p$-adic uniformization of the unitary special cycles, following \cite[\S 6]{KRunn}. 
Recall that we had fixed an embedding $\tau_0: o_k / (p) \to \F = \mathbb{F}_p^{alg}$, which lifts uniquely to an embedding $\tau_0: o_k \to W = W(\F)$; via these maps, we may view both $\F$ and $W$ as $o_k$-algebras.

Fix a triple $\underline{\mathbf E} = (\mathbf E, i_{\mathbf E}, \lambda_{\mathbf E}) \in \mathscr{E}^+(\F)$. We may also identify the $p$-divisible group $\underline{\mathbf E}[p^{\infty}]$, with its extra data, with the triple $\underline{\Y} = (\Y, i_{\Y}, \lambda_{\Y})$ of the previous section. 
We may further assume, by replacing $\mathbf E$ with an isogenous elliptic curve if necessary, that there is an $o_k$-linear isomorphism
\begin{equation} \label{etaEDefEqn}
 \eta_{\mathbf E}: Ta^p(\mathbf E) \isomto \widehat{o_k}{}^p
\end{equation}
such that the pullback of the symplectic form 
\[ (a,b) \mapsto tr(a b' (\sqrt{\Delta})^{-1} ) \]
on $\widehat{o_k}{}^p$ is equal to the Weil pairing
\[ e_{\mathbf E}\colon Ta^p(\mathbf E)_{\Q} \times Ta^p(\mathbf E)_{\Q} \ \to \ \A_f^p(1) \simeq \A_f^p. \]
defined by $\lambda_{\mathbf E}$. 
The two base points $\underline{\mathbf E} = (\mathbf E, i_{\mathbf E}, \lambda_{\mathbf E})$ and $\underline{\mathbf A} = (\mathbf A, \iota_{\mathbf A})$ allow us to define the global analogues of the spaces of special homomorphisms of the previous section, as follows. For any optimal embedding $\phi: o_k \hookrightarrow \OB,$
let 
\begin{align*}
&\mathcal{V}^+_{\phi} := \left\{ \beta \in \Hom(\mathbf E, \mathbf A)_{\Q} \ | \ \beta \circ i_{\mathbf E}(a) = \iota_{\mathbf A} \left( \phi(a) \right) \circ \beta \text{ for all }  a \in o_{k,p} \right\} \\
\text{and } \quad &\mathcal{V}^-_{\phi} := \left\{ \beta \in \Hom(\mathbf E, \mathbf A)_{\Q} \ | \ \beta \circ i_{\mathbf E}(a) = \iota_{\mathbf A} \left( \phi(a') \right) \circ \beta \text{ for all }  a \in o_{k,p} \right\}. 
\end{align*}
We view these spaces as $k$-vector spaces via the action 
\[ a \cdot \beta := \beta \circ i_{\mathbf E}(a), \qquad \text{for all } a \in k, \ \beta \in \mathcal{V}^{\pm}_{\phi}.  \]
Recall that for any optimal embedding $\phi$, we had defined a (non-principal) polarization $\lambda_{\mathbf A, \phi}$ on $\mathbf A$, as in \eqref{lambdaADefEqn}. The spaces $\mathcal{V}^{\pm}_{\phi}$ can then be equipped with Hermitian forms $h^{\pm}$, defined by the formulas:
\begin{align} \label{V+HermFormEqn}
h^+(\beta_1, \beta_2) :=& \lambda_{\mathbf E}^{-1} \circ \beta_2^{\vee} \circ \lambda_{\mathbf A, \phi} \circ \beta_1 \in \End_{o_k}(\mathbf E)_{\Q} \simeq k, \qquad \beta_1, \beta_2 \in \mathcal{V}^+_{\phi}
\intertext{and}  \label{V-HermFormEqn}
h^-(\beta_1, \beta_2) :=& \lambda_{\mathbf E}^{-1} \circ \beta_1^{\vee} \circ \lambda_{\mathbf A,\phi} \circ \beta_2 \in \End_{o_k}(\mathbf E)_{\Q} \simeq k, \qquad \beta_1, \beta_2 \in \mathcal{V}^-_{\phi}
\end{align}
respectively. Let $q^{\pm}(\beta) = h^{\pm}( \beta, \beta)$ denote the associated quadratic forms. 
Finally, we recall the decomposition $\mathscr{E} = \mathscr E^+ \coprod \mathscr E^-$, which induces a decomposition of the (global) cycles
\[ \Zed(m, \phi) \ = \ \Zed^+(m, \phi)\  \coprod \  \Zed^-(m, \phi). \]
We may now state their $p$-adic uniformization: 
\begin{theorem} \label{ZuLocalUnifThm}
 Suppose $m \in \Z_{>0}$. 
For any fractional ideal $\lie{a}$ in $k$, and any optimal embedding $\phi: o_k \to \OB$, let 
\begin{align}
 \Omega^{\pm}(m, \lie{a}, \phi) & := \notag \\
 &\Bigg\{ \beta \in \mathcal{V}^{\pm}_{\phi} \ | \ q^{\pm}(\beta) = m \cdot  \frac{p^{ord_pN(\lie{a})}}{ N(\lie{a})} , \ \  \eta_{\mathbf A} \circ \beta^p \circ \eta_{\mathbf E}^{-1} (\widehat{\lie{a}}{}^p) \subset \widehat{\OB}{}^p \Bigg\},
\end{align}
where, for any $\beta \in \mathcal{V}^{\pm}_{\phi}$, we denote by $\beta^p\colon Ta^p(\mathbf E)_{\Q} \to Ta^p(\mathbf A)_{\Q}$ the induced map on rational prime-to-$p$ Tate modules. Note $\Gamma'$ acts on these sets by composition, and we let $o_k^{\times}$ act via the $k$-vector space structure on $\mathcal{V}_{\phi}^{\pm}$. 
Then 
\begin{equation} \label{Z+unifEqn}
 \widetilde{\Zed}{}^{\pm}(m, \phi) \simeq \left[ ( \Gamma' \ \times  \ o_k^{\times}) \Big\backslash \ \left( \coprod_{[\lie{a}] \in Cl(k)} \coprod_{\beta \in \Omega^{\pm}(m, \lie{a}, \phi)} \ Z(\beta[p^{\infty}]) \right) \right], 
\end{equation}
where $\beta[p^{\infty}] \in \V^{\pm}_{\phi}$ is the quasi-morphism of $p$-divisible groups induced by $\beta$. Here $\lie a$ ranges any set of representatives of the class group $Cl(k)$ of $k$.
%
\begin{proof} This theorem is a reformulation of \cite[Proposition 6.3]{KRunn} in the present context.
\end{proof}
\end{theorem}
\begin{remark} Note  that the right hand side of \eqref{Z+unifEqn} is independent of the choices of representatives $\lie{a} \in Cl(k)$. Indeed, for any $a \in k^{\times}$, we have a bijection
\begin{equation} \label{scaleRelEqn}
	\Omega^{\pm}(m, a\cdot \lie{a}, \phi) \ \isomto \ \Omega^{\pm}(m, \lie{a} , \phi), \qquad \beta \mapsto  p^{-ord_pn(a)/2}a \ \cdot \ \beta.
\end{equation}
However, the cycles $Z(\beta[p^{\infty}])$ on $\D$ depend only on (i) the image of $\beta[p^{\infty}]$ in $\mathbb P(\V^{\pm}_{\phi_0})$ (which determines the horizontal part), and (ii) the p-adic valuation of $q^{\pm}(\beta)$, (which determines the vertical part).  By construction, both of these are the same for elements $\beta$ appearing on either side of \eqref{scaleRelEqn}. 
\end{remark}

\end{subsection}
\end{section}

\begin{section}{Proof of the main theorem}

\begin{subsection}{The (formal) Shimura lift } \label{ShimLiftSect}

The Shimura lift is a classical operation that takes modular forms of half-integral weight to modular forms of even weight. In its original formulation \cite{Shim}, the lift of a modular form $F$ is realized by forming a formal $q$-expansion involving the Fourier coefficients of $F$, and then proving that the result is indeed the $q$-expansion of a modular form with the desired properties.
This recipe inspires the following definition:
\begin{definition}[The formal Shimura lift]  \label{FShimLiftDef}
Let $M$ be any $\Z$-module, and 
\[ F \ = \ \sum_{n \geq 0} \ a(n) \ q^n \ \in \ M \llbracket q \rrbracket \]
a formal power series with coefficients in $M$. Suppose that
\begin{compactenum}[(i)]
\item $\kappa\geq 3$ is an odd integer and $\lambda := (\kappa - 1)/2$;
\item $N$ and $t$ are positive integers with $t$ squarefree;
\item and $\chi$ is a Dirichlet character modulo $4N$.
\end{compactenum} 
We define a new character
\[ \chi_t(a) := \chi(a) \ \left( \frac{ -1}{a} \right)^{(\kappa - 1)/2}  \left( \frac{t}{a} \right), \]
where $(\cdot / \cdot)$ denotes Shimura's modification of the Kronecker symbol, cf.\ \cite[Appendix A]{Cip}. 
Then the Shimura lift  of $F$, with respect to the parameters $(\kappa, N, t, \chi)$, is by definition the formal power series 
\begin{equation}  \label{ShLiftExpEqn}
 Sh(F) \ := \ \sum_{m \geq 0 } \ b(m) \ q^{tm} \ \in \ M  \llbracket q \rrbracket \otimes \C
\end{equation}
whose coefficients are given as follows. For $m>0$, we set 
\begin{equation} \label{shLiftFourierEqn}
 b(m) := \sum_{n | m} \ \chi_t(n) \ n^{(\kappa - 3)/2} \ a \left( t \frac{m^2}{n^2} \right).
\end{equation}
The constant term $b(0)$ is given as follows: let 
\[ \check{\chi_t}(a) = \sum_{h = 0}^{4Nt - 1} \ \chi_t(h) \exp( 2 \pi i a h /4Nt) \]
denote the Gauss sum, and set 
\begin{equation} \label{holShLiftCoeffZero} 
b(0) = a(0) \ \frac{(-1)^{\lambda}}{4} (2 \pi i)^{\lambda}  (\pi)^{-2 \lambda} (Nt)^{ \lambda - 1} \ \Gamma(\lambda) \ \sum_{m>0} m^{-\lambda} \check{\chi}_t(m) .
\end{equation}
Note that if $\chi^2 = 1$ then $b(m) \in M$ for all $m > 0$. 
\end{definition}
As alluded to above, if $F$ is the $q$-expansion of a holomorphic modular form of weight $\kappa/2$, level $\Gamma_0(4N)$ and nebencharacter $\chi$, then for each squarefree integer $t$, the Shimura lift $Sh(F)$   is the $q$-expansion of a modular form of weight $\kappa - 1$ for $\Gamma_0(4Nt)$ with character $\chi^2$, cf.\ \cite[Proposition 2.17]{Cip}.

\end{subsection}




\begin{subsection}{The Shimura lift formula in fibres of bad reduction} \label{CalcSec}
In this section, which is the technical heart of the paper, we show that a relation in the spirit of \eqref{shLiftFourierEqn} holds for the orthogonal and unitary special cycles in a formal neighbourhood of the fibre at $p$ for $p|D_B$ by comparing their $p$-adic uniformizations. Let $\F = \overline{\F}_p$, and fix a base point $\underline{\mathbf A} \in \CB(\F)$ as in Section \ref{padicUnifsSec}.

Our first aim is to relate the sets $\Omega^o(n)$ and $\Omega^{\pm}(m, \lie{a}, \phi)$, in the case of interest: the squarefree part of $n$ is equal to $|\Delta|$.
 Suppose $\xi \in \Omega^o(|\Delta| a^2)$, for some $a \in \Z_{>0}$. Then, by assumption, the endomorphism 
\[ \eta_{\mathbf A}^{-1} \circ \xi^p \circ \eta_{\mathbf A} \in End_{\OB}(\widehat{\OB}{}^p) \]
is given by right-multplication by some finite adele $(x_{\ell})_{\ell \neq p} \in \widehat{\OB}{}^p$, such that $(x_{\ell})^2 = a^2 \Delta$. Therefore, for every $\ell \neq p$, we obtain an embedding
\[ \varphi_{\ell}: k_{\ell} = k \otimes_{\Q} \Q_{\ell} \to B_{\ell}, \qquad a \sqrt{\Delta} \mapsto x_{\ell}. \]
We define the \emph{conductor} $c = c(\xi)$ of $\xi$ to be the smallest (rational) positive integer such that for all $\ell \neq p$, we have
\[ \varphi_{\ell} \left( \Z_{\ell}[ c \sqrt{\Delta}] \right) \subset \OB_{,\ell}; \]
in other words, $c(\xi)$ is the smallest integer such that $\varphi_{\ell}$ maps the unique $o_k$-order of conductor $c(\xi)$ into $\OB_{,\ell}$, for all $\ell \neq p$. We note that since
\[ \varphi_{\ell}\left( \Z_{\ell}[ a \sqrt{\Delta} ] \right) =  \Z_{\ell} + \Z_{\ell} x_{\ell} \subset \OB_{,\ell}, \]
we have that $c(\xi)$ divides $a$, and by valuation considerations it is easy to see $(c(\xi), D_B) = 1$. 
In particular, we obtain a disjoint decomposition
\begin{equation} \label{OmegaoCondDecompEqn}
 \Omega^o(|\Delta| a^2) = \coprod_{\substack{c|a \\ (c, D_B) = 1}} \Omega^o(|\Delta|a^2,c), 
\end{equation}
where $\Omega^o(|\Delta|a^2,c)$ denotes the  subset of elements $\xi \in \Omega^o(|\Delta|a^2)$ with $c(\xi) = c$. Note also that if $\xi \in \Omega^o(|\Delta|a^2,c)$ and $t$ is any integer, then the conductor of $t \cdot \xi$ is again $c$,
and one checks easily that we have a bijection
\begin{equation} \label{OmegaoScalingEqn}
 \Omega^o(|\Delta| a^2, c) \ \isomto  \ \Omega^o(|\Delta| a^2 t^2, c), \qquad \xi \mapsto t \cdot \xi.
\end{equation}

Now suppose $\phi: o_k \to \OB$ is an optimal embedding, and $\ell \neq p$ is a prime dividing $D_B$; recall our standing assumption that such a prime is inert in $k$. Then, reducing modulo $\ell$, we obtain two maps
\[ \overline{\varphi}_{\ell},  \  \overline \phi_{\ell} : o_{k, \ell} / (\ell) \to \OB_{,\ell} / (\theta), \]
where $\theta \in \OB$ is a fixed element such that $\theta^2 = - D_B$. As both the source and target of the maps are  isomorphic to the field $\F_{\ell^2}$, the two maps $\overline{\phi}_{\ell}$ and $\overline{\varphi}_{\ell}$ are either equal, or they differ by the Frobenius automorphism on $\F_{\ell^2}$.
This observation allows to define the \emph{Frobenius type away from $p$} of $\xi$ as follows:
\begin{equation} \label{frobAwayPDefEqn}
 \nu^p(\xi, \phi) := \prod_{\substack{\ell | D_B \\ \ell \neq p }} \nu_{\ell}(\xi, \phi),
\qquad \text{where} \quad
\nu_{\ell}(\xi, \phi) := \begin{cases} 
						1,				&	\text{if } \overline \phi_{\ell} = \overline \varphi_{\ell} \\
						\ell, 				&	\text{otherwise} .
				\end{cases}
\end{equation}
Recall that we also had the notion of a Frobenius type 
\begin{equation} \label{frobAtPDefEqn}
 \nu_p(\xi, \phi) := \nu_p(\xi[p^{\infty}], \phi_0)
\end{equation}
 at $p$, which was defined in \eqref{frobLocDefEqn} in terms of the parity of the $p$-valuation of the norm of any eigenvector of the induced map $\xi[p^{\infty}]$ on Dieudonn\'e modules. 

For an optimal embedding $\phi$, and an element $\beta \in \mathcal{V}^+_{\phi}$, we may also define the notion of a conductor, as follows. Let
\[ (h_{\ell})_{\ell \neq p} :=  \eta_{\mathbf A} \circ \beta^p \circ \eta_{\mathbf E}^{-1}(1) \in B(\A_f^p), \]
and for each $\ell \neq p$, define an embedding $ \varphi_{\ell}' : k_{\ell} \to B_{\ell}$ determined by the formula
\[ \varphi_{\ell}'(\sqrt{\Delta}) = (h_{\ell})^{-1} \cdot \phi(\sqrt{\Delta}) \cdot h_{\ell}. \]
As before, we define the \emph{conductor} $c(\beta)$ of $\beta$ to be the smallest integer $c$ such that $\varphi_{\ell}' (o_c) \subset \OB_{,\ell}$ for all $\ell \neq p$, where $o_c = \Z[ c \sqrt{\Delta}]$ is the unique order of conductor $c$. Note that this quantity depends on the choice of embedding $\phi$.

\begin{lemma} 
\label{NormLemma}
Suppose $\beta \in \mathcal{V}^+_{\phi}$, and let $h = (h_{\ell})_{\ell \neq p} :=  \eta_{\mathbf A} \circ \beta^p \circ \eta_{\mathbf E}^{-1}(1) \in B(\A_f^p)$. 
Then 
\begin{compactenum}
\item $q^+(\beta) =  \Delta \ Nrd_{B(\A_f^p)} (h)  \in \A_f^p$. 
\item If $\beta \in \Omega^+(m, \lie{a}, \phi)$, then the conductor $c(\beta)$ divides $m$.
\end{compactenum}

\begin{proof}
(i) Recalling our fixed trivialization $\A_f^p(1) \simeq \A_f^p$, consider the diagram 
\[
\begin{CD} 
Ta^p(\mathbf E)_{\Q} \times Ta^p(\mathbf E)_{\Q} @> \beta^p \times \beta^p>> Ta^p(\mathbf A)_{\Q} \times Ta^p(\mathbf A)_{\Q}  @>e_{\mathbf A}>> \A_f^p(1) \simeq \A_f^p \\
@V\eta_{\mathbf E} \times \eta_{\mathbf E}VV @V\eta_{\mathbf A} \times \eta_{\mathbf A}VV @| \\
\A_{k,f}^p \times \A_{k,f}^p @>(x,y) \mapsto \big(\phi(x) h, \phi(y)h \big)>> B(\A_f^p) \times B(\A_f^p) @>(x, y)  \mapsto  Trd \big(x y^{\iota} \phi(\sqrt{\Delta}) \big)>> \A_f^p
\end{CD}
\]
where $e_{\mathbf A}$ is the Weil pairing on $\mathbf A$ defined by the polarization $\lambda_{\mathbf A, \phi}$. It is an immediate consequence  of the definitions that both squares commute. 
Let 
\[ e_{\mathbf E} : Ta^p(\mathbf E)_{\Q} \times Ta^p(\mathbf E)_{\Q} \to \A_f^p \] 
denote the Weil pairing on $\mathbf E$ defined by $\lambda_{\mathbf E}$. Then on the one hand, by taking adjoints, we have  
\[ e_{\mathbf A} \left( \beta^p(x), \beta^p(y) \right) = q^+(\beta) \cdot e_{\mathbf E}(x, y), \qquad x, y \in Ta^p(\mathbf E)_{\Q}, \]
while on the other hand, if we write $s = \eta_{\mathbf E}(x)$ and $t = \eta_{\mathbf E}(y)$, then the commutative diagram above tells us that
\begin{align*}
 e_{\mathbf A} \left(\beta^p(x), \beta^p(y) \right)  =& \ Trd\left[ \ \phi(s) h  \cdot \left( \phi(t) h\right)^{\iota} \cdot \phi(\sqrt{\Delta} )  \right] \\
 =&  \ Nrd(h) \cdot Trd \left( \phi( s \cdot t' \cdot \sqrt{\Delta}) \right) \\
 =& \ Nrd(h) \ \Delta \ e_{\mathbf E}(x, y) ,
 \end{align*}
 where the last line follows by the choice of $\eta_{\mathbf E}$, as in \eqref{etaEDefEqn}. This proves (i). 
 
(ii) Suppose $\beta \in \Omega^+(m, \lie{a}, \phi)$. It suffices to check that
\[ m \cdot (h_{\ell})^{-1} \phi(\sqrt{\Delta}) h_{\ell} \ \overset{?}{\in} \ \OB_{,\ell} , \qquad \text{for all } \ell \neq p. \]
Choose $a_{\ell} \in k_{\ell}^{\times}$ such that $\lie{a} \otimes \Z_{\ell} = a_{\ell} \cdot o_{k,\ell}$. Then by the definition of $\Omega^+(m, \lie{a}, \phi)$, we have $h_{\ell} \in \phi(a_{\ell}^{-1}) \cdot \OB_{,\ell} $, and in particular, 
\[ (h_{\ell})^{\iota} \ \phi(\sqrt{\Delta}) \ h_{\ell} \ \in \frac{1}{N(\lie{a})} \OB_{,\ell}. \]
Hence, by combing part (i) of the lemma with the assumption
 \[ q^+(\beta) = m \cdot \frac{p^{ord_pN(\lie{a})}}{N(\lie{a})}, \]
 we have
\[ m \cdot (h_{\ell})^{-1} \phi(\sqrt{\Delta}) h_{\ell} =  p^{-ord_p(N(\lie{a}))} \cdot  \Delta \cdot N(\lie{a} ) \cdot \left[ (h_{\ell})^{\iota} \ \phi(\sqrt{\Delta}) \ h_{\ell} \right] \in \OB_{,\ell} \]
as required. 

\end{proof}
\end{lemma}

\begin{lemma} \label{eigenvectorLemma}
 Let $\beta \in \Omega^+(m, \lie{a}, \phi)$. Then there is a unique element $\xi = \xi(\beta) \in \Omega^o(|\Delta|m^2)$ such that 
$\xi \circ \beta = m \sqrt{\Delta} \cdot \beta$.
Moreover, the conductor $c(\xi)$ of $\xi$ is equal to $c(\beta)$. 
\begin{proof}
Note that each $\xi \in B'(\Q) $ defines a $k$-linear endomorphism $[\xi]$ of $\mathcal{V}^+_{\phi}$ by composition:
\[ [\xi] \cdot \beta := \xi \circ \beta. \]
Moreover, it follows immediately from definitions that
\[ q^+([\xi] \cdot \beta, \ [\xi] \cdot \beta) = Nrd'(\xi) \cdot q^+(\beta, \beta) \qquad \text{for all } \beta \in \mathcal{V}^+_{\phi}.\]
Hence, if $\beta \in \Omega^+(m, \lie{a}, \phi)$ and $\xi_1, \xi_2$ both satisfy $[\xi_i] \cdot \beta  = m \sqrt{\Delta} \cdot \beta$, then $Nrd'( \xi_1 - \xi_2) = 0$, so $\xi_1 = \xi_2$, as $B'$ is division. This proves the uniqueness claim in the lemma.

To show existence, we choose an element $\vartheta \in \OB$ with $\vartheta^{\iota} = - \vartheta$, and such that for all $a \in k$, we have $\vartheta \phi(a) = \phi(a') \vartheta$.  
Then we may define a quasi-isogeny
\[ \psi_{\beta}: \mathbf E \times \mathbf E \to \mathbf A, \qquad (x, y) \mapsto \beta(x) + \iota_{\mathbf A}(\vartheta) \cdot \beta(y).  \]
Note that for all $a \in o_k$, we have
\begin{align} \label{psibetaEqn1}
\iota_{\mathbf A}\left( \phi(a) \right) \cdot \psi_{\beta}(x, y) \ =&  \ \psi_{\beta} \left( i_{\mathbf E}(a) \cdot x, \ i_{\mathbf E}(a') \cdot y \right) \\
\text{and} \qquad \qquad  \label{psibetaEqn2}
\iota_{\mathbf A}\left( \vartheta \right) \cdot \psi_{\beta}(x, y)\ =&  \ \psi_{\beta}(\vartheta^2 y, x).
\end{align}
We may then define an element $\xi \in End(\mathbf A)_{\Q}$ by the formula
\[ \xi := m \cdot \psi_{\beta}  \circ \big( i_{\mathbf E}(\sqrt{\Delta}) , i_{\mathbf E}(\sqrt{\Delta}) \big) \ \circ \psi_{\beta}^{-1}. \]
It follows from \eqref{psibetaEqn1} and \eqref{psibetaEqn2} that $\xi$ commutes with $\iota_{\mathbf A}(\vartheta)$, and $\iota_{\mathbf A} (\phi(a))$ for all $a \in o_k$, so in fact $\xi \in B'(\Q) = End_{\OB}(\mathbf A)_{\Q}$, and it is easily seen that
\begin{equation} \label{evalueeqn}
 \xi \circ \beta = m \cdot \beta \circ i_{\mathbf E}(\sqrt{\Delta}) = m \sqrt{\Delta} \cdot \beta. 
\end{equation}
It is also straightforward to verify that $\xi$ satisfies the conditions \eqref{OmegaOrthDefEqn} defining $\Omega^o(m^2|\Delta|)$, as well as the claim regarding conductors.
%
%
%
\end{proof}
\end{lemma}

\begin{proposition} \label{Omega+FibreProp}
Fix a set of representatives $\{ \lie{a}_1, \dots, \lie{a}_h \}$ for the class group $Cl(k)$ of $k$, such that each $\lie{a}_i$ is relatively prime to $(p)$. 
Consider the map
\begin{align*}
f\colon \coprod_{i=1}^h \Omega^+(m, \lie{a}_i, \phi)   \to \Omega^o(|\Delta|m^2), \qquad \beta \mapsto \xi(\beta)
 \end{align*}
where $ \xi(\beta)$ is the unique element which has $\beta$ as an eigenvector with eigenvalue $m \sqrt{\Delta}$, as in Lemma \ref{eigenvectorLemma}.  For any $\xi \in  \Omega^o(|\Delta|m^2)$, let $\nu^p = \nu^p(\xi, \phi)$ (resp. $\nu_p = \nu_p(\xi[p^{\infty}], \phi)$) denote its  Frobenius type away from (resp. at) $p$, and $c = c(\xi)$ its conductor. Then we have
\[ \#  \left(  f^{-1}(\xi) \right ) = |o_k^{\times}| \cdot \rho \left( \frac{m}{c |\Delta| \nu_p \nu^p } \right), \]
where for any rational number $N$, $\rho(N)$ denotes the number of integral ideals of $k$ of norm $N$. 

In particular, if $m / c |\Delta| \nu_p \nu^p$ is not an integer, then the fibre $f^{-1}(\xi)$ is empty. 

\begin{proof}
Let $\xi \in \Omega^o(|\Delta|m^2)$. Then, viewing $\xi$ as a $k$-linear endomorphism of $\mathcal{V}^+_{\phi}$, it has two distinct eigenvalues $\pm m \sqrt{\Delta}$ and so there certainly exists some $\beta_0 \in \mathcal{V}^+_{\phi}$ such that 
\[ \xi \cdot \beta_0 = m \sqrt{\Delta} \cdot \beta_0. \]
Moreover, any other eigenvector with the same eigenvalue differs from $\beta_0$ by a scalar in $k^{\times}$. 
Hence, we have 
\begin{align}
 \#f^{-1}(\xi)  \	=& \ \# \coprod_{i=1}^h \left\{ \beta \in \Omega^+(m, \lie{a}_i, \phi) \ | \ \xi(\beta) = \xi \right\} \notag\\
			=& \ \# \coprod_{i=1}^h \left\{ a \in k^{\times} \ | \ q^+(a \beta_0) = \frac{m}{N(\lie{a}_i)}, \ \eta_{\mathbf A} \circ \beta_0^p \circ \eta_{\mathbf E}^{-1}(a \widehat{\lie{a}_i}{}^p) \subset \widehat{\OB}{}^p \right\} \notag \\
			=& |o_k^{\times} | \cdot \# \Big\{ \lie{a} \subset k  \text{ a fractional ideal } | \notag \\
			& \qquad \quad  \qquad N(\lie{a}) = \frac{m}{q^+(\beta_0)}, \quad  \eta_{\mathbf A} \circ \beta_0^p \circ \eta_{\mathbf E}^{-1}(\widehat{\lie{a}}{}^p) \subset \widehat{\OB}{}^p \Big\} \label{fibrePropIdealEqn}.
\end{align} 
As before, set
\[ h = (h_{\ell})_{\ell \neq p}  :=  \eta_{\mathbf A} \circ \beta_0^p \circ \eta_{\mathbf E}^{-1} (1), \qquad \text{and} \qquad x = (x_{\ell})_{\ell \neq p}  :=  \eta_{\mathbf A} \circ \xi^p \circ \eta_{\mathbf A}^{-1} (1) . \]
Then the condition $\xi \cdot \beta_0 = m \sqrt{\Delta} \cdot \beta_0$ implies
\begin{equation*}
x_{\ell} \ = \ m \cdot h_{\ell}^{-1} \ \phi(\sqrt{\Delta}) \ h_{\ell}, \qquad \text{for all } \ell \neq p.
\end{equation*}
For each $\ell \neq p$, let $\varphi_{\ell}: k_{\ell} \to B_{\ell}$ be the embedding determined by the relation
\begin{equation} \label{varphiDefEqn}
 \varphi_{\ell} (\sqrt{\Delta})\  = \ m^{-1} x_{\ell} \ = \ h_{\ell}^{-1} \ \phi(\sqrt{\Delta}) \ h_{\ell}. 
\end{equation}

We shall translate the conditions on the right hand side of \eqref{fibrePropIdealEqn} to a collection of local ones. First, suppose that $\ell$ is a prime not dividing $D_B$. We may fix an isomorphism $\OB_{,\ell} \simeq M_2(\Z_{\ell})$ that identifies 
\begin{equation} \label{choiceOfLocIsoEqn}
 \phi(\sqrt{\Delta}) \qquad \text{with} \qquad \begin{pmatrix} & 1 \\ \Delta & \end{pmatrix};
\end{equation}
note that this is possible when $\ell = 2$ on account of the assumption that $|\Delta|$ is even. 
Let $c = c(\xi)$ be the conductor of $\xi$, and define
\[ w(c)_{\ell}:= \begin{pmatrix} c & \\ & 1 \end{pmatrix} \in \OB_{,\ell}. \]
Then the map
\[ \varphi_c \colon k_{\ell} \to B_{\ell}, \qquad \varphi_c( \sqrt{\Delta}) = w(c)_{\ell}^{-1} \ \phi(\sqrt{\Delta}) \ w(c)_{\ell} \] 
is a local embedding of conductor $c$; that is, we have 
\[ \varphi_c ( k_{\ell}) \cap \OB_{,\ell} = \varphi_c \left( \ \Z_{\ell}[c\sqrt{\Delta}] \ \right). \] 
By the definition of the conductor of $\xi$, the same is true for the embedding $\varphi_{\ell}: k_{\ell} \to B_{\ell}$ as in \eqref{varphiDefEqn}. 
However, by \cite[Theoreme II.3.2]{Vig}, any two embeddings of conductor $c$ necessarily differ by an inner automorphism determined by an element of $(\OB_{,\ell})^{\times}$. Furthermore, for $x, y \in B_{\ell}^{\times}$, we have
\[ Ad_{x^{-1}} \circ \phi = Ad_{y^{-1}} \circ \phi \qquad \iff \qquad x = \phi(a) y \ \text{for some  } a \in k_{\ell}^{\times}.\] 
Hence, for each $\ell$ not dividing $D_B$, we may write
\begin{equation} \label{hellUnram}
h_{\ell} = \phi(a_{\ell}) \cdot w(c)_{\ell} \cdot u_{\ell} 
\end{equation}
for some $a_{\ell} \in k_{\ell}^{\times}$ and $u_{\ell} \in (\OB_{,\ell})^{\times}$. Note moreover that $h_{\ell} \in (\OB_{,\ell})^{\times}$ for almost all $\ell$.

Now consider a prime $\ell | D_B$, $\ell \neq p$. Fix a uniformizer $\Pi_{\ell} \in \OB_{,\ell}$ such that $\Pi_{\ell} \phi(a) = \phi(a') \Pi_{\ell}$ for all $a \in k_{\ell}$. Recalling our assumption that $\ell$ is inert in $k$, we may write
\begin{equation} \label{hellRam}
 h_{\ell} = \phi(a_{\ell} )\cdot  (\Pi_{\ell})^{\epsilon_{\ell}}, 
\end{equation}
for some $a_{\ell} \in k_{\ell}^{\times}$ and $\epsilon_{\ell} \in \{0,1\}$. One checks that for the reductions 
\[ \overline{\varphi}_{\ell}, \ \overline{\phi}  \ \in \ Hom \left( \ o_{k,\ell} / (\ell) \ , \ \OB_{,\ell} / (\Pi_{\ell}) \right),\]
 we have
\[ \overline{\varphi}_{\ell} = \overline{\phi} \ \iff \ \epsilon = 0, \]
where $\varphi_{\ell}$ is the embedding determined by $h_{\ell}$ as in \eqref{varphiDefEqn}. Hence, by definition of the Frobenius type \eqref{frobAwayPDefEqn}, we have
\[ \epsilon_{\ell} = ord_{\ell} \left( \nu_{\ell}( \xi, \phi) \right). \]

At this point, we have amassed a list of elements $(a_{\ell}) \in (\A_{k,f}^p)^{\times}$, as they appear in  \eqref{hellUnram} and \eqref{hellRam}.
We supplement this list with an element $a_p \in k_p^{\times}$ as follows: note that by the definition of the Frobenius type \emph{at p} (cf.\ \eqref{frobAtPDefEqn} and \eqref{frobLocDefEqn}), we have $ord_p \left( \nu_p(\xi, \phi) \right) \equiv ord_p q^+(\beta_0) \pmod{2} .$
We then set $ a_p := p^{r/2},$  where $ r = ord_p q^+(\beta_0) - ord_p \nu(\xi, \phi)$, and let $\lie{a}_0$ denote the fractional ideal defined by $(a_{\ell}) \in \A_{k,f}^{\times}$. 

Abbreviating $\nu_{\ell} = \nu_{\ell}(\xi, \phi)$, the product formula for $\Q$ implies
\begin{align}
\notag	q^+(\beta_0) 	=& \ 	\left(	\prod_{\ell \neq p} \ | q^+(\beta_0) |^{-1}_{\ell}\right)  |q^+(\beta_0)|^{-1}_p  \\
\notag 					=& \ \left( \prod_{\ell \neq p} |\Delta|^{-1}_{\ell} \cdot |Nrd(h_{\ell})|^{-1}_{\ell} \right) \ |q^+(\beta_0)|^{-1}_p 
												\qquad \qquad   [ \text{by Lemma } \ref{NormLemma} (i) ]\\
\notag					=& \	 \left( \prod_{\ell \nmid D_B} |\Delta|^{-1}_{\ell} |n(a_{\ell})|^{-1}_{\ell} \cdot |Nrd(w(c)_{\ell})|^{-1}_{\ell}   \right) \\
\notag					& \qquad \qquad \times \left( \prod_{\substack{ \ell | D_B \\ \ell \neq p}} |n(a_{\ell})|^{-1}_{\ell} \cdot \nu_{\ell} \right)
									 |n(a_p)|^{-1}_p \cdot \nu_p \\
\notag					=& \ |\Delta| \cdot  N(\lie{a}_0)  \left( \prod_{\ell \neq p} |Nrd(w(c)_{\ell})|^{-1}_{\ell} \cdot \nu_{\ell} \right)  \cdot (\nu_p) \\
\notag					=& \ |\Delta| \cdot N(\lie{a}_0) \cdot c\cdot \nu^p(\xi, \phi) \cdot \nu_p(\xi, \phi), 
\end{align}
where we have used the facts: (i) $Nrd(w(c)_{\ell}) = c$ for all $\ell$ and  (ii) $|\Delta|$ is relatively prime to $D_B$.

At long last, we return to the quantity we wish to compute. Note that for a fractional ideal $\lie a$ appearing in the right hand side of \eqref{fibrePropIdealEqn}, we have the equivalence
\[ \eta_{\mathbf A} \circ \beta^p \circ \eta_{\mathbf E}^{-1} (\widehat{\lie{a}}{}^p) \subset \widehat{\OB}{}^p \qquad \iff \qquad \phi(\widehat{\lie{a}}{}^p) \subset \widehat{\OB}{}^p \cdot h^{-1}. \]
This in turn is equivalent to the collection of local statements, for all $\ell \neq p$:
\[  \phi(\lie{a}_{\ell} ) \ \subset  \ {\OB}_{,\ell} \cdot h_{\ell}^{-1} \ = 
	\begin{cases} \OB_{,\ell} \cdot \  \Pi_{\ell}^{-\nu_{\ell}}   \phi(a_{\ell})^{-1}, &	 \text{if  } \ell | D_B \\ 
			\OB_{,\ell} \cdot w(c)_{\ell}^{-1} \cdot \phi(a_{\ell})^{-1}, & \text{if } \ell \nmid D_B. 
	\end{cases}
\]
Hence, replacing the ideals $\lie{a}$ appearing in \eqref{fibrePropIdealEqn} by $\lie{a}_0\cdot \lie{a}$, we obtain
\begin{align}
 \#  f^{-1} & (\xi) =    |o_k^{\times}| \cdot \# \Big\{  \lie{a} \subset k  \text{ a fractional ideal} \ | \ N(\lie{a}) = \frac{m}{c|\Delta| \nu^p \nu_p } , \notag \\
& \phi(\lie{a}_{\ell})  \subset \OB_{,\ell} \cdot w(c)_{\ell}^{-1} \text{ for } \ell \nmid D_B , \ \text{and }  \phi(\lie{a}_{\ell})  \subset \OB_{,\ell} \Pi_{\ell}^{-\nu_{\ell}} \text{ for } \ell | D_B  \Big\}. \label{fibrePropIdealEqn2}
\end{align}

To conclude the proof, we show that an ideal $\lie{a}$ appearing on the right hand side above is necessarily integral. Note that for $\ell | D_B$, including the case $\ell = p$, the condition on the norm of $\lie{a}$ implies $\phi(\lie{a}_{\ell}) \subset \OB_{,\ell}$. If $\ell \nmid D_B$, then with respect to the isomorphism $\OB_{,\ell} \simeq M_2(\Z_{\ell})$ as in \eqref{choiceOfLocIsoEqn}, we compute
\[ \OB_{,\ell} \cdot w(c)_{\ell}^{-1} = \left\{ \begin{pmatrix} c^{-1} x & y \\ c^{-1} w & z \end{pmatrix}, \ x, y, w, z \in \Z_{\ell} \right\}. \]
Recall that for $a ,b \in  \Q_{\ell}$, we have 
\[ \phi(a + b \sqrt{\Delta}) = \begin{pmatrix} a & b \\ b \Delta & a \end{pmatrix}. \]
Therefore,
\[ \Big( \phi( k_{\ell} ) \ \cap \ \OB_{,\ell} w(c)_{\ell}^{-1} \Big) \ \subset \ \phi \left( \Z_{\ell}[\sqrt{\Delta}] \right) \ \subset \ \OB_{,\ell}. \]
and so, for all fractional ideals $\lie{a}$ appearing in \eqref{fibrePropIdealEqn2}, and all finite primes $\ell$, we have $\phi(\lie{a}_{\ell}) \subset \OB_{,\ell}$. As $\phi:o_k \to \OB$ is an \emph{optimal} embedding, it follows that $\lie{a} \subset o_k$ is integral.

\end{proof}
\end{proposition}

Before stating the main result of this section, we need a few more lemmas.
\begin{lemma}  \label{varpiLemma}
There exists a prime $q$ that is split in $k$, and an element 
\[ \varpi \in \End(\mathbf E)_{\Q}, \qquad  - Nm(\varpi) = \varpi^2 = -pq, \]
such that $\varpi \circ i_{\mathbf E}(a) = i_{\mathbf E}(a') \circ \varpi $ for all $a \in o_k$. 
\begin{proof}
See \cite[p.\ 144]{Mann}. The idea is that since $\End(\mathbf E)_{\Q}$ is the quaternion algebra ramified at exactly $p$ and $\infty$, the existence of such an element $\varpi$ is equivalent to the existence of a prime $q$ such that 
\[ (-pq, \Delta)_{\ell} = \begin{cases} -1, & \text{if }\ell \in \{ \infty, p \} \\ 1, & \text{otherwise}, \end{cases} \]
where $(\cdot, \cdot)_{\ell}$ is the Hilbert symbol.
This imposes a finite set of congruence conditions on $q$, for which (infinitely many) solutions exists by Dirichlet's theorem, and furthermore such a solution is necessarily split in $k$.
\end{proof}
\end{lemma}

\begin{lemma} \label{conjNotEquivLemma} Suppose $\phi: o_k \to \OB$ is optimal, and let $\phi'$ denote the conjugate embedding. Then $\phi$ and $\phi'$ are not $\OBx$-equivalent.
\begin{proof}
Write $\phi' = Ad_t \circ \phi$ for some $t \in B^{\times}$, which is always possible by Noether-Skolem. Let $\ell$ be a prime dividing $D_B$; then there is a uniformizer $\Pi_{\ell}$ such that 
\[ \phi' = Ad_{\Pi_{\ell}} \circ \phi \in Hom(o_{k,{\ell}}, \OB_{,\ell}). \]
Hence $\Pi_{\ell}^{-1} \cdot t \in \phi(k_{\ell}^{\times})$, and so $ord_{\ell}Nrd(t)$ is necessarily odd. In particular, $t \notin \OBx$. 
\end{proof}
\end{lemma}

\begin{lemma} \label{AdtSwitchLemma} Let $m \in \Z_{>0}$ be a positive integer, $\phi_1: o_k \to \OB$ an optimal embedding, and $\lie{a}$ a fractional ideal of $k$. Suppose $t \in \OB$ such that $Nrd(t)$ divides $gcd(D_B ,m)$. Note that $t$ normalizes $\OB$, and in particular 
\[ \phi_2 := Ad_{t^{-1}} \circ \phi_1 \]
is again an optimal embedding. Then we have a bijection
\begin{align*} 
\Omega^{\pm}(m, \lie{a}, \phi_1) \ \isomto  \ \Omega^{\pm}( m/ Nrd(t), \lie{a}, \phi_2), \qquad \beta   \mapsto   \iota_{\mathbf A}\left(t^{-1}\right)  \circ \beta.
\end{align*}

\begin{proof} 
Suppose $\beta \in \Omega^{\pm}(m, \lie{a}, \phi_1)$, and set $\beta' := \iota_{\mathbf A}\left( t^{-1}\right) \circ \beta$. We first verify 
\[ \beta' \stackrel{?}{\in} \Omega^{\pm}(m/Nrd(t), \lie{a}, \phi_2). \] 
The condition on the norm of $\beta'$ follows immediately from definitions, and so we only need to check the inclusion
\begin{equation} \label{TModCondCheckEqn}
\eta_{\mathbf A} \circ (\beta')^p \circ \eta_{\mathbf E}^{-1} \left( \widehat{\lie{a}}{}^p  \right) \stackrel{?}{\subset} \widehat{\OB}{}^p 
\end{equation}
holds. Write $\widehat{\lie{a}}{}^p = (a_{\ell}) \cdot \widehat{o_k}{}^p$ for some prime-to-$p$ idele $(a_{\ell})$, and set 
\[ h = (h_{\ell})_{\ell \neq p} :=  \eta_{\mathbf A} \circ (\beta)^p \circ \eta_{\mathbf E}^{-1} \Big(( a_{\ell}) \Big) \in \widehat{\OB}{}^p, \]
and
\[ h' =  (h'_{\ell}) :=  \eta_{\mathbf A} \circ (\beta')^p \circ \eta_{\mathbf E}^{-1} \Big( (a_{\ell}) \Big); \]
note that to prove \eqref{TModCondCheckEqn}, it suffices to show that $h' \in \widehat{\OB}{}^p$. 

 Recall that we had chosen $\eta_{\mathbf A}: Ta^p(\mathbf A) \to \widehat{\OB}{}^p$ to be an $\OB$-linear isomorphism, where the action of $\OB$ on $\widehat{\OB}{}^p$ is given by left-multiplication. Hence we have
\begin{align*}
h' = t^{-1} \cdot h \ \in \ t^{-1} \cdot \widehat{\OB}{}^p.
\end{align*}

For primes $\ell$ not dividing $Nrd(t)$, note that  $t^{-1} \cdot \OB_{,\ell} = \OB_{,\ell}$. 
On the other hand, suppose that $\ell$ divides $Nrd(t)$; in particular, $\ell$ divides $gcd(D_B,m)$. 
Then, by Lemma \ref{NormLemma},  
\[ ord_{\ell} Nrd(h_{\ell}) = ord_{\ell} (m)  \geq 1. \]
As $B_{\ell}$ is division, and $ord_{\ell}Nrd(t) =1$, it follow that $h'_{\ell} = t^{-1} \cdot h_{\ell} \in \OB_{,\ell} $ by valuation considerations. 

Thus $h'_{\ell} \in \OB_{,\ell}$ for all $\ell \neq p$, as required, and so we have shown that the assignment $\beta \mapsto \iota_{\mathbf A}(t^{-1}) \circ \beta$ indeed defines a map $\Omega^{\pm}(m, \lie{a}, \phi_1) \to \Omega^{\pm}(m/Nrd(t), \lie{a}, \phi_2)$. By a similar argument, there is an inverse map
\[ \Omega^{\pm}(m/Nrd(t), \lie{a}, \phi_2) \to \Omega^{\pm}(m, \lie{a}, \phi_1), \qquad \beta' \mapsto \beta := \iota_{\mathbf A}(t) \circ \beta', \]
which concludes the proof of the lemma.
\end{proof}
\end{lemma}

Recall that $\widetilde{\Zed}{}^{\pm}(m, \phi)$ denotes the formal completion  of $\Zed^{\pm}(m, \phi)$ along its fibre at $p$, and that we have
 $p$-adic uniformizations 
\begin{equation} \label{padicUnifLaterEqn}
 \widetilde{\Zed}{}^{\pm}(m, \phi) \simeq \left[ o_k^{\times} \times \Gamma' \ \Big\backslash \coprod_{[\lie{a}]} \coprod_{\beta \in \Omega^{\pm}(m, \lie{a}, \phi)} Z(\beta[p^{\infty}])\right] 
\end{equation}
as in Theorem \ref{ZuLocalUnifThm}.

\begin{lemma} \label{uCycleTwistLemma}
Let $m \in \Z_{>0}$ be a positive integer, $\phi_1: o_k \to \OB$ an optimal embedding, and $\lie{a}$ a fractional ideal of $k$. Suppose  $t \in \OB$ such that $Nrd(t)$ divides $gcd(D_B/p, m)$, and let $\phi_2 := Ad_{t^{-1}} \circ \phi_1$. Then we have an equality of cycles on $\widetilde{\CB}$:
\[ \widetilde{\Zed}{}^{\pm}(m, \phi_1) = \widetilde{\Zed}{}^{\pm}(m/ Nrd(t), \phi_2). \]
In particular, the cycle $\widetilde{\Zed}{}^{\pm}(m, \phi)$ only depends on the equivalence class $[\phi] \in Opt / \OBx$ of $\phi$. 

\begin{proof}
Let $\beta \in \Omega^{\pm}(m, \lie{a}, \phi_1)$, set $\beta' = \iota_{\mathbf A}(t^{-1}) \circ \beta$, and let $\mathbf b = \beta[p^{\infty}]$ and $\mathbf b' = \beta'[p^{\infty}]$ denote the corresponding maps of $p$-divisible groups. Recall that for a scheme $S \in \mathbf{Nilp}$, with special fibre $\overline S = S \times \F$, we had defined the $S$-points of the cycle $Z(\mathbf b)$ on $\D$ to be the locus of points $(X, \iota_X, \rho_X) \in \D(S)$ such that the map 
\[ \rho_X^{-1} \circ \mathbf b: \Y_{\overline S} \to X_{\overline S} \]
lifts to a map $\Y_S \to X$. 

Similarly, the cycle $Z(\mathbf b')$ parametrizes tuples $(X, \iota_X, \rho_X)$ such that $ \rho_X^{-1} \circ \mathbf b': \Y_{\overline S} \to X_{\overline S}$ lifts to a map $\Y_S \to X$. Note that the image of $t$ in $B_p$ in fact lies in $(\OBp)^{\times}$ by valuation considerations, and the endomorphism 
\[ \iota_X(t^{-1})_{\overline{S}} : X \times \overline S \to X \times \overline S \]
evidently admits a lift to $S$, namely $\iota_X(t^{-1}) \in End_S(X)^{\times}$. 
Combining these observations with the fact that the quasi-isogeny $\rho_X$ is assumed to be $\OBp$-linear, we have that
\begin{align*}
\rho_X^{-1} \circ \mathbf b' = \rho_X^{-1} \circ \iota_{\X}(t^{-1})_{\overline S} \circ \mathbf b \ \text{ lifts} \ &\iff \ \iota_{X}(t^{-1})_{\overline S} \circ \rho_X^{-1} \circ \mathbf b \ \text{ lifts} \\
&\iff \rho_X^{-1} \circ \mathbf b \ \text{ lifts}.
\end{align*}
Hence, we find that \emph{as cycles on $\D$}, we have $Z(\mathbf b) = Z(\mathbf b').$
The lemma then follows from Lemma \ref{AdtSwitchLemma} and the $p$-adic uniformization \eqref{padicUnifLaterEqn}.
\end{proof}
\end{lemma}

We now state and prove our key result:
\begin{theorem}  \label{mainCompletionsThm}

Let $m \in \Z_{>0}$. Then we have the following equalities of cycles:

(i) If $|\Delta| $ does not divide $m$, then $   \widetilde{\Zed}(m, \phi) = 0$. 

(ii) If $m = |\Delta| m'$, then we have 
\begin{align} \label{FourierCoeffEqn}
 \sum_{[\phi] \in Opt/ \OBx}  & \widetilde{\Zed}(m, \phi) + \widetilde{\Zed}   \left(\frac{m}{ gcd(D_B,m)}, \phi \right) \notag \\
 &= 2 h(k) \sum_{\substack{\alpha | m' \\ (\alpha, D_B)=1 }}  \  \chi_k(\alpha)  \ \widetilde{\Zed}{}^o \left( |\Delta| \frac{(m')^2}{\alpha^2} \right), 
\end{align}
where $\chi_k$ is the quadratic character associated to $k$, cf.\ \eqref{chiKNotSecEqn}. 
\begin{proof}
(i) Suppose $|\Delta|$ does not divide $m$. Then for all possible values of $c$, $\nu_p$, and $\nu^p$, the quantity $m/c |\Delta| \nu_p \nu^p$
 is not an integer, and so by Proposition \ref{Omega+FibreProp}, we have
 \[ \Omega^+(m, \lie{a}_i, \phi) = \emptyset \]
 for any fractional ideal $\lie{a}_i$ and any embedding $\phi$. Hence we also have
 \[ \Omega^-(m, \lie{a}_i, \phi) = \Omega^+(m, \lie{a}_i, \phi') = \emptyset, \]
where $\phi'$ is the conjugate embedding. 
By the $p$-adic uniformization \eqref{padicUnifLaterEqn}, it follows that $\widetilde{\Zed}{}^{\pm}(m,\phi)=\emptyset,$ and so as a cycle, $\widetilde{\Zed}(m, \phi) = 0$. 

(ii) Suppose $|\Delta|$ divides $m$, and $\phi: o_k \to \OB$ is an embedding with conjugate $\phi'$. By Lemma \ref{conjNotEquivLemma}, the embeddings $\phi$ and $\phi'$ are $\OBx$-inequivalent. Moreover, for any $n$, we have $\Zed^+(n, \phi) = \Zed^-(n, \phi')$ essentially by definition. Summing over classes of optimal embeddings, it follows that 
\begin{align}
  \sum_{[\phi] \in Opt/ \OBx} \widetilde{\Zed}(n, \phi) &= 2 \sum_{[\phi] \in Opt / \OBx} \widetilde{\Zed}{}^+(n, \phi) \notag \\
&= 2 \sum_{[\phi] \in Opt / \OBx} \widetilde{\Zed}{}^-(n, \phi)  . \label{signSwitchEqn}
\end{align}
We proceed by cases:

\fbox{ \textbf{Case 1:} $ord_p(m) > 0$:}  
For any cycle $Z$ on $\D$, we let $[Z]$ denote the corresponding cycle on $\widetilde{\CB} = [ \Gamma' \backslash \D]$, and for convenience, we let
\[ A\  := \ \sum_{[\phi]}  \ \widetilde{\Zed} (m, \phi)  \ + \  \widetilde{\Zed}(m/gcd(m,d), \phi) \]
denote the cycle we wish to compute (the left hand side of \eqref{FourierCoeffEqn}). 

Write $gcd(m, D_B) = p \cdot \mu$, and fix an element $t \in \OB$ such that $Nrd(t) = \mu$. Suppose $\phi_1: o_k \to \OB$ is an optimal embedding. By Lemma \ref{uCycleTwistLemma}, we have a equality of cycles
\[ \widetilde{\Zed}{}^-(m / p, \phi_1) = \widetilde{\Zed}{}^- \left( \frac{m}{gcd(m, D_B)}, Ad_{t^{-1}} \circ \phi_1 \right), \]
and so, upon taking the sum over all optimal embeddings and using \eqref{signSwitchEqn}, 
\begin{align*}
A &=2 \sum_{[\phi] \in Opt/\OBx} \widetilde{\Zed}{}^+(m, \phi) + \widetilde{\Zed}{}^-  \left( \frac{m}{ gcd(m, D_B)}, \phi \right) \\
 &= 2 \ \sum_{[\phi]} \widetilde{\Zed}{}^+(m, \phi) + \widetilde{\Zed}{}^- (m /p, Ad_{t} \circ \phi)  \\
&=  2  \ \sum_{[\phi]} \widetilde{\Zed}{}^+(m, \phi) + \widetilde{\Zed}{}^- (m /p, \phi) . 
 \end{align*}
Now suppose $\beta \in \Omega^+(m, \lie{a}, \phi)$ and set 
\[ \beta' := \beta \circ \varpi, \]
where $\varpi \in End(\mathbf E)$ satisfies $\varpi^2 = -pq$ for some split prime $q$, and $\varpi \circ i_{\mathbf E}(a) = i_{\mathbf E}(a') \circ \varpi$ for all $a \in o_k$, cf. Lemma \ref{varpiLemma}. Then it is easily verified that 
\[ \beta' \in \Omega^-(m / p, \ \lie{a}' \cdot \lie{q}, \ \phi),\]
 where $\lie{a}'$ is the conjugate of $\lie{a}$, and $\lie{q}$ is one of the ideals above $q$. 
 
Furthermore, let $\xi = \xi(\beta) \in \Omega^o(m^2|\Delta|)$ denote the unique element such that $\beta$ is an eigenvector with eigenvalue $m \sqrt{\Delta}$, cf.\ Lemma \ref{varpiLemma}. Then $\beta'$ is also an eigenvector, with eigenvalue $-m \sqrt{\Delta}$. 

Moreover, the corresponding maps $\beta[p^{\infty}]$ and $\beta'[p^{\infty}]$ on $p$-divisible groups are (up to scalars in $\Z_p^{\times}$) \emph{precisely} the special homomorphisms $\mathbf b^+$ and $\mathbf b^-$ appearing in Theorem \ref{locCycCompThm}, with $\xi[p^{\infty}]$ playing the role of $j$. The conclusion of that theorem then reads
\[ Z^o(\xi[p^{\infty}])^{pure} \ =  \ Z(\beta[p^{\infty}]) \ + \ Z(\beta'[p^{\infty}]),\]
as cycles on $\D$. Hence
\begin{align}
A &= 2 \sum_{[\phi]} \sum_{[\lie{a}]} \frac{1}{|o_k^{\times}|} \Bigg\{ \sum_{\substack{\beta \in \Omega^+(m, \lie{a}, \phi) \\ \text{mod } \Gamma'}} [ Z(\beta[p^{\infty}])]  + \sum_{\substack{\beta' \in \Omega^-(m/p, \lie{a}, \phi) \\ \text{mod } \Gamma'}} [Z(\beta'[p^{\infty}])] \Bigg\} \notag \\
&= 2 \sum_{[\phi]} \sum_{[\lie{a}]} \frac{1}{|o_k^{\times}|} \Bigg\{\sum_{\substack{\beta \in \Omega^+(m, \lie{a}, \phi) \\ \text{mod } \Gamma'}} [ Z(\beta[p^{\infty}])]  + \sum_{\substack{\beta' \in \Omega^-(m/p, \lie{q}\lie{a'}, \phi) \\ \text{mod } \Gamma'}} [Z(\beta'[p^{\infty}])] \Bigg\}  \notag \\
&= 2 \sum_{[\phi]} \sum_{[\lie{a}]} \frac{1}{|o_k^{\times}|} \Bigg\{ \sum_{\substack{\beta \in \Omega^+(m, \lie{a}, \phi) \\ \text{mod } \Gamma'}} [Z^o(\xi[p^{\infty}])] \Bigg\}\  \ \ \text{where } \xi = \xi(\beta). \label{mainThmIntExp1}
\end{align}
Given integers $c$, $\nu_p$, and $\nu^p$ such that $c|m$, $\nu_p \in \{ 1, p \}$ and $\nu^p | (D_B / p)$,  we set 
\begin{align*}
\Omega^o(|\Delta|m^2,& c, \nu, \phi)  := \\
& \left\{ \xi \in \Omega^o(|\Delta|m^2) \ | \ c(\xi) = c, \ \nu_p(\xi, \phi) = \nu_p, \ \nu^p(\xi, \phi) = \nu^p  \right\};
\end{align*}
that is, the set of elements of $\Omega^o(|\Delta|m^2)$ whose conductor and Frobenius types relative to $\phi$ are as specified, and where, for ease of notation, we have written $\nu = \nu_p \cdot \nu^p$. Note that the action of $\Gamma'$ by conjugation preserves these sets. 
By Proposition \ref{Omega+FibreProp}, we may continue:
\begin{align}
A= 2 \ \sum_{[\phi]} \sum_{\substack{c | m \\ (c, D_B) = 1}} \sum_{\nu | D_B } \rho \left( \frac{m}{|\Delta| c \nu^p \nu_p} \right) \left\{ \sum_{\substack{ \xi \in \Omega^o(m^2 |\Delta|,c, \nu,  \phi) \\ \text{mod } \Gamma'}} [Z^o(\xi[p^{\infty}])]\right\}. \notag
\end{align}
Writing $m = |\Delta|m'$, only terms with $c $ dividing $m'$ contribute to the above sum. Morevover, recall that all  primes dividing $D_B$ are inert in $k$ by assumption. Let $\nu^*$ be the unique integer dividing $D_B$ such that $ord_{\ell}( m' /  \nu^*)$ is even for all $\ell$ dividing $D_B$: only terms involving $\nu = \nu^*$ can contribute to the above sum, since $o_k$  has no ideals with norm an odd power of an inert prime. 
Thus 
\begin{align*}
A =2   \sum_{\substack{c | m' \\ (c, D_B) = 1}}  \rho \left( \frac{m'}{ c \nu^*} \right) \left\{\sum_{[\phi]} \sum_{\substack{ \xi \in \Omega^o(m^2 |\Delta|,c, \nu^*,  \phi) \\ \text{mod } \Gamma'}} [Z^o(\xi[p^{\infty}])]\right\}.
\end{align*}
Note that by \cite[Corollaire III.5.12]{Vig}, 
\[ \# Opt / \OBx = h(k) \cdot 2^{o(D_B)}, \]
where $o(D_B)$ denotes the number of prime factors of $D_B$. Now consider  the action of the normalizer $N_{B^{\times}}(\OB)$ on the set $Opt/\OBx$, acting by conjugation. For a fixed element $\xi \in \Omega^o(m^2|\Delta|)$, the various values of the Frobenius types $\nu(\xi, \phi)$, as $\phi$ varies in an $N(\OB)$-orbit of optimal embeddings, will cover all $2^{o(D_B)}$ possibilities. Thus, for fixed $\xi$, there are exactly $h(k)$ classes $[\phi]$ of optimal embeddings such that
\[ \nu(\xi, \phi) = \nu^p(\xi, \phi) \cdot \nu_p(\xi, \phi) = \nu^*. \]
Hence, it follows that
\[ A =  2 h(k)  \sum_{\substack{c | m' \\ (c, D_B) = 1}}  \rho \left( \frac{m'}{ c \nu^*} \right)  \sum_{\substack{ \xi \in \Omega^o(m^2 |\Delta|) \\ c(\xi) = c \\ \text{mod } \Gamma'}} [Z^o(\xi[p^{\infty}])].
\]
Next, we claim that for any integer $N>0$, we may write
\[ \rho(N) = \sum_{a | N} \chi_k(a); \]
indeed, both sides of the above formula are multiplicative, and for $N = \ell^n$ a prime power, the fomula can immediately be verified by considering the cases $\ell$ split, inert, and ramified separately. Moreover, as $ord_{\ell}( m' / c \nu^*)$ is even for all $\ell | D_B$, it follows that 
\[ \rho \left( \frac{m'}{ c \nu^*} \right) = \sum_{\substack{a | (m' / c \nu^*)\\ (a, D_B) = 1}}  \chi_k(a) =  \sum_{\substack{a | (m' / c)\\ (a, D_B) = 1}}  \chi_k(a). \]
Substituting, we obtain
\begin{align}
A \ &= \ 2 h(k) \sum_{\substack{c | m' \\ (c, D_B) = 1}} \  \sum_{\substack{a | (m' / c)\\ (a, D_B) = 1}} \  \chi_k(a) \  \sum_{\substack{ \xi \in \Omega^o(m^2 |\Delta|) \\ c(\xi) = c \\ \text{mod } \Gamma'}} [Z^o(\xi[p^{\infty}])] \notag \\
&=2 h(k) \sum_{\substack{a | m' \\ (a, D_B) = 1}}  \chi_k(a)  \cdot \sum_{\substack{c | (m' / a) \\ (c, D_B) = 1}}    \sum_{\substack{ \xi \in \Omega^o(m^2 |\Delta|) \\ c(\xi) = c \\ \text{mod } \Gamma'}} [Z^o(\xi[p^{\infty}])]. \label{AFinalExp}
\end{align}
Applying \eqref{OmegaoScalingEqn}, we have that for each $a | m'$ with $(a, D_B) = 1$,
\begin{align*}
\sum_{\substack{c | (m' / a) \\ (c, D_B) = 1}} \ &  \sum_{\substack{ \xi \in \Omega^o(|\Delta|m^2) \\ c(\xi) = c \\ \text{mod } \Gamma'}}  [Z^o(\xi[p^{\infty}])]  \\
&=  \sum_{\substack{c | (m' / a) \\ (c, D_B) = 1}} \ \sum_{\substack{ \xi \in \Omega^o( |\Delta|m'^2/ a^2) \\ c(\xi) = c \\ \text{mod } \Gamma'}} \left[ Z^o\left( a|\Delta| \cdot \xi[p^{\infty}] \right) \right]. \\
&= \sum_{\substack{c | (m' / a) \\ (c, D_B) = 1}} \ \sum_{\substack{ \xi \in \Omega^o( |\Delta|m'^2/ a^2) \\ c(\xi) = c \\ \text{mod } \Gamma'}} \left[ Z^o\left(  \xi[p^{\infty}] \right) \right] \qquad (\text{since } a|\Delta| \in \Z_p^{\times}) \\
&=  \sum_{\substack{\xi \in \Omega^o( |\Delta|m'^2/ a^2) \\ \text{mod } \Gamma'}} \left[ Z^o\left(  \xi[p^{\infty}] \right) \right] \qquad (\text{by } \eqref{OmegaoCondDecompEqn})\\ 
&= \widetilde{\Zed}{}^o \left( |\Delta| \frac{(m')^2}{a^2} \right).
\end{align*}
Substituting this back into \eqref{AFinalExp} concludes the proof of the theorem in the case $ord_p(m) > 0$. 

\fbox{\textbf{Case 2:} $ord_p(m) = 0$.}  As the proof is along similar lines as the previous case, we shall only indicate the necessary modifications. Fix an element $t \in \OB$ with $Nrd(t) = p \cdot gcd(m, D_B)$,  and such that the image of $t$ in $\OBp$ is the uniformizer $\Pi$ as in Section 2. 
Let $\varpi \in End(\mathbf E)_{\Q}$ be as in Lemma \ref{varpiLemma}, so that in particular $\varpi^2 = - p q$ for a split prime $q$. Then if $\beta \in \Omega^+(m, \lie{a}, \phi)$, it follows that 
\begin{align*} 
\beta' := \iota_{\mathbf A}(t) \circ \beta  \circ \varpi^{-1} \in \Omega^- & \Big( \frac{m}{gcd(m,D_B)}, \  \lie{q}^{-1} \cdot \lie{a}', \ Ad_{t} \circ \phi \Big) \\
&= \Omega^+ \left( \frac{m}{gcd(m,D_B)}, \ \lie{q}^{-1} \cdot \lie{a}', \ Ad_{t} \circ \phi' \right), 
\end{align*}
where $\lie{q}$ is one of the prime ideals above $q$, and the embedding $\phi'$ is the conjugate of $\phi$.

Let $\xi = \xi(\beta) \in \Omega^o(|\Delta| m^2)$ denote the special endomorphism corresponding to $\beta$ as in Lemma \ref{eigenvectorLemma}. If $j = \xi[p^{\infty}]$ is the corresponding map on $p$-divisible groups, then $\beta[p^{\infty}]$ and $\beta'[p^{\infty}]$ are equal (up to scaling by $\Z_p^{\times}$) to the elements $\mathbf b_1$ and $\mathbf b_2$ described in Theorem \ref{locCycCompHorThm}; hence, we have 
\[ Z^o(\xi[p^{\infty}])^{pure}  \ = \ Z(\beta[p^{\infty}]) \ + \ Z(\beta'[p^{\infty}]). \]
Therefore, we have
\begin{align*}
	\begin{pmatrix} \text{left hand} \\ \text{side} \\ \text{of \eqref{FourierCoeffEqn}} \end{pmatrix} &= 
	 \sum_{[\phi]} \ \widetilde{\Zed}(m, \phi) \ + \ \widetilde{\Zed}(m/gcd(m,D_B), \phi)  \\
	&=2 \sum_{[\phi]}\  \widetilde{\Zed}{}^+(m, \phi) \ + \ \widetilde{\Zed}{}^+(m/gcd(m,D_B), \phi) \quad [ \text{by } \eqref{signSwitchEqn}]  \\
	&=2 \sum_{[\phi]} \sum_{[\lie{a}]} \frac{1}{|o_k^{\times}|} \Bigg\{ \sum_{\substack{\beta \in \Omega^+(m, \lie{a}, \phi) \\ \text{mod } \Gamma'}} 
		[ Z(\beta[p^{\infty}])]  \\
	 & \qquad \qquad \qquad + \sum_{\substack{\beta' \in \Omega^+(m/(m,D_B), \lie{a}, \phi) \\ \text{mod } \Gamma'}} [Z(\beta'[p^{\infty}])] \Bigg\}  \\
	&=2 \sum_{[\phi]} \sum_{[\lie{a}]} \frac{1}{|o_k^{\times}|} \Bigg\{\sum_{\substack{\beta \in \Omega^+(m, \lie{a}, \phi) \\ \text{mod } \Gamma'}} [ Z(\beta		[p^{\infty}])] \\
	&\qquad \qquad  \qquad \qquad + \sum_{\substack{\beta' \in \Omega^+(m/(m,D_B), \lie{q}^{-1} \lie{a'}, Ad_t \circ \phi') \\ \text{mod } \Gamma'}} [Z(\beta'[p^{\infty}])] \Bigg\}   \\
	&= 2 \sum_{[\phi]} \sum_{[\lie{a}]} \frac{1}{|o_k^{\times}|} \Bigg\{ \sum_{\substack{\beta \in \Omega^+(m, \lie{a}, \phi) \\ \text{mod } \Gamma'}} [Z^o(\xi[p^{\infty}])] \Bigg\},
\end{align*}
where in the last line $\xi = \xi(\beta)$. 
The proof proceeds from this point exactly as in the previous case, cf.\  \eqref{mainThmIntExp1}. 
\end{proof}
\end{theorem}
\end{subsection}

\begin{subsection}{The main theorem and applications}
Let $\omega$ denote the relative dualizing sheaf of $\CB$, which we view as a divisor class in the Chow group $CH^1(\CB)$, and choose any divisor $K$ in this class.
Define the \emph{orthogonal generating series}
\begin{equation} \label{orthGenSeriesDef}
 \Phi^o(\tau) \ :=  - K +  \sum_{n > 0} \Zed^o(n) \  q_{\tau}^n \ \in Div(\CB) \llbracket q_{\tau} \rrbracket ,
\end{equation}
and consider the base change
\[ \Phi_{/o_k}^o(\tau) \ :=  - K_{/o_k} +  \sum_{n > 0} \Zed^o(n)_{/o_k} \  q_{\tau}^n \ \in Div(\CBok) \llbracket q_{\tau} \rrbracket. \]
Let $\chi_k = (\cdot / \Delta)$ denote the quadratic character attached to $k$, and $\chi_k'$ denote its induction to level $4D_B |\Delta|$. Define the Gauss sum
\[ \check{\chi}_k'(a) \ := \ \sum_{h = 0}^{4 D_B |\Delta| - 1} \chi_k'(h) \ \exp (2 \pi i a h / 4 D_B |\Delta|), \]
and the `L-function'
\[ L(s, \check{\chi}_k' ) \ := \sum_{m>0} m^{-s} \ \check{\chi}'_k(m), \]
which is analytic (as written) in the half-plane $\Re(s)>0$. 

We then define the unitary generating series:
\begin{align} \label{unGenSeriesDef} \Phi^u(\tau) :=  \frac{i}{2\pi} L(1, \check{\chi}_k') \cdot K_{/o_k}   +  \frac{1}{2h(k)} \sum_{m>0}& \sum_{[\phi]} \left[ \Zed(m, \phi) + \Zed^*(m , \phi) \right] q_{\tau}^m \\
& \in Div(\CBok) \llbracket q_{\tau} \rrbracket \otimes_{\Z} \Q \notag
 \end{align}
where, as we recall, 
\[\Zed^*(m, \phi) = \Zed( m /  \gcd(m, D_B), \phi) ,\]
 and the sum on $[\phi]$ is over the set $Opt / \OBx$ of optimal embeddings taken up to $\OBx$-conjugacy. Note that
\[ \frac{i}{2 \pi} L( 1, \check{\chi}_k')\ =  \ - \frac{h(k)}{|o_k^{\times}|} \prod_{\ell | D_B} \left( \frac{2 \ell^2 + \ell - 1}{\ell^2} \right) \ \in \  \Q, \]
which can be seen by evaluating the $L$-function on the left via Dirichlet's class number formula. 

\begin{theorem}[Main theorem] \label{mainThm}
Suppose $\Delta<0$ is squarefree and even, and every prime dividing $D_B$ is inert in $k$. Then we have an equality 
\begin{equation} \label{mainThmEqn}
  Sh (\Phi^o_{/o_k})(\tau) \ = \ \Phi^u(\tau)  \ \in Div(\CBok) \llbracket q_{\tau} \rrbracket \otimes \Q 
\end{equation}
of formal generating series, where $Sh$ is the formal Shimura lift with parameters $\kappa = 3, N = D_B, t = |\Delta|$ and $\chi = 1$ in the notation of Definition \ref{FShimLiftDef}.
\begin{proof}
We apply the formulas \eqref{shLiftFourierEqn} and \eqref{holShLiftCoeffZero} for the formal Shimura lift. Note that the constant terms match by design, and so -- keeping in mind the shift by $t=|\Delta|$ in the exponent in \eqref{ShLiftExpEqn} --   it suffices to prove that
\begin{compactenum}[(i)]
\item if $|\Delta|$ does not divide $m$, then $\sum_{[\phi]} \Zed(m, \phi) = 0$;
\item and  if $m = m'|\Delta|$, then
\begin{align*} 
  \sum_{[\phi] \in Opt/ \OBx}  & {\Zed}(m, \phi) + {\Zed}   \left(\frac{m}{ gcd(D_B,m)}, \phi \right) \notag \\
 &= 2 h(k)\sum_{\substack{\alpha | m' \\ (\alpha, D_B)=1 }}  \  \chi_k(\alpha)  \ {\Zed}{}^o \left( |\Delta| \frac{(m')^2}{\alpha^2} \right)_{/o_k}, 
\end{align*}
\end{compactenum}
in $Div(\CBok)$.

By \cite[Proposition 3.4.5]{KRYbook}, each orthogonal special cycle decomposes as 
\[ \Zed^o(n)_{/o_k} = \Zed^o(n)^{hor} \ + \ \sum_{q | D_B} \Zed^o(n)^{ver}_q, \]
where $\Zed^o(n)^{hor}$ is the closure of the generic fibre in $\Zed^o(n)_{/o_k}$, and $\Zed^o(n)^{ver}_q$ is the sum of the vertical irreducible components supported in characteristic $q$. By Proposition \ref{UCycleIsDivisorProp}, we have  an identical decomposition for the unitary  cycles. 

Now for each $q|D_B$, Theorem \ref{mainCompletionsThm} implies that (i) and (ii) hold for the vertical components at $q$.  Moreover, since $\CBok $ is proper over $\Spec(o_k)$, the same proposition implies (i) and (ii) hold in the generic fibre, by Grothendieck's existence theorem \cite[Thm.\ 5.1.4]{EGA3}; more precisely, if we choose a prime $p|D_B$, then  the existence theorem asserts that the desired relations hold over $o_{k,p}$, and then we observe that $o_{k,p}$ is faithfully flat over the localization $(o_k)_{(p)}$.

\end{proof}
\end{theorem}

Our next step is to recall how certain numerical quantities involving the orthogonal special cycles $\Zed^o(n)$ arise as Fourier coefficients of actual modular forms of weight 3/2, as in \cite[\S 4]{KRYbook}. 

First, we consider the `rational degree' generating series
\[ \Phi^o_{deg}(\tau) \ := \ - \deg_k(\omega) \ + \ \sum_{n>0} \deg_k \Zed^o(n) \ q_{\tau}^n \ \in \ \Q\llbracket q_{\tau} \rrbracket; \]
obtained by taking the degrees of the generic fibres. 
Next, we may consider the generating series 
\[ \Phi^o_{/k} (\tau) \ = \ - \omega_{/k}  \ + \ \sum_{n > 0} \Zed^o(n)_{/k} \ q_{\tau}^n \ \in \ CH_{\C}^1(\CB_{/k}) \llbracket q_{\tau} \rrbracket  
\] 
formed by taking generic fibres and then passing to the  Chow group $CH_{\C}^1(\CB_{/k}) = CH^1(\CB_{/k}) \otimes_{\Z} \C$. 

Finally, let  $\mathscr Y$ denote an irreducible component of the fibre $(\CBok)_p$ of the Shimura curve, where $p | D_B$ is a prime of bad reduction. For
a closed substack $\Zed$ of $\CBok$, we define the pairing
\[ \la \Zed,  \mathscr Y \ra := 2 \log(p) \cdot \chi( \mathcal{O}_{\mathscr{Y}} \otimes^{\mathbb L} \mathcal{O}_{\Zed} ), \]
where $\chi$ is the stack version of the Euler-Poincar\'e characteristic, which takes into account the automorphism groups of points, cf. \cite[\S VI.4]{DeRa}. 
We then form the generating series
\begin{equation} \label{PhiODefEqn}
	 \Phi^o_{\mathscr Y}(\tau) \  := \ - \la \omega, \mathscr{Y}\ra +  \sum_{n > 0} \ \la {\Zed}{}^o(n), \ \mathscr{Y} \ra \ q_{\tau}^{n}
\end{equation}

We may form the analogues  $\Phi^u_{deg}(\tau)$, $\Phi^u_{/k}(\tau)$ and $\Phi^u_{\mathscr Y}(\tau)$ for the unitary cycles, by applying the appropriate operators to the coefficients.

\begin{corollary}
(i) Each of $\Phi^o_{\deg}(\tau)$, $\Phi^o_{/k}(\tau)$ and $\Phi^o_{\mathscr Y}(\tau)$ is (the $q$-expansion of ) a holomorphic modular form of weight 3/2, level $4D_B$ and trivial character. 

(ii) Under the hypotheses of Theorem \ref{mainThm}, the Shimura lift of each of these modular forms is equal to its unitary counterpart. In particular, the series $\Phi^u_{\deg}$, $\Phi^u_{/k}$ and $\Phi^u_{\mathscr Y}$ are the $q$-expansions of modular forms of weight 2. 
\begin{proof}
(i) The modularity of these three generating series is proved in \cite[\S 4]{KRYbook}; specifically, see Equation (4.2.12),  Theorem 4.5.1 and Theorem  4.3.4.

(ii) This follows immediately by applying the appropriate map to both sides of \eqref{mainThmEqn}. 
\end{proof}
\end{corollary}

\begin{remark} Consider the more natural generating series 
\[ \Phi_{\phi}^{naive} \ = \ C + \frac{1}{2h(k)}\sum_{n>0} \Zed(n, \phi) q^n \ \in \ Div(\CBok) \llbracket q \rrbracket \otimes \Q \]
which omits the $\Zed^*$ terms appearing in $\Phi^u$;
on account of the general philosophy of such generating series, one may be tempted to conjecture that  for an appropriate constant term $C$, the $\Phi_{\phi}^{naive}$ is already  `modular' (to fix ideas, we may take this to mean that applying a linear functional that factors through the Chow group yields the $q$-expansion of a modular form). If this is indeed the case, we may explain the modularity of $\Phi^u$ in the following way. For an integer $d$, let $U_d$ and $B_d$ denote the following operators on formal $q$-expansions: if $F = \sum_{n \geq 0} a(n) q^n$, then
\[ U_d(F) = \sum a(dn) q^n, \qquad B_d(F) = \sum a(n) q^{dn} .\]
When $F$ is the $q$-expansion of a modular form, these operators are in fact induced by the maps $U_d = U_d^1$ and $B_d $ in \cite{Li}; both are maps between spaces of modular forms that preserve the weight, but may change the level. 
Set
\[ \varphi_d := B_d \circ \left( 1 \ - \ U_d \right), \qquad \varphi_d(F) \ = \ \sum \left( a(n) - a(dn) \right) \ q^{dn}\]
and note that if $(d, d') = 1$, then $\varphi_d$ and $\varphi_{d'}$ commute.
Next, let $\mathcal{P}(D_B)$ denote the set of primes dividing $D_B$, and for a non-empty subset $I = \left\{ p_1, \dots, p_n \right\} \subset \mathcal{P}(D_B)$,  we define
\[ \varphi_I \ := \ \varphi_{p_1}  \circ \varphi_{p_2} \circ \cdots \circ \varphi_{p_n}. \]
Then a direct calculation reveals that at the level of formal generating series,
\[ \Phi^u = C' \ + \sum_{[\phi]} \left( 2 + \sum_{I \subset \mathcal{P}(D_B)  } \varphi_I\right) (\Phi_{\phi}^{naive}), \]
where 
\begin{align*} C' \ &= \ \frac{i}{2\pi}  L(1, \check{\chi}_k') \cdot K\ - \  2 \cdot \# | Opt / \OBx | \cdot C. \\
\end{align*}
In particular, if $C$ happens to take on the  value
\[ C \  =  \ \frac{i}{4\pi} \cdot \frac{L(1, \check{\chi}_k')}{\# | Opt / \OBx|} \cdot K \ = \ - \frac{2^{-o(D_B)}} {2| o_k^{\times}|}  \prod_{\ell | D_B} \left( \frac{2 \ell^2 + \ell - 1}{\ell^2}  \right) \cdot K,  \]
then $C' = 0$; this would imply that whenever the application of a suitable functional to the coefficients of $\Phi_{\phi}^{naive}$ yields a modular form, the same is true for $\Phi^u$. It would be interesting to know if this implication can be reversed, and also if there is a geometric interpretation to the value of $C$ given above.  $\diamond$

\end{remark}

\end{subsection}


\end{section}

\Addresses

\small
\begin{thebibliography}{1}
\bibitem[BC]{BC} J.-F. Boutot and H. Carayol, \emph{Uniformisation p-adique des courbes de Shimura: les th\'eor\`emes de \v{C}erednik et de Drinfeld}. Courbes modulaire et courbes de Shimura, Ast\'erisque \textbf{196--197} (1991), pp. 45--158. 
\bibitem[Cip]{Cip} B. Cipra, \emph{On the Niwa-Shintani theta-kernel lifting of modular forms}. Nagoya Math. J., \textbf{91} (1983), pp. 49-117. 
\bibitem[DeRa]{DeRa} P. Deligne and M. Rapoport, \emph{Les sch\'emas de modules de courbes elliptiques}. Modular Functions of One Variable II (Proc. Intl. Summer School, Univ. Antwerp 1972). Lecture Notes in Mathematics, \textbf{349}. Springer-Verlag, (1973).
\bibitem[EGA3]{EGA3} A. Grothendieck, \emph{\'El\'ements  de g\'eometrie alg\'ebriques III}. Pub.\ Math.\ IH\'ES \textbf{11} (1961), pp. 5 -- 167.
\bibitem[Gro]{Gross}  B. Gross, \emph{On canonical and quasi-canonical liftings}. Invent. math., \textbf{84} (1986), pp. 321--326. 
\bibitem[How]{How} B. Howard, \emph{Intersection theory on Shimura surfaces}. Compos. Math. 145, \textbf{2} (2009), 423--475. 
\bibitem[How2]{How2} B. Howard, \emph{Complex multiplication cycles and Kudla-Rapoport divisors}, Annals of Math., \textbf{176} (2012), pp. 1097--1171 
\bibitem[KR1]{KRpadic} S. Kudla and M. Rapoport, \emph{Height pairings on Shimura curves and $p$-adic uniformization}. Invent. math., \textbf{142} (2000), pp. 153--222. 
\bibitem[KR2]{KRloc} S. Kudla and M. Rapoport, \emph{Special cycles on unitary Shimura varieties I: unramified local theory}, Invent. math., \textbf{184} (2010), pp. 629--682.
\bibitem[KR3]{KRunn} S. Kudla and M. Rapoport, \emph{Special cycles on unitary Shimura varieties II: global theory}, to appear in J. Reine Angew. Math.
\bibitem[KR4]{KRdrin} S. Kudla and M. Rapoport, \emph{An alternative description of the Drinfeld upper half-plane}, preprint. (arXiv:1108.5713)
\bibitem[KRY]{KRYbook} S. Kudla, M. Rapoport, and T. Yang, Modular Forms and Special Cycles on Shimura Curves. Annals of Mathematics Studies, \textbf{161}. Princeton University Press. 2006
\bibitem[Li]{Li} W. Li, \emph{Newforms and functional equations}. Math. Ann. \textbf{212} (1975), pp. 285--315.
\bibitem[Mann]{Mann} W.R. Mann, appendix to: ``Gross-Zagier revisited" by B. Conrad. Appearing in the volume: \emph{Heegner points and Rankin L-series}, Math. Sci. Res. Inst. Publ., \textbf{49}. Cambridge Univ. Press, Cambridge, 2004. 
\bibitem[Mes]{Messing} W. Messing, The Crystals Associated to Barsotti-Tate Groups: with Applications to Abelian Schemes. Lecture Notes in Mathematics, \textbf{264}. Springer-Verlag, Berlin. 1972.
\bibitem[Mum]{Mum} D. Mumford. Abelian Varieties. Tata Institute of Fundamental Research Studies in Mathematics \textbf{5}, Oxford University Press. 1970.
\bibitem[MZ]{MierendorffZiegler} E. Mierendorff and K. Ziegler, \emph{Formal moduli of formal $\mathcal O_K$ modules}. Appearing in: ARGOS Seminar on intersections of modular correspondences. Ast\'erisque \textbf{312} (2007). 
\bibitem[RZ]{RZ} M. Rapoport and Th. Zink. Period Spaces for $p$-divisible Groups. Annals of Mathematics Studies, \textbf{141}. Princeton University Press. 1996. 
\bibitem[San]{San} S. Sankaran, \emph{Special cycles on Shimura curves and the Shimura lift}. Ph.\ D.\ thesis, University of Toronto (2012).
\bibitem[San2]{San2} S. Sankaran, \emph{Unitary cycles on Shimura curves and the Shimura lift II}, preprint. 
Available at: \url{http://www.math.uni-bonn.de/~sankaran/shlift2.pdf}
\bibitem[Shim]{Shim}G. Shimura, \emph{On modular forms of half integral weight}. Ann. of Math., \textbf{97} (1973), pp. 440--481.
\bibitem[Ter]{terstiegeAnti} U. Terstiege, \emph{Antispecial cycles on the Drinfeld upper half plane and degenerate Hirzebruch-Zagier cycles}. Manuscripta Mathematica, \textbf{125} (2008), pp. 191--223. 
\bibitem[Ter2]{Ter2} U. Terstiege, \emph{Intersections of special cycles on the Shimura variety for GU(1,2)}. to appear in J. Reine. Angew. Math. 
\bibitem[Vig]{Vig} M.-F. Vigneras, Arithmetique des Algebres de Quaternions. Lecture Notes in Mathematics, \textbf{800}. Springer-Verlag, Berlin. 1980
\bibitem[Wew]{Wew} S. Wewers, \emph{Canonical and quasi-canonical liftings}. Appearing in: ARGOS Seminar on intersections of modular correspondences. Ast\'erisque \textbf{312} (2007). 


\end{thebibliography}
\end{document}